\documentclass[10pt]{article}

\usepackage{amssymb,amsfonts,amsmath,amscd} 
\usepackage{amsthm}
\usepackage{stmaryrd}
\usepackage{graphicx,graphics}
\usepackage{latexsym}
\usepackage{color} 
\usepackage[all]{xy}
\usepackage{pstricks}
\usepackage{pst-node}
\usepackage{setspace} 
\usepackage{verbatim}
\usepackage{hyperref} 
\usepackage{fullpage}
\newtheorem{thm}{Theorem}[section]
\newtheorem{prop}[thm]{Proposition}
\newtheorem{lem}[thm]{Lemma}

\newtheorem{col}[thm]{Corollary}
\newtheorem{prop-defn}[thm]{Proposition-Definition}

\theoremstyle{definition}
\newtheorem{defn}[thm]{Definition}
\newtheorem{exa}[thm]{Example}
\newtheorem{rem}[thm]{Remark}

\numberwithin{equation}{subsection}

\newcommand{\Span}[1]{\left<#1\right>}
\newcommand{\Hilb}[2]{#1\text{-Hilb}(#2)}

\newcommand{\C}{\mathbb C}
\newcommand{\R}{\mathbb R}
\newcommand{\N}{\mathbb N}
\newcommand{\Z}{\mathbb Z}

\newcommand{\OO}{\mathcal O}

\DeclareMathOperator{\GL}{GL}
\DeclareMathOperator{\SL}{SL}

\DeclareMathOperator{\Irr}{Irr}
\DeclareMathOperator{\HILB}{Hilb}

\DeclareMathOperator{\BD}{BD}

\DeclareMathOperator{\Spec}{Spec}

\DeclareMathOperator{\GCD}{gcd}
\DeclareMathOperator{\Aut}{Aut}

\DeclareMathOperator{\Sym}{Sym}

\begin{document}

\title{$G$-graphs and special representations \\for binary dihedral groups in $\GL(2,\C)$} 
\author{\'{A}lvaro Nolla de Celis} 
\date{}
\maketitle

\begin{abstract} Given a finite subgroup $G\subset\GL(2,\C)$ it is know that the minimal resolution of the singularity $\C^2/G$ is the moduli space $Y=\Hilb{G}{\C^2}$ of $G$-clusters $\mathcal{Z}\subset\C^2$. The explicit description of $Y$ can be obtained by calculating every possible distinguished basis for $\mathcal{O}_\mathcal{Z}$ as vector spaces. These basis are the so called $G$-graphs. In this paper we classify the $G$-graphs for any small binary dihedral subgroup $G$ in $\GL(2,\C)$ and, in the context of the special McKay correspondence, we use this classification to give a combinatorial description of the special representations of $G$ appearing in $Y$ in terms of its maximal normal cyclic subgroup $H\unlhd G$.  \end{abstract}



\section{Introduction}
Motivated by the McKay correspondence \cite{McK80} (read \cite{Rei02} for a survey) $G$-graphs were introduced by Nakamura \cite{Nak01} to construct a crepant resolution of the quotient $\C^3/G$ for abelian subgroups $G\subset\SL(3,\C)$. This resolution is the moduli space $G$-Hilb of $G$-clusters introduced by Ito and Nakamura \cite{IN96}. Among other properties, it has been proved to be the minimal resolution of $\C^2/G$ for small finite subgroups $G\subset\GL(2,\C)$ in \cite{Ish02}, and a crepant resolution for subgroups $G\subset\SL(3,\C)$ in \cite{BKR01}. 

Since then, $G$-graphs have been a useful tool in studying $G$-Hilb when $G$ is abelian (for instance \cite{Kidoh}, \cite{Ito02} for $G\subset\GL(2,\C)$ and \cite{CR02} for $G\subset\GL(3,\C)$). The first attempt to extend the notion of $G$-graph to calculate $G$-Hilb for non-abelian subgroups was made by Leng in \cite{Leng} for binary dihedral subgroups $G\subset\SL(2,\C)$ and some binary trihedral subgroups $G\subset\SL(3,\C)$ (see also \cite{Seb05}). In this paper we calculate every possible $G$-graph for small binary dihedral subgroups $G\subset\GL(2,\C)$. 

The key idea in the construction of these $G$-graphs is to consider the action of a dihedral group $G$ on $\C^2$ as the cyclic action by its maximal normal abelian subgroup $H$ of index two, followed by a dihedral involution. This interpretation allows us to construct $G$-graphs from the union of two $H$-graphs identified by the action of $G/H$, that we call $qG$-{\em graphs}. From a $qG$-graph, we use the representation theory of $G$ to construct a $G$-graph in a unique way. In Propositions \ref{TypeA}, \ref{TypeB1} and \ref{TypeB2} we classify the $G$-graphs arising from $qG$-graphs into types $A$ and $B$. In propositions \ref{TypeD} and \ref{TypeC} we classify the $G$-graphs which do not come from $qG$-graphs into types $C$ and $D$.

Given any $G$-graph $\Gamma$ there exists a corresponding affine open set $U_{\Gamma}$ in $\Hilb{G}{\C^2}$ which consists of all $G$-clusters $\mathcal{Z}$ for which $\Gamma$ is a basis of the vector space $\mathcal{O}_{\mathcal{Z}}$. In Theorem \ref{ABCD} and the following Corollary \ref{cover} we prove that the number of distinct $G$-graphs gives us a covering of $\Hilb{G}{\C^2}$ with the minimum number of open sets. The explicit equations of these open sets are given in \cite{NdC2} using the moduli space of representations of the McKay quiver. 

In the last part of the paper we apply the classification of $G$-graphs to list the special representations of any small binary dihedral group $G\subset\GL(2,\C)$. We start by using the explicit description of the ideals corresponding to the $G$-graphs that we give in Section \ref{SecGraphs}, to produce a 1-parameter family of ideals connecting any $G$-cluster at the exceptional divisor $E$ in $\Hilb{G}{\C^2}$. Thanks to the {\em special} McKay correspondence \cite{Rie04} the exceptional curves appearing in the minimal resolution correspond one-to-one to the special irreducible representations, which now can be expressed in terms of $G$-graphs. This minimal resolution is $\Hilb{G}{\C^2}$ by Ishii's Theorem (see \cite{Ish02} and Theorem \ref{thm!Ishii}), which lead us to give in Theorem \ref{thm:Special} a combinatorial description of the special representations of a small binary dihedral group $G=\BD_{2n}(a)$ in terms of the continued fraction $\frac{2n}{a}$. We note that the classification of the special representations was discovered independently by Wemyss and Iyama in \cite{IW08} using Cohen-Macaulay modules. 

The paper is distributed as follows: In Section \ref{Cyc-quot} we recall the background material on cyclic quotient singularities that is needed. In Section \ref{DihedralGroups} we describe the binary dihedral groups in $\GL(2,\C)$ and the resolution of singularities of $\C^2/G$. In Section \ref{GHilbGgraphs} we define $G$-graphs and recall the explicit construction of $\Hilb{H}{\C^2}$ in terms of $H$-graphs when $H$ is an abelian subgroup in $\GL(2,\C)$. Section \ref{SecGraphs} is dedicated to the calculation of the $G$-graphs for the dihedral groups and their classification into the types $A$, $B$, $C$ and $D$. In Section \ref{Sect:Walk} we prove that any $G$-cluster $\mathcal{Z}$ admits as basis for $\mathcal{O}_{\mathcal{Z}}$ a $G$-graph of type $A$, $B$, $C$ or $D$, thus obtaining an open cover for $\Hilb{G}{\C^2}$. Finally Section \ref{SpecialRepr} is devoted to the special representations and their description in terms of the continued fraction $\frac{2n}{a}$. 

I would like to thank M. Reid for introducing me to this topic and for his support during my PhD.

\section{Background on cyclic quotient singularities}\label{Cyc-quot}

In this section we introduce some notation needed for the rest of the paper about cyclic quotient singularities. The material is taken from \cite{Reid}. 

Let $H=\Span{\frac{1}{k}(1,a)}:=\Span{\left(\begin{smallmatrix}\varepsilon&0\\0&\varepsilon^{a}\end{smallmatrix}\right)|\varepsilon^k=1\text{ primitive}}$ be a cyclic group in $\GL(2,\C)$ with $(k,a)=1$. The quotient singularity $X:=\C^2/H=\Spec\C[x,y]^H$ and the minimal resolution $Y\to X$ are toric and are completely determined by the continued fractions $\frac{k}{k-a}$ and $\frac{k}{a}$ as follows. Let $a_i$ for $i=1,\ldots,l$ be the entries of the Hirzebruch--Jung continued fraction 
\[\frac{k}{k-a}=a_1-\frac{1}{a_2-\ldots\frac{1}{a_{l}}}=[a_{1},\ldots,a_{l}]\]

The ring $\C[x,y]^H$ of invariants is generated by the monomials $u_i$ for $i=0,\ldots,l$ that satisfy 
\begin{align}
u_{i-1}u_{i+1}&=u_i^{a_i} \text{ for }i=1,\ldots,l \label{Cyc-eqns}
\end{align}
where $u_0=x^k$ and $u_1=x^{k-a}y$. Then $X\subset\C^{l+2}$ is determined set-theoretically by the relations (\ref{Cyc-eqns}).

Now consider the lattice $L:=\Z^2+\frac{1}{k}(1,a)\cdot\Z\subset\R^2$ and its dual lattice of invariant monomials $M$. Define the Newton polygon of $L$ as the convex hull in $\R^2$ of all nonzero lattice points in the positive quadrant. Then the resolution $Y$ is determined explicitly by the continued fraction 
\[\frac{k}{a}=[b_{1},\ldots,b_{m}]\]
in the following way; let $e_0=\frac{1}{k}(0,k)=(0,1)$, $e_1=\frac{1}{k}(1,a)$ and $e_{i+1}+e_{i-1}=b_ie_i$ for $i=1,\ldots,m$, be the nonzero lattice points on the boundary of the Newton polygon of $L$. Then the exceptional divisor $E\subset Y$ consists of $m$ exceptional curves $E_{1},\ldots,E_{m}$ where each $E_{i}\cong\mathbb{P}^1$ with selfintersections $-b_{1},\ldots,-b_{m}$ respectively. These rational curves intersect according to the following dual graph of type $A$:

\begin{center}
\begin{pspicture}(0,0)(2.5,0.5)
\scalebox{0.8}{
	\psset{nodesep=3pt}
	\rput(0,0){\rnode{p1}{$\bullet$}}\rput(0,0.35){$E_1$}
	\rput(1,0){\rnode{p2}{$\bullet$}}\rput(1,0.35){$E_2$}
	\rput(2,0){\rnode{p3}{$\cdots$}}
	\rput(3,0){\rnode{p4}{$\bullet$}}\rput(3.1,0.35){$E_m$}
	\ncline{-}{p1}{p2}\ncline{-}{p2}{p3}\ncline{-}{p3}{p4}
	}
\end{pspicture}
\end{center}

Furthermore, $Y$ is covered by $m+1$ affine open sets $Y=Y_{0}\cup\ldots\cup Y_{m}$, where each $Y_i\cong\C^2$ has coordinates $\lambda_i,\mu_i$ defined to be the dual basis in $M$ of the consecutive lattice points $e_i,e_{i+1}$ of $L$. 

The relation between the entries of the continued fractions $\frac{k}{k-a}$ and $\frac{k}{a}$ is given by an algorithm due to Riemenschneider (see \cite{Rie74}): given the entries $a_1,\ldots,a_{l}$ of the continued fraction $\frac{k}{k-a}$, form $l$ rows with $a_i-1$ points in each row as follows
\begin{center}
\begin{pspicture}(0,0)(3.5,1.75)
\rput(0,1.25){$\underbrace{\times\times\ldots\times\times}_{\tiny a_1-1}$}
\rput(1.55,0.75){$\underbrace{\times\times\ldots\times\times}_{\tiny a_2-1}$}
\rput(3.1,0.25){$\underbrace{\times\times\ldots\times\times}_{\tiny a_3-1}$}
\rput(4.25,0){\rput(0,0){$\cdots$}}
\end{pspicture}
\end{center}
where for $i=1,\ldots,l$ the first point in row $i$ is placed in the same column as the last point in row $i-1$. Then for $j=1,\ldots,m$ the number of points in column $j$ is $b_j-1$. Vice versa, given $b_1,\ldots,b_m$ we can recover $a_1,\ldots,a_l$ in the same way.

\begin{exa} Let $k=46$ $a=17$ with $\frac{46}{17}=[3,4,2,3]$. Then we write 
\begin{pspicture}(0,0)(1.5,0.5)
\rput(0.4,-0.325){
\rput(0,0.75){$\times\times$}
\rput(0.43,0.5){$\times\!\!\times\!\!\times\ $}
\rput(0.6125,0.25){$\times$}
\rput(0.75,0){$\times\times$}
}
\end{pspicture}
so that $\frac{46}{29}=[2,3,2,4,2]$. 
\end{exa}

\section{Dihedral groups in $\GL(2,\C)$} \label{DihedralGroups}

We consider the following representation of binary dihedral subgroups in $\GL(2,\C)$ in terms of their action on the complex plane $\C^2_{x,y}$:
\[
\BD_{2n}(a)=\left< \alpha=
\begin{pmatrix}
\varepsilon & 0 \\
0 & \varepsilon^{a} \\
\end{pmatrix},~
\beta= 
\begin{pmatrix}
0 & 1 \\
-1 & 0 \\
\end{pmatrix}
 :~ \varepsilon^{2n}=1\text{ primitive, $a^2\equiv 1$ (mod $2n$)} \right>
\]

\noindent In other words, $\BD_{2n}(a)$ is the group of order $4n$ generated by the cyclic group $H:=\Span{\alpha}=\Span{\frac{1}{2n}(1,a)}\unlhd G$ and the dihedral symmetry $\beta$ which sends the coordinates $(x,y)$ to $(-y,x)$. The subgroup $H$ is a choice of maximal cyclic index 2 subgroup of $G$ (note that $\beta^2\in H$). The condition $a^2\equiv 1$ (mod $2n$) is equivalent to the classical dihedral condition of $\alpha\beta=\beta\alpha^{a}$. 

We start by giving the definitions of the integers $q$ and $k$ which appear frequently throughout the paper.
\begin{defn}\label{defn:qk} Let $q:=\frac{2n}{(a-1,2n)}$, and $k$ such that $n=kq$. 
\end{defn}

An element $g\in G$ is a {\em quasireflection} if it fixes a hyperplane, and a group $G$ is called {\em small} if it does not contain any quasireflection. A theorem of Chevalley, Shephard and Todd \cite{Ch} states that if $H\subset\GL(n,\C)$ is generated by reflections then $\C^n/H\cong\C^n$, which traditionally reduces the study of these quotients to small groups. 

\begin{prop}\label{CritSmall} The group $G=\BD_{2n}(a)$ is small $\Longleftrightarrow\GCD(a+1,2n)\nmid n$.
\end{prop}
\begin{proof}
The elements of $G$ are of the form $\alpha^i$ and $\alpha^i\beta$ for $i=0,\ldots,2n-1$. Since $(2n,a)=1$ the subgroup $H$ is small, so quasireflections can only occur among the elements of the form $\alpha^i\beta$ (with $i\neq n$). Then  
\[
\alpha^{i}\beta \text{ is a quasireflection } \Longleftrightarrow \text{det}(\alpha^{i}\beta-I)=0  \Longleftrightarrow 1+\varepsilon^{(a+1)i} = 0 \Longleftrightarrow (a+1)i\equiv n \text{ (mod $2n$).} 
\]
Therefore, $G$ has no quasireflections if and only if do not exist any solutions to the equation $(a+1)x\equiv n \text{ (mod }2n)$. As a linear congruence, it only has solutions if the $\GCD(a+1,2n)$ divides $n$. 
\end{proof}

\begin{rem}\label{Brie}{\bf Brieskorn classification.} Small binary dihedral groups in $\GL(2,\C)$ were originally classified by Brieskorn in \cite{Bri} as follows:
\[
D_{N,q}:= \left\{ \begin{array}{ll} 
\Span{\psi_{2q}, \tau, \phi_{2k}}, & \text{if }k:=N-q\equiv1\text{ (mod 2)} \\ 
\Span{\psi_{2q}, \tau\circ\phi_{4k}}, & \text{if }k\equiv0\text{ (mod 2)} 
\end{array} \right.
\]
\[
\text{with }\psi_{r}=\begin{pmatrix}\varepsilon_{r}&0\\0&\varepsilon_{r}^{-1}\end{pmatrix}, \tau=\begin{pmatrix}0&i\\i&0\end{pmatrix}, \phi_{r}=\begin{pmatrix}\varepsilon_{r}&0\\0&\varepsilon_{r}\end{pmatrix}, \varepsilon_{r}=\exp\frac{2\pi i}{r} \text{ and $|D_{N,q}|=4kq$}
\]
where $q$ and $k$ are as in \ref{defn:qk}. The groups $\BD_{2n}(a)$ which are small correspond to the case $k$ odd. The case $k$ even is obtained by taking $\beta=\left(\begin{smallmatrix}0&1\\\varepsilon^q&0\end{smallmatrix}\right)$ (See \cite{thesis}, $\S3$ for more details). For simplicity, in this paper we only treat $\BD_{2n}(a)$ groups but we emphasize that the methods used here apply to the groups with $k$ even as well. 
\end{rem}

The group $\BD_{2n}(a)$ has irreducible 1-dimensional representations $\rho_{j}^+$ and $\rho_j^-$ of the form
\[
\begin{array}{cc}
\rho_{j}^\pm(\alpha)=\varepsilon^j, & \rho_{j}^\pm(\beta)={\small{\left\{\begin{array}{ll}\pm i&\text{ if $n,j$ odd}\\\pm1&\text{ otherwise}\end{array}\right.}}
\end{array}
\]
where $\varepsilon$ is a $2n$-th primitive root of unity and $j$ is such that $j\equiv aj$ (mod $2n$). The values $r$ for which $r\not\equiv ar$ (mod $2n$) form in pairs the irreducible 2-dimensional representations $V_{r}$ of the form
\[
\begin{array}{cc}
V_{r}(\alpha)=\begin{pmatrix}\varepsilon^r&0\\0&\varepsilon^{ar}\\\end{pmatrix}, & V_{r}(\beta)=\begin{pmatrix}0&1\\(-1)^r&0\end{pmatrix}
\end{array}
\]
By definition, the natural representation is $V_{1}$. 

Notice that the number of 1-dimensional representations coincides with twice the number of scalar diagonal elements in $\frac{1}{2n}(1,a)$. Since $n=kq$, then the number of 1-dimensional representations is $4k$. If we call $d$ the number of 2-dimensional irreducible representations, using the formula $|G| = \sum_{\rho\in\Irr G}(\dim(\rho))^2$ we have that $4n = 4k+4d$, which gives $d=n-k$. 

\begin{rem}\label{McKG}{\bf $\Irr G$ from $\Irr H$.} Irreducible representations of $H$ give rise to irreducible representations of $G$ in the following way. Let $\Irr H=\{\rho_0,\ldots,\rho_{2n-1}\}$. The group $G$ acts on $H$ by conjugation $g\cdot h:=ghg^{-1}$, for $g\in G, h\in H$. The induced action of $G$ on the characters is given by $g\cdot\chi_{\rho_j}(h):=\chi_{\rho_j}(g^{-1}hg)$. Since the character is a function on the conjugacy classes in $H$, this action is constant in the cosets $gH$. Therefore $G/H=\Span{\beta}\cong\Z/2\Z$ acts on the characters $\chi_{\rho_j}$ of $H$, which induces an action on $\Irr H$ by 
\[ 
\beta\cdot\rho_j := \rho_{aj}.
\]
The free orbits under the action of $G/H$ are $\{\rho_{j},\rho_{aj}\}$ with $aj\not\equiv j$ mod $2n$, and they combine to produce the 2-dimensional representation $V_{j}\in\Irr G$. Every fixed representation $\rho_{j}$ with $aj\equiv j$ (mod $2n$) splits into the two 1-dimensional representations $\rho_{j}^+$ and $\rho_{j}^-$ in $\Irr G$, corresponding to the two characters of $G/H\cong\Z/2\Z$.
\end{rem}

In what follows we take the notation as in \cite{Yos}, $\S10$. Let $V(=V_1)$ a vector space with basis $\{x,y\}$ where $G$ acts naturally. Define $S=\Sym V:=\C[V^*]$ the polynomial ring in the variables $x$ and $y$. Then the action of $G$ extends to $S$ by $g\cdot f(x,y):=f(g(x),g(y))$ for $f\in S$, $g\in G$.

\begin{defn}\label{def-action} Let $G=\BD_{2n}(a)$ and $f\in S$. A polynomial $f$ is said to belong to $\rho_j^\pm$ if
\[
\alpha(f)=\varepsilon^jf \text{ and } \beta(f)=\left\{\begin{array}{ll}\pm if&\text{ if $n,j$ odd}\\\pm f&\text{ otherwise}\end{array}\right.
\]
A pair of polynomials $(f,g)$ is said to belong to $V_k$ if $g=\beta(f)$ and $\alpha(f,\beta(f))=(\varepsilon^kf,\varepsilon^{ak}\beta(f))$.
\end{defn}

Let $S_{\rho}:=\{f\in\C[x,y]: f\in\rho\}$ the $S^G$-module of $\rho$-invariants. We say that a polynomial $f$ belongs to $\rho$, or simply $f\in\rho$ if it belongs to the corresponding module $S_\rho$. Note that these are precisely the Cohen Macaulay $S^G$-modules $S_\rho=(S\otimes\rho^*)^G$ where $G$ acts on $S$ as above and $G$ acts on a representation $\rho$ by the inverse transpose. 

\begin{exa} In Table \ref{SInvBD12} we present the irreducible representations for the group $\BD_{12}(7)$ together with some elements belonging to the corresponding modules $S_{\rho}$.
\end{exa}

{\renewcommand{\arraystretch}{1.1}
\begin{table}[htdp]
\begin{center}
\begin{small}
\begin{tabular}{|c|c|c|l|}
\hline
 & $\alpha$ & $\beta$ & ~~~~~~~~~~~~~~~~~$S_{\rho}$ \\ 
\hline
$\rho_0^+$ & 1 & 1 & $1$, $x^{12}+y^{12}$, $x^5y-xy^5$, $x^6y^6$ \\
\hline
$\rho_0^-$ & 1 & $-1$ &  $x^{12}-y^{12}$, $x^5y+xy^5$, $x^3y^3$ \\
\hline
$\rho_1^+$ & $\varepsilon^2$ & 1 &  $x^2+y^2$, $x^7y-xy^7$ \\
\hline
$\rho_1^-$ & $\varepsilon^2$ & $-1$ &  $x^2-y^2$, $x^7y+xy^7$ \\
\hline
$\rho_2^+$ & $\varepsilon^4$ & 1 &  $x^4+y^4$, $x^9y-xy^9$, $x^2y^2$ \\
\hline
$\rho_2^-$ & $\varepsilon^4$ & $-1$ & $x^4-y^4$, $x^9y+xy^9$, $x^5y^5$ \\
\hline
$\rho_3^+$ & $-1$ & 1 &  $x^6+y^6$, $x^{11}y-xy^{11}$ \\
\hline
$\rho_3^-$ & $-1$ & $-1$ &  $x^6-y^6$, $x^{11}y+xy^{11}$ \\
\hline
$\rho_4^+$ &  $\varepsilon^8$ & 1 & $x^8+y^8$, $x^6y^2+x^2y^6$, $x^4y^4$ \\
\hline
$\rho_4^-$ & $\varepsilon^8$ & $-1$ &  $x^8-y^8$, $x^6y^2-x^2y^6$, $xy$ \\
\hline
$\rho_5^+$ & $\varepsilon^{10}$ & 1 &  $x^{10}+y^{10}$, $x^3y-xy^3$ \\
\hline
$\rho_5^-$ & $\varepsilon^{10}$ & $-1$ &  $x^{10}-y^{10}$, $x^3y+xy^3$ \\
\hline
$V_1$ & $\left(\begin{smallmatrix}\varepsilon&0\\0&\varepsilon^7\end{smallmatrix}\right)$ & $\left(\begin{smallmatrix}0&1\\-1&0\end{smallmatrix}\right)$ & $(x,y)$, $(y^7,-x^7)$, \\
	& & & $(x^6y,-xy^6)$, $(x^2y^5,-x^5y^2)$ \\
\hline
$V_2$ & $\left(\begin{smallmatrix}\varepsilon^3&0\\0&\varepsilon^9\end{smallmatrix}\right)$ & $\left(\begin{smallmatrix}0&1\\-1&0\end{smallmatrix}\right)$ &  $(x^3,y^3)$, $(y^9,-x^9)$,  \\
	& & & $(xy^2,x^2y)$, $(x^8y,-xy^8)$ \\
\hline
$V_3$ & $\left(\begin{smallmatrix}\varepsilon^5&0\\0&\varepsilon^{11}\end{smallmatrix}\right)$ & $\left(\begin{smallmatrix}0&1\\-1&0\end{smallmatrix}\right)$ &  $(x^5,y^5)$, $(y^{11},-x^{11})$,  \\
	& & & $(xy^4,x^4y)$, $(x^{10}y,-xy^{10})$ \\
\hline
\end{tabular}
\end{small}
\end{center}
\caption{Some semi-invariant elements in each $S_{\rho}$ for $\BD_{12}(7)$}
\label{SInvBD12}
\end{table}}

\subsection{Resolution of dihedral singularities}\label{Resoln}

Let $G=\BD_{2n}(a)$ be a small binary dihedral group with cyclic maximal subgroup $H=\Span{\alpha}\unlhd G$. Let us now look at the geometric construction of the resolution of a dihedral singularity $\C^2/G$.  

Consider first the action of $H$ on $\C^2$. The quotient affine variety $X=\C^2/H$ is the toric quotient singularity of type $\frac{1}{2n}(1,a)$, with an isolated singular point at the origin. Recall from Section \ref{Cyc-quot} that the resolution of singularities $Y=\Hilb{H}{\C^2}\to X$ is determined by the continued fraction $\frac{2n}{a}=[b_1,\ldots,b_m]$.

\begin{lem}\label{lem:symm}
If $a^2\equiv1$  $(\mathrm{mod}\ 2n)$ then the entries of the continued fraction $\frac{2n}{a}$ are symmetric with respect to the middle term, that is $b_i=b_{m+1-i}$ for $i=1,\ldots,m$.
\end{lem}

\begin{proof}
Let $L$ be the lattice of weights and $M$ the dual lattice of monomials, and consider the continued fractions $\frac{2n}{2n-a}=[a_{1},\ldots,a_{l}]$ and $\frac{2n}{a}=[b_{1},\ldots,b_{m}]$. If a monomial $x^{i}y^j$ is $H$-invariant then $i+aj\equiv0$ (mod $2n$), and by the assumption $a^2\equiv1$ (mod $2n$), $x^jy^{i}$ is also invariant. Therefore, the continued fraction $\frac{2n}{2n-a}$ is symmetric, i.e.\ $a_{i}=a_{l+1-i}$. Indeed, let $u_{i-1}=x^{t'}y^{w'}$, $u_{i}=x^{t}y^{w}$, $u_{i+1}=x^{t''}y^{w''}$ be three consecutive invariant monomials for some integers $t$ and $w$. Let also $u_{l+1-(i+1)}=x^{w''}y^{t''}$, $u_{l+1-i}=x^{w}y^{t}$, $u_{l+1-(i-1)}=x^{w'}y^{t'}$ be their symmetric partners. Since $u_{i-1}u_{i+1}=u_{i}^{a_{i}}$ we have $t'+t''=a_{i}t=a_{l+1-i}t$ and $w'+w''=a_{i}w=a_{l+1-i}w$, thus $a_{i}=a_{l+1-i}$ for all $i$. 

The symmetry in the entries of $\frac{2n}{2n-a}$ implies the symmetry of the entries of $\frac{2n}{a}$ by the algorithm explained in Section \ref{Cyc-quot}, thus $b_{i}=b_{m+1-i}$. 
\end{proof}

To complete the action of $G$ on $\C^2$ we act on $Y=\Hilb{H}{\C^2}$ with $G/H=\Span{\bar{\beta}}\cong\Z/2\Z$. Notice that the symmetry in the continued fraction $\frac{2n}{a}$ given in Lemma \ref{lem:symm} implies that the coordinates along the exceptional curves $E_{i}$ in the resolution $Y\to\C^2/H$ are also symmetric, i.e.\ $\beta$ identifies the affine subsets $Y_{i}\cong\C^2_{(x^r/y^s,y^{u}/x^v)}$ with $Y_{m-i}\cong\C^2_{(x^u/y^v,y^r/x^s)}$ as well as the rational exceptional curves in $E$ covered by these affine patches. 

\begin{lem}\label{prop!even}
If $G$ is small then the continued fraction expansion of $\frac{2n}{a}$ has an odd number of entries.
\end{lem}

\begin{proof}
Suppose that the expansion of $\frac{2n}{a}$ has an even number of elements, and let $\{Y_{i}\}_{i=1,\ldots,2h-1}$ be the open affine covering of $\Hilb{H}{\C^2}$. Since $(2n,a)=1$ the abelian action has no quasireflections so we need to check only the action of $\beta$ on $\Hilb{H}{\C^2}$. This action identifies $Y_i$ with $Y_{m-i}$ thus the only $\beta$-fixed part can only occur at the middle open set $Y_{h}$ which is sent to himself. This open set is defined by the lattice points $e_{h-1}=\frac{1}{2n}(u,v)$ and $e_h=\frac{1}{2n}(v,u)$ for some $0<u<v$, so $Y_{h}\cong\C^2_{\lambda,\mu}$ where $\lambda=x^{u}/y^v$ and $\mu=y^{u}/x^v$. Then, the action of $\beta$ on $Y_h$ is of the form $(\lambda,\mu)\mapsto((-1)^v\mu,(-1)^{u}\lambda)$. 

Now notice that $u$ and $v$ have the same parity. Indeed, if we call $e_2=\frac{1}{2n}(c,d)$ we have that $c=b_1$ and $d=a-2n$. Thus, since $a$ is odd we have that $c$ and $d$ have the parity as $b_1$. Using the formula $e_{i+1}=b_ie_i-e_{i-1}$ we can use induction to conclude that the entries of every lattice point $e_i$ for the action $\frac{1}{2n}(1,a)$ with $(2n,a)=1$ have the same parity. In particular, since $u$ and $v$ have the same parity $\beta$ fixes the line $\lambda=(-1)^u\mu$, i.e.\ $\beta$ is a quasireflection so the group is not small. 
\end{proof}

Since we are only interested in the case when $G$ is small, from now on we suppose that $\frac{2n}{a}$ has an odd number of elements, i.e.\ $\frac{2n}{a} = \lbrack b_{1},\ldots,b_{h-1}, b_{h}, b_{h-1},\ldots,b_{1}\rbrack$. Then, the exceptional divisor $E\subset Y$ has an odd number of irreducible components $E_{i}$, and there exists a middle rational curve $E_{h}\cong\mathbb{P}^1$. This rational curve has coordinate ratio $(x^q:y^q)$ and it is covered by $Y_{h-1}$ and $Y_{h}$. Then $\beta$ identifies $Y_{h-1}$ with $Y_{h}$ and it is an involution on $E_{h}$ with two fixed points since there are no quasireflections. In the quotient $Y/\!\Span{\beta}$, these fixed points become two $A_{1}$ singularities, and blowing them up we obtain the Dynkin diagram of type $D$ that we were looking for. 

\begin{rem} For a small finite subgroup $G\subset\GL(2,\C)$, by a result of Ishii  (\cite{Ish02}) $\Hilb{G}{\C^2}$ is the minimal resolution of $\C^2/G$. By the uniqueness of minimal models in dimension 2, the previous construction can be expressed in the following diagram:

\begin{center}
\begin{pspicture}(0,0)(7,2.5)
	\psset{nodesep=5pt}
	\rput(-0.75,0){\rnode{p00}{$\Hilb{G}{\C^2}$}}
	\rput(1.5,0){\rnode{p20}{$Y/\!\Span{\beta}$}}
	\rput(-0.15,1.25){\rnode{p12}{$G/H\curvearrowright$}} 
	\rput(1.5,1.25){\rnode{p22}{$\Hilb{H}{\C^2}$}}
	\rput(3.75,1.25){\rnode{p42}{$\C^2/H$}}
	\rput(3.75,2.5){\rnode{p44}{$\C^2$}}
	\rput(5.8,2.5){\rnode{p54}{$\curvearrowleft\ H\unlhd G\subset\GL(2,\C)$}} 	
	\ncline{->}{p00}{p20}\ncline{->}{p22}{p20}
	\ncline{->}{p22}{p42}\ncline{->}{p44}{p42}
\end{pspicture}
\end{center}

\noindent In other words, we have that $\Hilb{G}{\C^2}\cong\Hilb{G/H}{\Hilb{H}{\C^2}}$.
\end{rem}

\begin{prop}\label{BD2na} 
Let $G=\BD_{2n}(a)$ be a small binary dihedral group and $E$ be the exceptional divisor in $G$-$\HILB(\C^2)$, the minimal resolution of $\C^2/G$. Then $E$ has the following Dynkin diagram of type $D$:
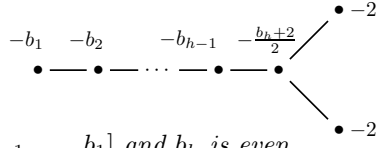
\begin{center}
\begin{pspicture}(0,-0.5)(6,1)
\scalebox{0.8}{
	\psset{nodesep=3pt}
	\rput(0,0){\rnode{p1}{$\bullet$}}\rput(-0.2,0.5){$-b_{1}$}
	\rput(1,0){\rnode{p2}{$\bullet$}}\rput(0.8,0.5){$-b_{2}$}
	\rput(2,0){\rnode{p3}{$\cdots$}}
	\rput(3,0){\rnode{p4}{$\bullet$}}\rput(2.5,0.5){$-b_{h-1}$}
	\rput(4,0){\rnode{p5}{$\bullet$}}\rput(3.8,0.5){$-\frac{b_{h}+2}{2}$}
	\rput(5,1){\rnode{p6}{$\bullet$}}\rput(5.4,1){$-2$}
	\rput(5,-1){\rnode{p7}{$\bullet$}}\rput(5.4,-1){$-2$}
	\ncline{-}{p1}{p2}\ncline{-}{p3}{p4}
	\ncline{-}{p2}{p3}\ncline{-}{p4}{p5}
	\ncline{-}{p5}{p6}\ncline{-}{p5}{p7}
	}
\end{pspicture}
\end{center}	
where $\frac{2n}{a}=[b_{1},\ldots,b_{h-1},b_{h},b_{h-1},\ldots,b_{1}]$ and $b_{h}$ is even.
\end{prop}

\begin{proof}
Let $Y=\Hilb{H}{\C^2}$, $\pi\colon Y\to Y/\!\Span{\beta}$ be the quotient map, and $E_{i}$, $i=1,\ldots,2h-1$, be the exceptional curves in $Y$ with $E_{i}^2=E_{2h-i-1}^2=-b_{i}$. Denote by $E'_{j}$ the exceptional curves in $Y/\!\Span{\overline{\beta}}$ for $j=1,\dots,h$ and by $\widetilde{E}_{i}$ the strict transform of $E'_{i}$. 

For $i\neq h$ the action of $\beta$ on $E_i$ has no fix points we have $\widetilde{E}_i^2=-b_i$. It remains to check what happens on the middle rational curve $E_h$. The curve $E_{h}$ is covered by the affine charts $Y_{h}\cong\C^2_{(\lambda,\mu)}$ and $Y_{h+1}\cong\C^2_{(\lambda',\mu')}$, where $\lambda=x^{i}/y^j$, $\mu=y^q/x^q$, $\lambda=x^{q}/y^q$ and $\mu=y^i/x^j$. Using the fact that $\frac{1}{2n}(i,j)$ and $\frac{1}{2n}(q,q)$ are consecutive points in the Newton polygon of the lattice $L:=\Z^2+\frac{1}{2n}(1,a)\cdot\Z$, we know that $i+j=b_{h}q$, and the action of $\beta$ is $(\lambda,\mu)\mapsto((-1)^j\lambda\mu^{b_{h}},(-1)^{q}1/\mu)$.

Let us study this action in detail. First note that since $i$ and $j$ have the same parity and $i+j=b_{h}q$, we therefore have that $b_{h}$ and $q$ cannot be both odd. In addition, $j$ and $q$ cannot be both even. Indeed, let $(0,2n)$, $(1,a), \ldots, (r,s)$, $(u,v),\ldots, (j,i)$, $(q,q)$, $(i,j), \ldots, (a,1)$, $(2n,0)$ be the sequence of points in the boundary of the Newton polygon (all of them divided by $\frac{1}{2n}$). If $q$ and $j$ (and therefore $i$) are even then $u=b_{h-1}j-q$ and $r=b_{m-2}u-j$ are also even. By induction we deduce that $1$ is even, which is absurd. Hence, the only possibilities for fixed locus of the action of $\beta$ on $Y$ are shown in Table \ref{table!qreflec}. 

\begin{table}[htdp]
\begin{center}
\begin{tabular}{|c|c|c|l|c|}
\hline
$q$ & $j$ & $b_{h}$ & Fixed locus & $\BD_{2n}(a)$ \\
\hline
even & odd 	& even 	& Points $(0,1)$ and $(0,-1)$ fixed & small \\
	 &	     	& odd   	& Point $(0,1)$ and line $\mu=-1$ fixed & not small  \\
\hline
odd	& even	& $\equiv0$ (4) & Lines $\mu=i$ and $\mu=-i$ fixed & not small \\
	& 		& $\equiv2$ (4) & Points $(0,i)$ and $(0,-i)$ fixed & small \\
\hline
odd	& odd	& $\equiv0$ (4) & Points $(0,i)$ and $(0,-i)$ fixed & small \\
	& 		& $\equiv2$ (4) & Lines $\mu=i$ and $\mu=-i$ fixed & not small 	 \\
\hline
\end{tabular}
\end{center}
\caption{Fixed locus in $Y=\Hilb{H}{\C^2}$ by the action of $\beta$.}
\label{table!qreflec}
\end{table}%
In the case when the group $\BD_{2n}(a)$ is small, we see that $b_{h}$ is even and  $Y/\!\Span{\beta}$ has two singular $A_{1}$ points along $E'_{h}\cong\mathbb{P}^1$ and $E'_{h}=-b_{h}/2$. Let $f:\widetilde{Y}\to Y/\!\Span{\beta}$ the resolution of these two $A_{1}$ singularities denoting by $C_{1}$, $C_{2}$ the corresponding rational curves in $\widetilde{Y}$. Then $-b_{h}/2 = (\widetilde{E}_{h}+C_{1}/2+C_{2}/2)^2$=$(\widetilde{E}_{h})^2+C_{1}^2/4+C_{2}^2/4+\widetilde{E}_{h}C_{1}+\widetilde{E}_{h}C_{2}$=$(\widetilde{E}_{h})^2+1$, so that $(\widetilde{E}_{h})^2=-(b_{h}+2)/2$.
\end{proof}

\begin{rem} (i) If $a=1$ the group $\BD_{2n}(1)$ is abelian, and by Proposition \ref{CritSmall} it is small if and only if $n$ is odd. The continued fraction $\frac{2n}{1}=[2n]$ has only one term $b_h=2n$ with $-\frac{b_h+2}{2}=-(n+1)$. Then $E$ has a type $A$ Dynkin diagram of the form
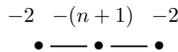
\begin{center}
\begin{pspicture}(0,0)(2,0.5)
\scalebox{0.8}{
	\psset{nodesep=3pt}
	\rput(0,0){\rnode{p1}{$\bullet$}}\rput(-0.3,0.5){$-2$}
	\rput(1,0){\rnode{p2}{$\bullet$}}\rput(0.9,0.5){$-(n+1)$}
	\rput(2,0){\rnode{p3}{$\bullet$}}\rput(2.1,0.5){$-2$}
	\ncline{-}{p1}{p2}\ncline{-}{p2}{p3}
	}
\end{pspicture}
\end{center}					
(ii) If we consider the non-small group $\BD_{12}(5)$, the quotient variety $X/\!\Span{\beta}$ is non-singular and the exceptional divisor is of the form
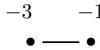
\begin{center}
\begin{pspicture}(0,0)(1.25,0.2)
\scalebox{0.8}{
	\psset{nodesep=3pt}
	\rput(0,0){\rnode{p1}{$\bullet$}}\rput(-0.2,0.5){$-3$}
	\rput(1,0){\rnode{p2}{$\bullet$}}\rput(1,0.5){$-1$}
	\ncline{-}{p1}{p2}
	}
\end{pspicture}
\end{center}
To obtain the minimal resolution we need to contract the $-1$-curve, obtaining a single $\mathbb{P}^1$ with selfintersection $-2$. Explicitly, the subgroup of quasireflections is generated by $H=\Span{\alpha\beta,\alpha^3\beta,\alpha^5\beta,\alpha^7\beta,\alpha^9\beta,\alpha^{11}\beta}$ which has order 12. The quotient is $\C^2_{x,y}/H\cong\C^2_{u,v}$ where $u=x^6-y^6$ and $v=xy$. Now the action of $\alpha$ in the new coordinates is $(u,v)\mapsto(\varepsilon^6u,\varepsilon^6v)$, i.e.\ it is of type $\frac{1}{12}(6,6)$. Since $\C^2_{u,v}/\frac{1}{12}(6,6)\cong\C^2_{u,v}/\frac{1}{2}(1,1)$, the exceptional divisor in the minimal resolution of $\C^2/\BD_{12}(5)$ consists of a single $\mathbb{P}^1$ with selfintersection $-2$. 
\end{rem}

\section{$G$-Hilb and $G$-graphs}\label{GHilbGgraphs}

We start by describing the $G$-invariant Hilbert scheme $G$-Hilb which motivates the definition of $G$-graph.

\begin{defn}\label{GHilb} Let $G\subset\GL(n,\C)$ be a finite subgroup. A $G$-{\em cluster} is a $G$-invariant zero dimensional subscheme $\mathcal{Z}\subset\C^n$ for which $\mathcal{O}_{\mathcal{Z}}$ is isomorphic to the regular representation of $G$ as a $\C G$-module. The $G$-{\em Hilbert scheme} $\Hilb{G}{\C^n}$ is the moduli space parametrising $G$-clusters. 
\end{defn}

Recall that the regular representation $\C[G]\cong\bigoplus_{\rho\in\text{Irr}G}\rho^{\dim(\rho)}$ and $\mathcal{O}_{\mathcal{Z}}\cong\C[x_1.\ldots,x_n]/I_\mathcal{Z}$ where $I_\mathcal{Z}$ is the ideal defining $\mathcal{Z}$. Thus we may pick a vector space basis of $\OO_\mathcal{Z}$ that contains $\dim(\rho)$ elements in each $\rho\in\Irr G$. To describe a distinguished basis of $\mathcal{O}_{\mathcal{Z}}$ with this property, it is convenient to use the notion of $G$-graph. 

\begin{defn}\label{defnGgraph} Let $G\subset\GL(n,\C)$ be a finite subgroup. A {\em G-graph} is a subset $\Gamma\subset\C[x_{1},\ldots,x_{n}]$ satisfying the following conditions:
\begin{enumerate}
\item It contains $\dim(\rho)$ number of elements in each irreducible representation $\rho$. 
\item If a monomial $x_{1}^{\lambda_1}\cdots x_{n}^{\lambda_n}$ is a summand of a polynomial $P\in\Gamma$, then for every $0\leq \mu_j\leq\lambda_j$ the monomial $x_{1}^{\mu_1}\cdots x_{n}^{\mu_n}$ must be a summand of some polynomial $Q_{\mu_{1},\ldots,\mu_{n}}\in\Gamma$.
\end{enumerate}
\end{defn}

\noindent Note that for any $G$-cluster $\mathcal{Z}$ we can choose a basis for the vector space $\mathcal{O}_{\mathcal{Z}}$ which is a $G$-graph. In other words, we can always find a basis of $\mathcal{O}_\mathcal{Z}$ which is minimal in the sense of condition 2 in Definition \ref{defnGgraph}. Indeed, let $\rho\in\Irr G$ of dimension $d$ and let $g_1,\ldots,g_d\in\rho$ be the basis elements in $\Gamma$. Now let $f\in\rho$ and suppose that $f\notin\Gamma$ but $pf\in\Gamma$ for some $p\in\C[x_1,\ldots,x_n]$. Then we have a relation of the form $f\equiv a_1g_1+\ldots+a_dg_d$ modulo $I_\mathcal{Z}$ for some $a_i\in\C$, which implies that $pf\equiv a_ipg_1+\ldots+a_dpg_d$ modulo $I_\mathcal{Z}$. Since $pf\in\Gamma$ we have $a_j\neq0$ for at least one $j$, which allows us to consider the expression $g_j\equiv(1/a_j)f-(a_i/a_j)g_1-\ldots-(a_d/a_j)g_d$ modulo $I_\mathcal{Z}$. Thus we may choose $f$ to be the basis element in $\Gamma$ instead of $g_j$(compare with \cite{Leng}, Chapter 2).

For any $G$-graph $\Gamma$ there exists a set $U_{\Gamma}\subset$ G-Hilb($\C^n$) consisting of all $G$-clusters $\mathcal{Z}$ such that $\mathcal{O}_{\mathcal{Z}}$ admits $\Gamma$ for basis as a vector space. Since being a basis of an open set is an open condition, the set $U_\Gamma$ is open. Therefore, given the set of all possible $G$-graphs $\Gamma$, their union covers $G$-Hilb($\C^n$). 

\begin{exa}\label{exa-graph} If we consider the cyclic group $G=\Span{\frac{1}{5}(1,3)}$ then $\Gamma=\{1,x,x^2,y,xy\}$ is a $G$-graph. For the non-abelian binary dihedral group $D_{4}=\Span{\frac{1}{4}(1,3),\left(\begin{smallmatrix}0&1\\-1&0\end{smallmatrix}\right)}\subset\SL(2,\C)$, $\Lambda=\{1,x,y,x^2+y^2,x^2-y^2,y^3,-x^3,x^4-y^4 \}$ is a $D_{4}$-graph (note that $(x,y),(y^3,-x^3)\in V_1$). 
\end{exa}

\begin{defn} Let $I_\Gamma$ be the ideal generated by the polynomials $f\in\rho$ which are not in $\Gamma$, for any $\rho\in\Irr G$. We say that the $G$-graph $\Gamma$ is {\em represented by the ideal} $I_\Gamma$. Note that $I_\Gamma$ determines $\Gamma$ uniquely. We may also say that a $G$-cluster in an open set $U_\Gamma\subset\Hilb{G}{\C^n}$ is represented by the $G$-graph $\Gamma$. 

The {\em representation} of a $G$-graph $\Gamma$ is the Young diagram in the lattice $M$ consisting of monomials which are summands of polynomials in $\Gamma$ modulo $I_\Gamma$. This representation reflects the nature of a $G$-graph and is useful to describe visually how a $G$-graph varies through $\Hilb{G}{\C^2}$. 
\end{defn}

\begin{exa}\label{ex!twins} In Example \ref{exa-graph} the $G$-graph $\Gamma$ is represented by the ideal $I_\Gamma=\langle x^3,x^2y,y^2\rangle$, and $\Lambda$ is represented by the ideal $I_\Lambda=\langle xy,x^4+y^4\rangle$. The representations of $\Gamma$ and $\Lambda$ are shown in Figure \ref{Fig-graph}. 

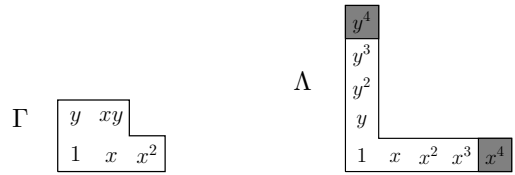
\begin{figure}[htbp]
\begin{center}
\Large{
\begin{pspicture}(0,0)(7,2)
\scalebox{0.85}{
\rput(5,-0.25){
	\scalebox{0.65}{
	\psline(0,0)(4,0)(4,0.8)(0.8,0.8)(0.8,4)(0,4)(0,0)
	\psline[fillstyle=solid,fillcolor=gray](3.2,0)(4,0)(4,0.8)(3.2,0.8)(3.2,0)
	\psline[fillstyle=solid,fillcolor=gray](0,3.2)(0,4)(0.8,4)(0.8,3.2)(0,3.2)
	\rput(0.4,0.4){$1$}
	  \rput(1.2,0.33){$x$}\rput(2,0.4){$x^2$}\rput(2.8,0.4){$x^3$}\rput(3.6,0.4){$x^4$}
	  \rput(0.4,1.2){$y$}\rput(0.4,2){$y^2$}\rput(0.4,2.8){$y^3$}\rput(0.4,3.6){$y^4$}
	  }}
\rput(0.5,-0.25){
	\scalebox{0.7}{
	\psline(0,0)(2.4,0)(2.4,0.8)(1.6,0.8)(1.6,1.6)(0,1.6)(0,0)
	\rput(0.4,0.4){$1$}
	  \rput(1.2,0.33){$x$}\rput(2,0.4){$x^2$}
	   \rput(1.2,1.2){$xy$}\rput(0.4,1.2){$y$}
	  }}
	  }
\rput(0,0.5){\normalsize $\Gamma$}
\rput(3.75,1){\normalsize $\Lambda$}
\end{pspicture}	}
\caption{Representation of the $G$-graphs $\Gamma$ and $\Lambda$.}
\label{Fig-graph}
\end{center}
\end{figure}
The relation $x^4+y^4\in I_\Lambda$ identifies $x^4$ and $y^4$ in $\C[x,y]/I_\Lambda$ and we say that $x^4$ and $y^4$ are {\em ``twins''}. 
\end{exa}

\subsection{The cyclic case}\label{seq!cyclic} 

In this section we recall the construction of the minimal resolution $Y:=\Hilb{H}{\C^2}$ of $\C^2/H$ for an abelian subgroup $H\subset\GL(2,\C)$ using Hirzebruch--Jung continued fractions and the relation with $H$-graphs (see \cite{Kidoh}). 

As it was introduced in Section \ref{Cyc-quot}, $Y$ is the union of $m+1$ open sets $Y_i\cong\C^2$ defined by two consecutive boundary lattice points $e_i,e_{i+1}$ of $L$. If we denote $e_i=\frac{1}{k}(r,s)$ and $e_{i+1}=\frac{1}{k}(u,v)$ then the corresponding open set in $Y$ is of the form $Y_i\cong\mathbb{C}^2_{(\xi_i,\eta_i)}$, where $\xi_i=x^{s}/y^{r}$ and $\eta_i=y^{u}/x^{v}$. Every point $(\xi_{i},\eta_{i})\in Y_i$ corresponds to the $H$-cluster $\mathcal{Z}_{\xi_{i},\eta_{i}}$ defined by the ideal 
\[ I_{\xi_{i},\eta_{i}}=\langle x^{s}-\xi_iy^{r}, ~ y^{u}-\eta_ix^{v},~ x^{s-v}y^{u-r}-\xi_i\eta_i\rangle \]

The $H$-graph $\Gamma_i$ corresponding to the open set $Y_i$ is determined by setting $\xi_{i}=\eta_{i}=0$. In other words, $\Gamma$ is represented by the ideal $I_\Gamma=I_{0,0}=\langle x^s, y^u, x^{s-v}y^{u-r}\rangle$ which pictorially is given by the following ``stair'' shape:
 
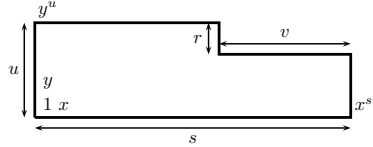
\begin{center}
\begin{pspicture}(0,-0.2)(5,1.25)
\scalebox{0.7}{
	\psline[linewidth=1.5pt](0,0)(6,0)(6,1.2)(3.5,1.2)(3.5,1.8)(0,1.8)(0,0)
	\psline{<->}(-0.2,0)(-0.2,1.8)\rput(-0.4,0.9){$u$}
	\psline{<->}(0,-0.2)(6,-0.2)\rput(3,-0.4){$s$}
	\psline{<->}(3.3,1.2)(3.3,1.8)\rput(3.1,1.5){$r$}
	\psline{<->}(3.5,1.4)(6,1.4)\rput(4.75,1.6){$v$}
	 \rput(0.25,0.25){1}\rput(0.55,0.22){$x$}\rput(0.25,0.65){$y$}
	\rput(6.25,0.25){$x^s$}\rput(0.25,2.05){$y^u$}
	}
\end{pspicture}
\end{center}	

For $(\xi_{i},\eta_{i})\in\mathbb{C}^2$, any monomial in $\mathbb{C}[x,y]$ can be written in terms of elements in $\Gamma$ modulo the ideal $I_{\mathcal{Z}_{\xi_{i},\eta_{i}}}$. In other words, $\Gamma$ is a basis for the vector space $\mathbb{C}[x,y]/I_{\mathcal{Z}_{\xi_{i},\eta_{i}}}$. Note also that $k=su-rv$, so the number of elements in $\Gamma$ agrees with the order of the group. Thus $\C^2_{\xi_{i},\eta_{i}}$ is an open set in $\Hilb{H}{\C^2}$. 

\begin{exa}\label{exa-A12}  Consider the group $H=\Span{\frac{1}{12}(1,7)}$. We have $\frac{12}{7}$=[2,4,2] and therefore $Y=\Hilb{H}{\C^2}$ is of the form $Y=Y_0\cup Y_{1}\cup Y_{2}\cup Y_3$, where $Y_i\cong\mathbb{C}^2_{(\xi_i,\eta_i)}$, $i=0,\ldots,3$ (see Figure \ref{Res12(7)}). The corresponding $H$-clusters for these affine pieces are defined by the ideals:
\begin{align*}
I_{\xi_{0},\eta_{0}}&=\langle x^{12}-\xi_0, y-\eta_{0}x^7\rangle &  
I_{\xi_{2},\eta_{2}}&=\langle x^{2}-\xi_2y^2, y^7-\eta_{2}x, xy^5-\xi_{2}\eta_{2}\rangle \\
I_{\xi_{1},\eta_{1}}&=\langle x^{7}-\xi_1y, y^2-\eta_{1}x^2, x^5y-\xi_{1}\eta_{1}\rangle &
I_{\xi_{3},\eta_{3}}&=\langle x-\xi_3y^7, y^{12}-\eta_{3}\rangle
\end{align*}
and the representation of the $H$-graphs is shown in Figure \ref{GGraphs12(7)}.
\begin{figure}[ht]
\begin{center}
\begin{pspicture}(0,0)(8,1)
	\psset{nodesep=3pt,arcangle=30}
\scalebox{0.8}{
	\rput(0,0){\rnode{1}{}}\rput(4,0){\rnode{2}{}}	
	\rput(3,0){\rnode{3}{}}\rput(7,0){\rnode{4}{}}
	\rput(6,0){\rnode{5}{}}\rput(10,0){\rnode{6}{}}
	\ncarc{-}{1}{2}\ncarc{-}{3}{4}
	\ncarc{-}{5}{6}
	\psline[linestyle=dotted](0.5,0)(0.1,1.5)
	\psline[linestyle=dotted](9.5,0)(9.9,1.5)
	\rput(-0.5,0.8){$\xi_{0}=x^{12}$}\rput(1.1,0.15){$\eta_{0}$}
	\rput(2.85,0.15){$\xi_{1}$}\rput(4.2,0.15){$\eta_{1}$}
	\rput(5.85,0.15){$\xi_{2}$}\rput(7.2,0.15){$\eta_{2}$}
	\rput(9,0.15){$\xi_{3}$}
	\rput(10.5,0.8){$\eta_{3}=y^{12}$}
	\rput(2,0.8){$x^7:y$}
	\rput(5,0.8){$x^2:y^2$}
	\rput(8,0.8){$x:y^7$}
	\psline{-}(1,0.41)(0.85,0.45)\psline{-}(1,0.41)(0.87,0.3)
	\psline{-}(3,0.41)(3.15,0.45)\psline{-}(3,0.41)(3.13,0.3)
	\psline{-}(4,0.41)(3.85,0.45)\psline{-}(4,0.41)(3.87,0.3)
	\psline{-}(6,0.41)(6.15,0.45)\psline{-}(6,0.41)(6.13,0.3)
	\psline{-}(7,0.41)(6.85,0.45)\psline{-}(7,0.41)(6.87,0.3)
	\psline{-}(9,0.41)(9.15,0.45)\psline{-}(9,0.41)(9.13,0.3)
	\psline{-}(0.32,0.7)(0.27,0.55)\psline{-}(0.32,0.7)(0.42,0.59)
	\psline{-}(9.68,0.7)(9.73,0.55)\psline{-}(9.68,0.7)(9.58,0.59)
	}
\end{pspicture}
\caption{Resolution of singularities $Y$ of the cyclic singularity of type $\frac{1}{12}(1,7)$.}
\label{Res12(7)}
\end{center}
\end{figure}
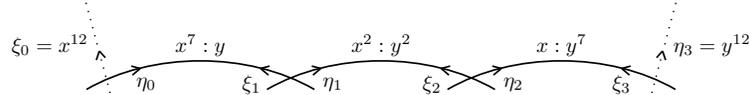
\begin{figure}[htbp]
\begin{center}
\begin{pspicture}(0,0)(9,3.5)
\scalebox{0.7}{
\rput(0,0){\rnode{I0}{\psline(0.2,0)(6.1,0)(6.1,0.5)(0.2,0.5)(0.2,0)
	\rput(0.37,0.23){$1$}\rput(0.75,0.19){$x$}
	\rput(1.25,0.25){$x^2$}\rput(1.75,0.25){$x^3$}
	\rput(2.25,0.25){$x^4$}\rput(2.75,0.25){$x^5$}
	\rput(3.25,0.25){$x^6$}\rput(3.75,0.25){$x^7$}
	\rput(4.25,0.25){$x^8$}\rput(4.75,0.25){$x^9$}
	\rput(5.25,0.25){$x^{10}$}\rput(5.8,0.25){$x^{11}$}}}
\rput(6.5,0){\rnode{I1}{\psline(0,0)(3.8,0)(3.8,0.5)(2.8,0.5)(2.8,1.05)(0,1.05)(0,0)
	\rput(0.2,0.23){$1$}\rput(0.6,0.19){$x$}
	\rput(1.15,0.25){$x^2$}\rput(1.78,0.25){$x^3$}
	\rput(2.45,0.25){$x^4$}\rput(3.1,0.25){$x^5$}
	\rput(3.6,0.25){$x^6$}\rput(0.20,0.69){$y$}			
	\rput(0.6,0.69){$xy$}\rput(1.15,0.75){$x^2y$}
	\rput(1.78,0.75){$x^3y$}\rput(2.45,0.75){$x^4y$}}}
\rput(10.7,0){\rnode{I2}{\psline(0,0)(1.1,0)(1.1,2.5)(0.5,2.5)(0.5,3.5)(0,3.5)(0,0)
	\rput(0.25,0.25){$1$}\rput(0.75,0.25){$x$}
	\rput(0.25,0.75){$y$}\rput(0.75,0.75){$xy$}
	\rput(0.25,1.25){$y^2$}\rput(0.75,1.25){$xy^2$}
	\rput(0.25,1.75){$y^3$}\rput(0.75,1.75){$xy^3$}
	\rput(0.25,2.25){$y^4$}\rput(0.75,2.25){$xy^4$}
	\rput(0.25,2.75){$y^5$}\rput(0.25,3.25){$y^6$}}}
\rput(12.2,0){\rnode{I3}{\psline(-0.1,0)(0.5,0)(0.5,6)(-0.1,6)(-0.1,0)
	\rput(0.21,0.25){$1$}\rput(0.21,0.75){$y$}
	\rput(0.25,1.25){$y^2$}\rput(0.25,1.75){$y^3$}
	\rput(0.25,2.25){$y^4$}\rput(0.25,2.75){$y^5$}
	\rput(0.25,3.25){$y^6$}\rput(0.25,3.75){$y^7$}
	\rput(0.25,4.25){$y^8$}\rput(0.25,4.75){$y^9$}
	\rput(0.25,5.25){$y^{10}$}\rput(0.25,5.75){$y^{11}$}}}
\rput(3.25,1){$\Gamma_{0}$}\rput(8.25,1.6){$\Gamma_{1}$}\rput(10.7,4){$\Gamma_{2}$}\rput(11.5,6){$\Gamma_{3}$}
	}
\end{pspicture}
\caption{$H$-graphs for the group $H=\Span{\frac{1}{12}(1,7)}$.}
\label{GGraphs12(7)}
\end{center}
\end{figure}
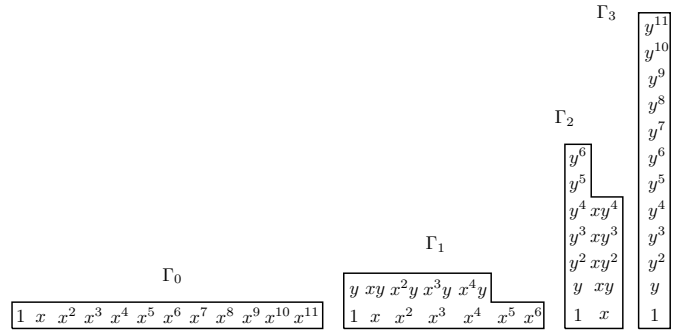	
\end{exa}

\section{$G$-graphs for $\BD_{2n}(a)$ groups}\label{SecGraphs}

Let $G=\BD_{2n}(a)$ be a small binary dihedral group and let $H$ be the maximal normal cyclic subgroup of $G$. As we have seen in Section \ref{Resoln}, the minimal resolution $Y$ of $\C^2/G$ is obtained by acting with $\beta$ on $\Hilb{H}{\C^2}$. The $G$-graphs are constructed in the same way by translating the action of $\beta$ into the $H$-graphs. 

The symmetry along the coordinates of the exceptional divisor $E\subset\Hilb{H}{\C^2}=\bigcup_{i=0}^{2h}Y_{i}$ implies that $\beta$ identifies $Y_{i}$ with $Y_{2h-i}$, as well as identifying the corresponding $H$-graphs $\Gamma_i$ and $\Gamma_{2h-i}$. The union $\Gamma\cup\beta(\Gamma)$ of two $H$-graphs identified by $\beta$ is what we call a {\em $qG$-graph}. 

\subsection{$qG$-graphs}\label{Sect:qGraphs}

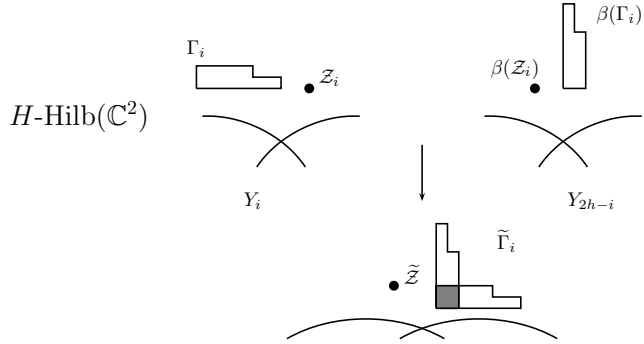
\begin{figure}[htbp]
\begin{center}
\begin{pspicture}(0,0)(6,4.5)
	\psset{nodesep=3pt,arcangle=30,dotsize=5pt}
\scalebox{0.75}{
	\rput(2,0){\rnode{P1}{}}\rput(5,0){\rnode{P2}{}}	
	\rput(4,0){\rnode{P3}{}}\rput(7,0){\rnode{P4}{}}
	\ncarc{-}{P1}{P2}\ncarc{-}{P3}{P4}
	\rput(4.75,0.6){
		\psline(0,0)(1.5,0)(1.5,0.2)(1,0.2)(1,0.4)(0,0.4)(0,0)
		\psline(0,0)(0.4,0)(0.4,1)(0.2,1)(0.2,1.5)(0,1.5)(0,0)
		\psline[fillstyle=solid,fillcolor=gray](0,0)(0.4,0)(0.4,0.4)(0,0.4)(0,0)
		}
	\psdot(4,1)\uput[30](4,1){$\widetilde{\mathcal{Z}}$}
	\rput(6,1.8){$\widetilde{\Gamma}_{i}$}
	
	\rput(0.5,4){\rnode{Q1}{}}\rput(2.5,3){\rnode{Q2}{}}	
	\rput(1.5,3){\rnode{Q3}{}}\rput(3.5,4){\rnode{Q4}{}}
	\ncarc{-}{Q1}{Q2}\ncarc{-}{Q3}{Q4}
	\rput(0.5,4.5){
		\psline(0,0)(1.5,0)(1.5,0.2)(1,0.2)(1,0.4)(0,0.4)(0,0)
		}
	\psdot(2.5,4.5)\uput[30](2.5,4.5){$\mathcal{Z}_{i}$}
	\rput(0.5,5.2){$\Gamma_{i}$}
	\rput(1.5,2.5){$Y_{i}$}
	
	\rput(5.5,4){\rnode{R1}{}}\rput(7.5,3){\rnode{R2}{}}	
	\rput(6.5,3){\rnode{R3}{}}\rput(8.5,4){\rnode{R4}{}}
	\ncarc{-}{R1}{R2}\ncarc{-}{R3}{R4}
	\rput(7,4.5){
		\psline(0,0)(0.4,0)(0.4,1)(0.2,1)(0.2,1.5)(0,1.5)(0,0)
		}
	\psdot(6.5,4.5)\uput[120](6.5,4.5){$\beta(\mathcal{Z}_{i})$}
	\rput(8,5.8){$\beta(\Gamma_{i})$}
	\rput(7.5,2.5){$Y_{2h-i}$}

	\psline{->}(4.5,3.5)(4.5,2.5)	
	\rput(-1.5,4){\Large $\Hilb{H}{\C^2}$
	}}
\end{pspicture}
\caption{Action of $\beta$ on $\Hilb{H}{\C^2}$ in terms of $H$-graphs.}
\label{invo}
\end{center}
\end{figure}

Let $\mathcal{Z}_{i}$ be an $H$-cluster in $Y_{i}$ and $\beta(\mathcal{Z}_{i})$ its image under $\beta$ in $Y_{2h-i}$, with ideals $I_{\mathcal{Z}_{i}}$ and $I_{\beta(\mathcal{Z}_{i})}$ respectively. Denote also by $\widetilde{\mathcal{Z}}$ the point in the quotient $\widetilde{X}:=H$-Hilb$(\mathbb{C}^2)$/$\Span{\beta}$ corresponding to the orbit $\{\mathcal{Z}_{i},\beta(\mathcal{Z}_{i})\}$ (see Figure \ref{invo}). Suppose that $\mathcal{Z}_{i}$ and $\beta(\mathcal{Z}_{i})$ are not one of the two points in $\Hilb{H}{\C^2}$ fixed by $\beta$, which means that $\widetilde{\mathcal{Z}}$ is not one of the singular $A_{1}$ points $P_1,P_2\in\widetilde{X}$. Now $\pi:\Hilb{G}{\C^2}\to\widetilde{X}$ is the minimal resolution with exceptional divisor $E$ consisting of two rational curves. Therefore, there is an isomorphism $\Hilb{G}{\C^2}\backslash E\cong\widetilde{X}\backslash\{P_1,P_2\}$, which means that there exists a unique $G$-cluster $\mathcal{Z}$ corresponding to $\widetilde{\mathcal{Z}}$. 

As clusters in $\mathbb{C}^2$, we have that $\mathcal{Z}\supset\mathcal{Z}_{i}\cup\beta(\mathcal{Z}_{i})$, or equivalently $I_{\mathcal{Z}} \subset I_{\mathcal{Z}_{i}} \cap I_{\beta(\mathcal{Z}_{i})}$. In terms of graphs, if we denote by $\Gamma$, $\Gamma_{i}$ and $\beta(\Gamma_{i})$ the graphs corresponding to the ideals $I_{\mathcal{Z}}$, $I_{\mathcal{Z}_{i}}$ and $I_{\beta(\mathcal{Z}_{i})}$ respectively, we have that
\[ \Gamma \supset \Gamma_{i} \cup \beta(\Gamma_{i})=\widetilde{\Gamma}_{i} \]
where (abusing the notation) by the inclusion we mean that every monomial in $\widetilde{\Gamma}_i$ is a summand of a polynomial in $\Gamma$. 

But notice that $\widetilde{\Gamma}_{i}$ is not a $G$-graph. Indeed, $\Gamma_{i}$ and $\beta(\Gamma_{i})$ have an overlap (common basis elements shaded in Figure \ref{invo}), so the number of elements of $\widetilde{\Gamma}_{i}$ is always smaller than $|G|=2\cdot|H|$. Nevertheless, we will see in the next section that there is a unique way of extending $\widetilde{\Gamma}_i$ into a $G$-graph $\Gamma$. We call these new graphs $\widetilde{\Gamma}_i$ {\em ``quasi $G$-graphs"} ($qG$-{\em graphs}).  

For the two fixed points $\mathcal{Z}_1,\mathcal{Z}_2\in E_h\subset\Hilb{H}{\C^2}$ where $\mathcal{Z}_{i}=\beta(\mathcal{Z}_{i})$, we have that $I_{\mathcal{Z}_i}=I_{\beta(\mathcal{Z}_i)}$ and we can choose an $H$-graph $\Gamma_i$ for $\mathcal{Z}_i$ such that $\Gamma_i=\beta(\Gamma_i)$ for $i=1,2$. This implies that the $qG$-graph $\Gamma_i\cup\beta(\Gamma_i)$ is all overlap, and the extension to a $G$-graph in this case is not unique. In fact, for each of the two $\beta$-fixed points there is a projective line of $G$-clusters corresponding to the exceptional curves of the blow-up which do not come from orbits of $H$-graphs. In other words, because of the presence of fixed points, to define an open cover of $G$-Hilb we need to treat the case where $\widetilde{\Gamma}$ is the union of the two middle $H$-graphs separately. This lead us to $G$-graphs of {\em type $C$ and $D$} (see Section \ref{typesCyD}).

\begin{rem} Since every $H$-graph is given by two consecutive points $e_{i}$, $e_{i+1}$ in the boundary of the Newton polygon of $L$, every $qG$-graph is given by the consecutive pair $e_{i}$, $e_{i+1}$ together with its symmetric pair with respect to the diagonal. Therefore, in order to calculate all possible $qG$-graphs, we just have to look at the $qG$-graphs coming from consecutive points of the list: $e_{0}=\frac{1}{2n}(0,2n)$, $e_{1}=\frac{1}{2n}(1,a)$, $\ldots$, $e_{h}=\frac{1}{2n}(q,q)$.
\end{rem}

\begin{figure}[h]
\begin{center}
\begin{pspicture}(0,0)(7.5,4.25)
		\scalebox{0.7}{
\rput(0,0){\rnode{I0}{\psline(0.2,0)(6.3,0)(6.3,0.5)(0.2,0.5)(0.2,0)
	\rput(0.5,0.23){$1$}\rput(1,0.19){$x$}\rput(1.45,0.25){$x^2$}\rput(1.9,0.25){$x^3$}
	\rput(2.4,0.25){$x^4$}\rput(2.9,0.25){$x^5$}\rput(3.4,0.25){$x^6$}\rput(3.9,0.25){$x^7$}
	\rput(4.4,0.25){$x^8$}\rput(4.9,0.25){$x^9$}\rput(5.4,0.25){$x^{10}$}\rput(6,0.25){$x^{11}$}}}
\rput(0.3,0){\rnode{I3}{\psline(-0.1,0)(0.5,0)(0.5,6.1)(-0.1,6.1)(-0.1,0)
	\rput(0.21,0.75){$y$}\rput(0.25,1.25){$y^2$}\rput(0.25,1.75){$y^3$}
	\rput(0.25,2.25){$y^4$}\rput(0.25,2.75){$y^5$}\rput(0.25,3.25){$y^6$}\rput(0.25,3.75){$y^7$}
	\rput(0.25,4.25){$y^8$}\rput(0.25,4.75){$y^9$}\rput(0.25,5.25){$y^{10}$}\rput(0.25,5.75){$y^{11}$}}}
\rput(7,0){\rnode{I1}{\psline(0,0)(4.1,0)(4.1,0.5)(3.1,0.5)(3.1,1)(0,1)(0,0)
	\rput(0.2,0.23){$1$}\rput(0.75,0.19){$x$}\rput(1.45,0.25){$x^2$}\rput(2.1,0.25){$x^3$}
	\rput(2.75,0.25){$x^4$}\rput(3.35,0.25){$x^5$}\rput(3.85,0.25){$x^6$}\rput(0.20,0.69){$y$}			\rput(0.75,0.69){$xy$}\rput(1.45,0.75){$x^2y$}\rput(2.1,0.75){$x^3y$}\rput(2.75,0.75){$x^4y$}}}
\rput(7,0){\rnode{I2}{\psline(0,0)(1.1,0)(1.1,2.9)(0.5,2.9)(0.5,3.9)(0,3.9)(0,0)
	\rput(0.25,1.4){$y^2$}\rput(0.75,1.4){$xy^2$}\rput(0.25,2){$y^3$}\rput(0.75,2){$xy^3$}
	\rput(0.25,2.6){$y^4$}\rput(0.75,2.6){$xy^4$}\rput(0.25,3.1){$y^5$}\rput(0.25,3.6){$y^6$}}}
\rput(-0.5,3){\Large $\widetilde{\Gamma}_{0}$}\rput(10,2){\Large $\widetilde{\Gamma}_{1}$}
	}
\end{pspicture}
\caption{$qG$-graphs for the group $\BD_{12}(7)$.}
\label{qG12(7)}
\end{center}
\end{figure}
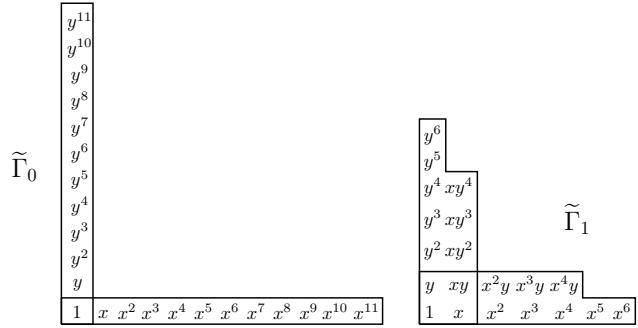

\begin{exa} Continuing with Example \ref{exa-A12}, the action of $\beta$ glues together $\Gamma_{0}$ with  $\Gamma_{3}$ and  $\Gamma_{1}$ with  $\Gamma_{2}$, obtaining the $qG$-graphs $\widetilde{\Gamma}_{0}$ and $\widetilde{\Gamma}_{1}$. See Figure \ref{qG12(7)}.
\end{exa}

Now we look at how two $H$-graphs identified by the action of $\beta$ can merge together into a $qG$-graph. The following proposition shows that there are only two types of gluing.

\begin{prop}\label{TypesAB} Let $G=\BD_{2n}(a)\subset\GL(2,\C)$ a small binary dihedral group. Let $\Gamma$ be the $H$-graph defined by the consecutive Newton polygon points $e_{i}=\frac{1}{2n}(r,s)$ and $e_{i+1}=\frac{1}{2n}(u,v)$. Then corresponding $qG$-graph $\widetilde{\Gamma}$ is either of type $A$ if $s-v>u$, or of type $B$ if $s-v=u-r$. 
\end{prop}
We may also say that the $H$-graph $\Gamma$ is of type $A$ or $B$. Their shape is shown in the following diagram, where the shaded area represents the overlap between the $H$-graphs:
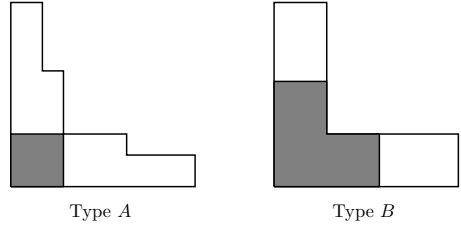
\begin{center}
\begin{pspicture}(0,-0.5)(5,2.75)
\scalebox{0.7}{
	\psline(0,0)(3.5,0)(3.5,0.6)(2.2,0.6)(2.2,1)(1,1)(1,2.2)(0.6,2.2)(0.6,3.5)(0,3.5)(0,0)
	\psline[fillstyle=solid,fillcolor=gray](0,0)(1,0)(1,1)(0,1)(0,0)
	\psline(5,0)(8.5,0)(8.5,1)(6,1)(6,3.5)(5,3.5)(5,0)
	\psline[fillstyle=solid,fillcolor=gray](5,0)(7,0)(7,1)(6,1)(6,2)(5,2)(5,0)
	\rput(1.7,-0.5){Type $A$}\rput(6.7,-0.5){Type $B$}
	}
\end{pspicture}
\end{center}

\begin{proof} Since any $qG$-graph $\widetilde{\Gamma}$ is the union of two $H$-graphs identified by the action of $\beta$, the shape of the $qG$-graph will depend on the relation between $s-v$ and $u$. 

Suppose that $s-v\leq u$. Since $s<v$ this means that $s-v= u-p$ for some $0\leq p<u$. Let $e_{i+2}=\frac{1}{2n}(t,w)$ the next point in the boundary of the Newton polygon. Then since $w=b_{i+1}v-s$ we have that $p=u+v-s=u+(1-b_{i+1})v+w$ which implies $b_{i+1}=2$ (otherwise $p<u-2v+w$, and since $v\geq u$ and $v>w$, then $p<0$ a contradiction). Thus $p=u-v+w$, and applying the same argument to $w$ we obtain that $b_{i+1}=b_{i+2}=\ldots=b_{h}=2$, that is, $\frac{2n}{a}=[b_0,\cdots,b_i,2,\ldots,2,b_i,\ldots,b_0]$. 

Finally, the chain of 2s in the middle of the continued fraction $\frac{2n}{a}$ gives us the value of $p$. Indeed, let $e_{h-1}=\frac{1}{2n}(c,d)$, $e_{h}=\frac{1}{2n}(q,q)$ and $e_{h+1}=\frac{1}{2n}(d,c)$ be the three middle Newton polygon points. Then proceeding with the previous argument we get 
$$\begin{array}{ccl}
p &=& u-d+q \\
  &=& u-(2q-c)+q = u-q+c \\
  &\vdots& \\
  &=& u-t+u = 2u-t \\
  &=& r \\
\end{array}$$
\noindent which gives the $qG$-graph of type B.
\end{proof}

From the previous proof, we can deduce that the distribution of the $qG$-graphs of type $A$ and $B$ in the exceptional locus depend on the number of 2s in the middle of the continued fraction $\frac{2n}{a}$. 

\begin{col}\label{NoTypeB} Let $\widetilde{\Gamma}_{0}, \ldots, \widetilde{\Gamma}_{h-1}$ be the sequence of $qG$-graphs for a given group $G=\BD_{2n}(a)$, and let $\frac{2n}{a}=[b_{0},\ldots,b_{h},\ldots,b_{0}]$.
\begin{itemize}
\item[(i)] If $b_{i}=2$ for $k\leq i\leq h-1$ and $b_{k-1}\neq2$, we have that $\widetilde{\Gamma}_{0}$, $\ldots$, $\widetilde{\Gamma}_{k-1}$ are of type $A$, and $\widetilde{\Gamma}_{k}$, $\ldots$, $\widetilde{\Gamma}_{h-1}$ are of type $B$.
\item[(ii)] There are no type B $qG$-graphs if and only if $b_{h}\neq2$.
\end{itemize}
\end{col}
	
For example, the continued fraction $[\ldots,b,2,2,2,2,2,b,\ldots]$ with $b\neq2$, gives three $qG$-graphs of type $B$. As a consequence we get the following corollary. In the case when $a=2n-1$, i.e.\ $\BD_{2n}(2n-1)\subset$ SL($2,\mathbb{C}$), the coefficients of the continued fraction $\frac{2n}{2n-1}$ are all 2, and every $qG$-graph is of type B.

\subsection{From $qG$-graphs to $G$-graphs}

In this section we construct the $G$-graph corresponding to a given $qG$-graph. First notice that every monomial of the $qG$-graph is in general a summand of a polynomial in the $G$-graph. More precisely, let $\widetilde{\Gamma}=\Gamma\cup\beta(\Gamma)$ be a $qG$-graph, and suppose that $x^ly^m\in\rho_j$ is an element in $\Gamma$ for some $j$. Then by Remark \ref{McKG} we have that $\beta(x^ly^m)=(-1)^mx^my^l\in\rho_{aj}$. If $j\equiv aj$ (mod $2n$) the polynomial that we obtain in the 1-dimensional representation $\rho_j^\pm$ is
\begin{align}
\text{either }~ x^ly^m\pm i(-1)^mx^my^l &\text{  if $n,j$ odd}\label{eqn:polyrep} \\ 
\text{or }~ x^ly^m\pm(-1)^mx^my^l &\text{  otherwise}\notag
\end{align}

If $j\not\equiv aj$ (mod $2n$) we obtain the pair $(x^ly^m,(-1)^mx^my^l)$ in $V_j$. Moreover, since $\Gamma$ is an $H$-graph there always exists another monomial $f\in\Gamma$ belonging to $\rho_{aj}$, which gives a second element in $V_j$, namely $(f,\beta(f))$. These two elements in $V_l$ are the same if and only if $f=x^ly^m$ is contained in the overlap $O:=\Gamma\cap\beta(\Gamma)$.

Therefore, for every representation $\rho_j\in\Irr H$ not contained in the overlap $O$, we obtain one element in every 1-dimensional representation $\rho_j^\pm$ of $G$, and two elements in every 2-dimensional irreducible representation $V_j$ of $G$ as desired. In $O$ we only get one element either in $\rho_j^+$ or $\rho_j^-$, and only one pair in $V_r$. Thus, to form a $G$-graph from a $qG$-graph we need to add new elements in the representations of $G$ arising from representations of $H$ contained in $O$. 

\begin{rem}\label{coeffeq0} Let $I$ be an ideal defining a $G$-cluster and let $\Gamma$ be a choice of $G$-graph represented by $I$. Suppose for simplicity that $\rho$ is a 1-dimensional irreducible representation of $G$ and let $g\in\rho$ be an element in $\Gamma$. Since modulo $I$ the module $S_\rho$ can be considered as a 1-dimensional vector space with basis $g$, for any other polynomial $f\in\rho$ not in $\Gamma$ we have a relation of the form $f\equiv a\cdot g$ (mod $I$) for some $a\in\mathbb{C}$. Doing the same for every irreducible representation of $G$ we obtain the relation of any semi-invariant polynomial as a $\C$-linear combination of elements in $\Gamma$ (mod $I$). We can (and do) choose to take all coefficients equal to zero. In other words, we take the ideal representing $\Gamma$, which geometrically corresponds to look at the intersection points of two of the exceptional curves in $\Hilb{G}{\C^2}$ plus the strict transform of the coordinate axis in $\C^2$. 
\end{rem}

If $\Gamma$ is a $G$-graph and $G$ is abelian we can always choose monomials for the elements in $\Gamma$. In our case we need to deal with combinations of monomials identified by $\beta$, which creates the new phenomenon of {\em twin elements} in $G$-graphs (see Example \ref{ex!twins}). For instance, let $f=x^{i}y^{j}-x^{j}y^{i}$ with $i>j$ be an element in some 1-dimensional representation $\rho$. If $f\notin\Gamma$ by Remark \ref{coeffeq0} we have that $x^iy^{j}=x^{j}y^{i}$ in $\C[x,y]/I$ so both monomials become identified as twins. Moreover, multiplying $f$ by $x^py^q$ for any positive integers $p$ and $q$, we have $x^{i+p}y^{j+q}=x^{j+p}y^{i+q}$, and we get a pair of symmetric {\em twin regions}. We call such an $f$ a {\em twin relation}. Note that for the purpose of counting basis elements, a pair of twin elements will count as a single basis element. 

Let $\widetilde{\Gamma}$ be the $qG$-graph defined by $e_i=\frac{1}{2n}(r,s)$ and $e_{i+1}=\frac{1}{2n}(u,v)$. We denote by $\Gamma(r,s;u,v)$ or simply $\Gamma$ the corresponding $G$-graph and $I$ its defining ideal. 

\subsection{$G$-graphs of type $A$}

Suppose that $\widetilde{\Gamma}$ is a $qG$-graph of type $A$. Then $\widetilde{\Gamma}=\Gamma_{I}\cup\Gamma_{\beta(I)}$ where $I=\langle x^s,y^u,x^{s-v}y^{u-r}\rangle$ and $\beta(I)=\langle x^u,y^s,x^{u-r}y^{s-v}\rangle$. In addition, we have the inequalities $r<u\leq v<s$ and $u<s-v$, coming from the lattice $L$ and the type $A$ condition in Proposition \ref{TypesAB} respectively. The general shape of $qG$-graph of type $A$ is the following:
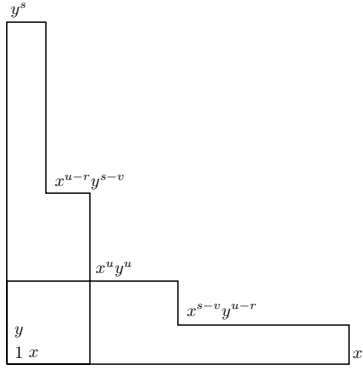
\begin{center}
\begin{pspicture}(0,-0.25)(5,4.5)
	\scalebox{0.65}{
	\psline{-}(0,0)(7,0)(7,0.8)(3.5,0.8)(3.5,1.7)(1.7,1.7)(1.7,3.5)(0.8,3.5)(0.8,7)(0,7)(0,0)
	\psline{-}(0,0)(1.7,0)(1.7,1.7)(0,1.7)(0,0)
	\rput(0.25,0.25){1}
	 \rput(0.55,0.22){$x$}
	 \rput(0.25,0.65){$y$}
	\rput(2.2,1.95){$x^uy^u$}
	\rput(7.25,0.25){$x^s$}\rput(0.25,7.25){$y^s$}
	\rput(4.4,1.1){$x^{s-v}y^{u-r}$}\rput(1.7,3.75){$x^{u-r}y^{s-v}$}
	}
\end{pspicture}
\end{center}

Clearly $\#\widetilde{\Gamma}=2\cdot\#\Gamma_{I}-\# O<2n$, with $O=\Gamma_{I}\cap\Gamma_{\beta(I)}$ the overlap. Then, $\widetilde{\Gamma}$ needs to be extended by $\# O=u^2$ elements to become a $G$-graph, where these extra elements belong to the irreducible representations appearing in $O$. The following lemmas show that we just need to look at some key representations. 

\begin{lem} The polynomial $x^{s-v}y^{u-r}+(-1)^{u-r}x^{u-r}y^{s-v}$ is $G$-invariant, i.e.\ it belongs to $I$. 
\end{lem}

\begin{proof} We need to show that $x^{s-v}y^{u-r}+(-1)^{u-r}x^{u-r}y^{s-v}$ is invariant under the action of $\alpha$ and $\beta$. We know that $ar\equiv s$ and $au\equiv v$ (mod $2n$), then $\alpha(x^{s-v}y^{u-r}+(-1)^{u-r}x^{u-r}y^{s-v})=\varepsilon^{s-v+au-ar}x^{s-v}y^{u-r}+(-1)^{u-r}\varepsilon^{u-r+as-av}x^{u-r}y^{s-v}=x^{s-v}y^{u-r}+(-1)^{u-r}x^{u-r}y^{s-v}$. 

Since $a$ is odd, the congruence $ar\equiv s$ (mod $2n$) also implies that $r$ and $s$ have the same parity (similarly for the pair $u,v$). Therefore $s-v$ and $u-r$ must have the same parity and then $\beta(x^{s-v}y^{u-r}+(-1)^{u-r}x^{u-r}y^{s-v})=x^{s-v}y^{u-r}+(-1)^{u-r}x^{u-r}y^{s-v}$ as required.
\end{proof}

\begin{lem} The monomial $x^{u}y^{u}$ is not in the $G$-graph $\Gamma$. 
\end{lem}

\begin{proof} Suppose that $x^{u}y^{u}\in\Gamma$. Notice that $x^{u}y^{u}, x^{u+v}+(-1)^{u}y^{u+v}\in\rho^{(-1)^u}_{u+v}$, so $x^{u+v}+(-1)^{u}y^{u+v}\in I$ and it forms a twin region. Observe also that $(x^{u}y^{u+1}, (-1)^{u}x^{u+1}y^{u})$ and $(y^{u+v+1}, -x^{u+v+1})$ belong to the same 2-dimensional representation $V_l$, and the monomials $x^{u+v+1}, y^{u+v+1}\in\widetilde{\Gamma}$ don't lie on the overlap $O$. Then there must be another monomial pair in the $qG$-graph $\widetilde{\Gamma}$ outside $O$ that lie in the same representation $V_l$, so it has to be included in $\Gamma$. Therefore $x^{u}y^{u+1}, x^{u+1}y^{u}\in I$, which implies that on the diagonal the $qG$-graph can only be extended by the element $x^{u}y^{u}$. 

We claim that along the sides $\widetilde{\Gamma}$ cannot be extended enough to obtain a $G$-graph. Indeed, by the previous lemma we know that $x^{s-v}y^{u-r}+(-1)^{u-r}x^{u-r}y^{s-v}\in I$ so it forms the pair of twin regions $T_1$. Also $x^r(x^{s-v}y^{u-r}+(-1)^{u-r}x^{u-r}y^{s-v})=x^{s+r-v}y^{u-r}+(-1)^{u-r}x^{u}y^{s-v}\in I$. Now, the type $A$ condition $u<s-v$ implies that $x^{u}y^{s-v}\in I$, and therefore $x^{s+r-v}y^{u-r}\in I$ (same for $x^{u-r}y^{s+r-v}$). Thus the twin regions $T_1$ have each block $r^2$ elements. In the same way, since $x^uy^u\in\Gamma$ the polynomial $x^{u+v}+(-1)^{u}y^{u+v}$ forms the twin region $T_2$, and because $u+u<s$, their size is at most $(u-r)^2-(u-r)$. Thus, the $qG$-graph is extended with at most $r^2+(u-r)^2-(u-r)+1<u^2=\#M$ elements, which is a contradiction.  
\end{proof}

As a consequence we have that $x^{s-v}y^{u-r}+(-1)^{u-r}x^{u-r}y^{s-v}$ creates the twin region $T_1$ of size $r^2$. The following lemma shows that there exists another polynomial creating a second twin region. 

\begin{lem}\label{lem:TwinT2}
The polynomial $x^{r+s}+(-1)^{r}y^{r+s}\in I$ and creates the twin region $T_2$ of size $(u-r)^2$.
\end{lem}

\begin{proof}
It is easy to see that $x^{r}y^{r}, x^{r+s}+(-1)^{r}y^{r+s}\in\rho^{(-1)^r}_{r+s}$. However, $x^{r}y^{r}$ belongs to the $qG$-graph so it must be in $\Gamma$, and therefore $x^{r+s}+(-1)^{r}y^{r+s}\in I$. Now combining the previous lemmas we see that $x^{s+u},y^{s+u}\in I$ and the size of the twin region $T_2$ is equal to $(u-r)^2$.
\end{proof}

\begin{figure}[h]
\begin{center}
\begin{pspicture}(0,0)(6,5.75)
	\scalebox{0.65}{
	\psline{-}(0,0)(7,0)(7,0.8)(3.5,0.8)(3.5,1.7)(1.7,1.7)(1.7,3.5)(0.8,3.5)(0.8,7)(0,7)(0,0)
	\psline{-}(0,0)(1.7,0)(1.7,1.7)(0,1.7)(0,0)
	\psline[fillstyle=solid,fillcolor=lightgray](3.5,0.8)(4.4,0.8)(4.4,1.7)(3.5,1.7)(3.5,0.8)
	\psline[fillstyle=solid,fillcolor=lightgray](0.8,3.5)(0.8,4.4)(1.7,4.4)(1.7,3.5)(0.8,3.5)
	\psline[fillstyle=solid,fillcolor=lightgray](0,0)(0.9,0)(0.9,0.9)(0,0.9)(0,0)
	\psline(7,0)(8.7,0)(8.7,0.8)(7,0.8)(7,0)
	\psline(0,7)(0,8.7)(0.8,8.7)(0.8,7)(0,7)
	\psline[fillstyle=solid,fillcolor=gray](0.9,0.9)(1.7,0.9)(1.7,1.7)(0.9,1.7)(0.9,0.9)
	\psline[fillstyle=solid,fillcolor=gray](7.9,0)(7.9,0.8)(8.7,0.8)(8.7,0)(7.9,0)
	\psline[fillstyle=solid,fillcolor=gray](0,7.9)(0.8,7.9)(0.8,8.7)(0,8.7)(0,7.9)
	\psline[linewidth=2pt]{-}(0,0)(8.7,0)(8.7,0.8)(4.4,0.8)(4.4,1.7)(1.7,1.7)(1.7,4.4)(0.8,4.4)(0.8,8.7)(0,8.7)(0,0)
	\rput(1.3,1.1){$x^r\!y^r$}
	\rput(2.15,1.92){$x^uy^u$}
	\rput(5.2,0.25){$x^{u+v}$}\rput(0.4,5.2){$y^{u+v}$}
	\rput(9.15,0.225){$x^{s+u}$}\rput(0.4,8.95){$y^{s+u}$}
	\rput(5.4,1.05){$x^{s-v+r}y^{u-r}$}\rput(1.8,4.65){$x^{u-r}y^{s-v+r}$}
	\rput(7.25,0.25){$x^s$}\rput(0.25,7.25){$y^s$}
	\rput(1.3,0.45){II}\rput(0.45,1.3){I}
	\rput(0.55,7.6){II}\rput(7.6,0.55){I}
	\rput(5.5,7.5){Twin Regions}
	\rput(5,6){\rnode{T1}{\psline[fillstyle=solid,fillcolor=lightgray](0,0)(0.9,0)(0.9,0.9)(0,0.9)(0,0)}}
	\rput(6.3,6.4){$T_{1}$}
	\rput(5.1,4.5){\rnode{T2}{\psline[fillstyle=solid,fillcolor=gray](0,0)(0.8,0)(0.8,0.8)(0,0.8)(0,0)}}
	\rput(6.3,4.9){$T_{2}$}
	}
\end{pspicture}
\caption{Representation of $\Gamma_A(r,s;u,v)$ from the extension of a $qG$-graph by the elements in the overlap.}
\label{FigTypeA}
\end{center}
\end{figure}
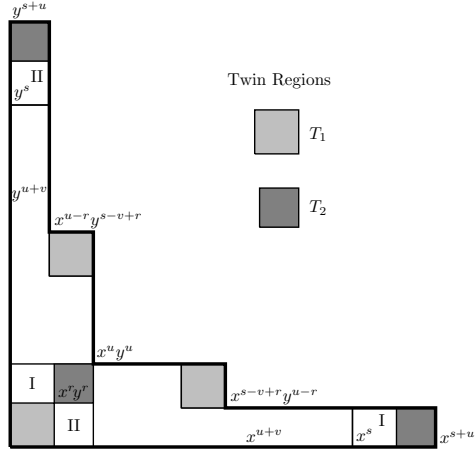

Hence, the extension to a $G$-graph from any type $A$ $qG$-graph has two twin regions $T_1$ and $T_2$, except for $\Gamma_{A}(0,2n;1,a)$, where only $T_2$ appears ($r$ is 0 in this case). Figure \ref{FigTypeA} represents the shape of the $G$-graph we were looking for. 

\begin{prop}\label{TypeA} Every free $G/H$-orbit of an $H$-cluster defined by an $H$-graph of type $A$ is represented by an ideal of the form
\[
I_{A}=\langle x^{u}y^{u}, x^{s-v}y^{u-r}+(-1)^{u-r}x^{u-r}y^{s-v}, x^{r+s}+(-1)^{r}y^{r+s}\rangle.
\]
$G$-graphs represented by ideals of the form $I_A$ are said to be of type $A$ and are denoted by $\Gamma_{A}(r,s;u,v)$.
\end{prop}

\begin{proof}
The previous lemmas show that the $qG$-graph only expands along its sides, and it does not grow further than the twin regions $T_1$ and $T_2$. Let I and II be the regions described in Figure \ref{FigTypeA}. Then, the number of elements added is
\[ \#T_1+\#T_2 + \#\text{I} + \#\text{II} = r^2+(r-u)^2+2r(u-r)=u^2 \]
exactly the number of elements in the overlap.  

It is remaining to prove that $\Gamma$ contains 1 polynomial for each 1-dimensional irreducible representation $\rho_{k}$ and 2 pairs of polynomials for each 2-dimensional irreducible representation $V_{l}$. By the symmetry of $\widetilde{\Gamma}$ every representation not appearing in the overlap will appear twice, as required. Therefore we need only check that the representations of the new extended blocks correspond exactly to the ones in the overlap.

Every representation contained in regions I and II is 2-dimensional, and we have one basis element coming from the overlap $O$ and a second from the extended region. The twin relations $x^{s-v}y^{u-r}+(-1)^{u-r}x^{u-r}y^{s-v}\in I_A$ and $x^{r+s}+(-1)^{r}y^{r+s}\in I_A$ give the correct number of elements in $T_{1}$ and $T_{2}$. We show in the diagram below the configuration of the representations contained in the twin region $T_{1}$. The case of $T_{2}$ is analogous.

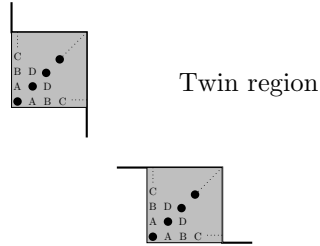
\begin{center}
\begin{pspicture}(0,0.5)(2,3.5)
{\scalebox{0.4}{
	\psset{dotsize=8pt}
\rput(-1,1){\rnode{F1}{
	\psline[linewidth=2pt](1,2.5)(4.5,2.5)(4.5,0)(5.5,0)
	\psline[fillstyle=solid,fillcolor=lightgray](2,2.5)(2,0)(4.5,0)(4.5,2.5)(2,2.5)
	\psdots(2.2,0.2)(2.7,0.7)(3.15,1.15)(3.6,1.6)
	\rput(2.7,0.2){A}\rput(3.2,0.2){B}\rput(3.7,0.2){C}
	\rput(2.2,0.7){A}\rput(2.2,1.2){B}\rput(2.2,1.7){C}
	\rput(3.2,0.7){D}\rput(2.7,1.2){D}
	\psline[linestyle=dotted](4,0.25)(4.5,0.25)
	\psline[linestyle=dotted](2.2,2)(2.2,2.5)
	\psline[linestyle=dotted](3.7,1.7)(4.5,2.5)}}

	\psline[linewidth=2pt](-1,4.5)(-1,8)(-3.5,8)(-3.5,9)
	\psline[fillstyle=solid,fillcolor=lightgray](-1,5.5)(-3.5,5.5)(-3.5,8)(-1,8)(-1,5.5)
	\psdots(-3.3,5.7)(-2.8,6.2)(-2.35,6.65)(-1.9,7.1)
	\rput(-2.8,5.7){A}\rput(-2.3,5.7){B}\rput(-1.8,5.7){C}
	\rput(-3.3,6.2){A}\rput(-3.3,6.7){B}\rput(-3.3,7.2){C}
	\rput(-2.3,6.2){D}\rput(-2.8,6.7){D}
	\psline[linestyle=dotted](-1.5,5.75)(-1,5.75)
	\psline[linestyle=dotted](-3.3,7.5)(-3.3,8)
	\psline[linestyle=dotted](-1.8,7.2)(-1,8)
}}
	\rput(1.75,2.5){Twin region}
\end{pspicture}
\end{center}

\noindent The elements on the diagonals (marked by dots) are in 1-dimensional representations. The rest (marked by letters) are in 2-dimensional representation pairs, that is, a monomial $x^iy^j$ and its partner with respect to the diagonal $x^jy^i$, will form an element in some 2-dimensional representation (we omit the sign here). Elements with the same letter are in the same representation. Since twin symmetric regions count as one, the representations in the overlap are fully present in the extension. 
\end{proof}

\subsection{$G$-graphs of type $B$}

Suppose that the $qG$-graph $\widetilde{\Gamma}=\Gamma_B(r,s;u,v)$ defined by $e_i=\frac{1}{2n}(r,s)$ and $e_{i+1}=\frac{1}{2n}(u,v)$ is of type $B$. Then we can define 
\[ m:=s-v=u-r\] 
which represents the width of $\widetilde{\Gamma}$. The value of $m$ remains constant for every $qG$-graph of type $B$. Indeed, if $e_{i+2}=\frac{1}{2n}(t,w)$ and we consider the next $qG$-graph $\Gamma_{B}(u,v;t,w)$ of type $B$, we have that $w=2v-s$ so that $v-w=s-v=m$. The general shape for a type B $qG$-graph is the following: 

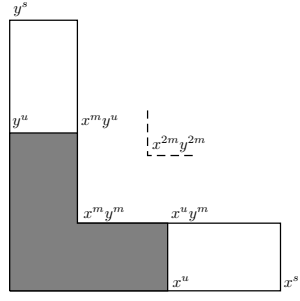
\begin{center}
\begin{pspicture}(0,0)(4,3.75)
\scalebox{0.6}{
	\psline(0,0)(6,0)(6,1.5)(1.5,1.5)(1.5,6)(0,6)(0,0)
	\psline[fillstyle=solid,fillcolor=gray](0,0)(3.5,0)(3.5,1.5)(1.5,1.5)(1.5,3.5)(0,3.5)(0,0)
	\rput(6.25,0.2){$x^s$}\rput(0.25,6.25){$y^s$}
	\rput(3.8,0.2){$x^u$}\rput(0.25,3.75){$y^u$}
	\rput(4,1.7){$x^uy^m$}\rput(2,3.75){$x^my^u$}
	\rput(2.1,1.7){$x^my^m$}
	\rput(-1,-1){	\rnode{P1}{\psline[linestyle=dashed](5,4)(4,4)(4,5)
				\rput(4.7,4.25){$x^{2m}y^{2m}$}}}
	}
\end{pspicture}
\end{center}

Notice that $(x^r,y^r),(y^s,(-1)^sx^s), (x^uy^m,(-1)^mx^my^u)\in V_r$, and since $x^r,y^r\in\widetilde{\Gamma}$ we have $(x^r,y^r)\in\Gamma$. Thus to pick the remaining element of $V_r$ in $\Gamma$ we have a choice between the other two elements in $V_r$. \\

\noindent\underline{Case 1:} Suppose that $(y^s,(-1)^sx^s)\in\Gamma$, that is, $x^{u}y^m,x^my^{u}\in I$. 

\begin{lem}\label{lem:modd}
Let $G=\BD_{2n}(a)$ and let $\Gamma=\Gamma_{B}(r,s;u,v)$ be a $qG$-graph of type $B$. Then the monomial $x^my^m$ is $H$-invariant with $m$ odd, and $x^{2m}y^{2m}$ is $G$-invariant.
\end{lem}

\begin{proof} We have that $\alpha(x^my^m)=\varepsilon^{s-v+a(u-r)}x^my^m=\varepsilon^{s-v+v-s}x^my^m=x^my^m$, so it is $H$-invariant. If we consider $\widetilde{\Gamma}(t,w;q,q)$ to be the last $qG$-graph, i.e.\ the one corresponding to $\Gamma_{h-1}$, we know by definition that $2m=w-t$. On the other hand, $\Gamma_{h-1}$ is an $H$-graph so it has $2n$ elements. Thus $2n=q(w-t)$, and since the group is of type $\BD_{2n}(a)$ by Remark \ref{Brie} we conclude that $m$ is odd. Therefore, $\beta(x^my^m)=-x^my^m$ and the $G$-invariant monomial is $x^{2m}y^{2m}$.
\end{proof}

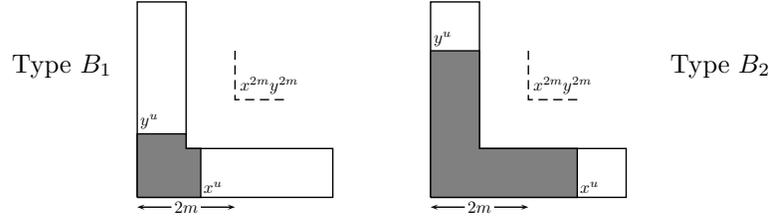
\begin{figure}[h]
\begin{center}
\begin{pspicture}(0,0)(6,2.5)
\scalebox{0.65}{
\rput(0,0){\rnode{P1}{
	\psline(0,0)(4,0)(4,1)(1,1)(1,4)(0,4)(0,0)
	\psline[fillstyle=solid,fillcolor=gray](0,0)(1.3,0)(1.3,1)(1,1)(1,1.3)(0,1.3)(0,0)}
	\psline[linestyle=dashed](3,2)(2,2)(2,3)
	\psline{<-}(0,-0.2)(0.7,-0.2)\rput(1,-0.2){$2m$}\psline{->}(1.3,-0.2)(2,-0.2)
	\rput(1.55,0.2){$x^u$}\rput(0.25,1.55){$y^u$}
	\rput(2.7,2.3){$x^{2m}y^{2m}$}
	}
\rput(6,0){\rnode{P1}{
	\psline(0,0)(4,0)(4,1)(1,1)(1,4)(0,4)(0,0)
	\psline[fillstyle=solid,fillcolor=gray](0,0)(3,0)(3,1)(1,1)(1,3)(0,3)(0,0)}
	\psline[linestyle=dashed](3,2)(2,2)(2,3)
	\psline{<-}(0,-0.2)(0.7,-0.2)\rput(1,-0.2){$2m$}\psline{->}(1.3,-0.2)(2,-0.2)
	\rput(3.25,0.2){$x^u$}\rput(0.25,3.25){$y^u$}
	\rput(2.7,2.3){$x^{2m}y^{2m}$}
	}}
	\rput(-1,1.75){Type $B_{1}$}
	\rput(7.75,1.75){Type $B_{2}$}

\end{pspicture}
\caption{$qG$-graphs of type $B_{1}$ and $B_{2}$ according to the size of the overlap.}
\label{typeB12}
\end{center}
\end{figure}

Now we split up the type $B$ graphs into two different cases as follows: 
\begin{defn} Let $\widetilde{\Gamma}=\widetilde{\Gamma}(r,s;u,v)$ be a $qG$-graph of type $B$. We say $\widetilde{\Gamma}$ is of type $B_{1}$ if $u<2m$, and it is of type $B_{2}$ if $u\geq 2m$. See Figure \ref{typeB12}.
\end{defn}

\begin{lem}\label{B1beforeB2} Let $\widetilde{\Gamma}_{0}, \widetilde{\Gamma}_{1}, \ldots, \widetilde{\Gamma}_{h-1}$ be the sequence of $qG$-graphs for a given group $G=\BD_{2n}(a)$. Let $\frac{2n}{a}=[b_{1},\ldots,b_{h},\ldots,b_{1}]$ and suppose that $b_{i}=2$ for $k\leq i\leq h-1$ and $b_{k-1}\neq2$, for some $0\leq k\leq h-1$. Then $\widetilde{\Gamma}_{k}$ is always of type $B_{1}$, while the rest $\widetilde{\Gamma}_{i}$ for $k<i\leq h-1$ are of type $B_{2}$.
\end{lem}

\begin{proof} By Corollary \ref{NoTypeB} the $qG$-graphs $\widetilde{\Gamma}_{i}$ are of type $A$ for $i\leq k-1$, and of type $B$ for $k\leq i\leq h-1$. Let $\widetilde{\Gamma}_{k-1}=\widetilde{\Gamma}_{A}(c,d;r,s)$ be a $G$-graph of type $A$ and $\widetilde{\Gamma}_{k}=\widetilde{\Gamma}_{B}(r,s;u,v)$ a $G$-graph of type $B$. Then $b_{k+1}=2$ and $d=2s-v$. Now, since $\widetilde{\Gamma}_{k-1}$ is of type $A$ we have $r < d-s = (2s-v)-s = s-v = m$ and then $u=m+r<2m$, so $\widetilde{\Gamma}_{k}$ is of type $B_{1}$. 

Suppose now that $k<h-1$ so that there exists $\widetilde{\Gamma}_{k+1}=\widetilde{\Gamma}_{B}(u,v;w,t)$ another $qG$-graph of type $B$. Then $b_{k+2}=2$ and $w=2u-r$. Moreover, since $m$ is constant for any $qG$-graph of type $B$, we have $m=w-u$. Then $2m = 2(w-u) = 2w - 2u = 2w-(w+r) = w-r$. Therefore $w\geq 2m$, i.e.\ a type $B_{2}$ $qG$-graph. By induction the rest of the $qG$-graphs are also of type $B_{2}$.
\end{proof}
 
\begin{prop}\label{TypeB1} Every free $G/H$-orbit of an $H$-cluster defined by an $H$-graph of type $B_1$ is represented by an ideal of the form
\[ 
I_{B_{1}}=\langle x^{r+s}+(-1)^{r}y^{r+s}, x^{m+s}y^{m-r}+(-1)^{m-r}x^{m-r}y^{m+s}, x^{u}y^{m}, x^{m}y^{u}\rangle.
\]
$G$-graphs represented by ideals of the form $I_{B_1}$ are said to be of type $B_1$ and are denoted by $\Gamma_{B_{1}}(r,s;u,v)$.
\end{prop} 

\begin{figure}[h]
\begin{center}
\begin{pspicture}(0,0)(5,4.5)
\scalebox{0.6}{
	\psline[linestyle=dashed](4,3)(3,3)(3,4)\rput(3.7,3.25){$x^{2m}y^{2m}$}
	\rput(1,1){\psline[fillstyle=solid,fillcolor=gray](0,0)(1.5,0)(1.5,0.5)(0.5,0.5)(0.5,1.5)(0,1.5)(0,0)\rput(0.25,0.25){T}}
	\rput(6.5,0){\psline[fillstyle=solid,fillcolor=gray](0,0)(1.5,0)(1.5,0.5)(0.5,0.5)(0.5,1.5)(0,1.5)(0,0)\rput(0.25,0.25){T}}
	\rput(0,6.5){\psline[fillstyle=solid,fillcolor=gray](0,0)(1.5,0)(1.5,0.5)(0.5,0.5)(0.5,1.5)(0,1.5)(0,0)\rput(0.25,0.25){T}}
	\rput(0,0){\psline[fillstyle=solid](0,0)(1,0)(1,1)(0,1)(0,0)}\rput(0.5,0.5){III}
	\rput(1.5,1.5){\psline[fillstyle=solid](0,0)(1,0)(1,1)(0,1)(0,0)}\rput(2,2){III}
	\rput(0,1){\psline(0,0)(1,0)(1,1.5)(0,1.5)(0,0)\rput(0.5,0.75){I}}
	\rput(5.5,0){\psline(0,0)(1,0)(1,1.5)(0,1.5)(0,0)\rput(0.5,0.75){I}}	
	\rput(1,0){\psline(0,0)(1.5,0)(1.5,1)(0,1)(0,0)\rput(0.75,0.5){II}}
	\rput(0,5.5){\psline(0,0)(1.5,0)(1.5,1)(0,1)(0,0)\rput(0.75,0.5){II}}
	\psline[linewidth=2pt](0,0)(8,0)(8,0.5)(7,0.5)(7,1.5)(2.5,1.5)(2.5,2.5)(1.5,2.5)(1.5,7)(0.5,7)(0.5,8)(0,8)(0,0)
	\rput(1.25,0.25){$x^r$}\rput(0.25,1.25){$y^r$}
	\rput(5.75,0.25){$x^s$}\rput(0.25,5.75){$y^s$}
	\rput(3.1,1.75){$x^uy^m$}\rput(2.1,2.75){$x^my^u$}
	\rput(8,0.75){$x^{m+s}y^{m-r}$}\rput(1.5,7.25){$x^{m-r}y^{m+s}$}
	}
\end{pspicture}
\caption{$G$-graph of type $B_{1}$.}
\label{GtypeB1}
\end{center}
\end{figure}
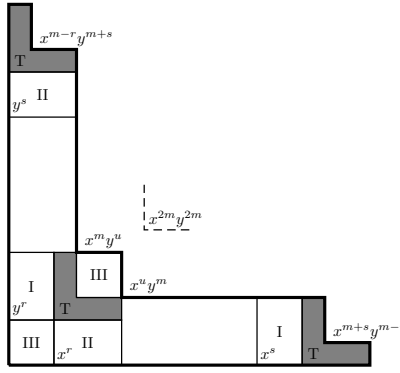

\begin{proof} Since $u<2m$ and $u=m+r$ we know that $r<m$, which implies that the monomial $x^ry^r$ is contained in the overlap $O$. Then $x^{r}y^{r}\in\Gamma$ and we have the relation $x^{r+s}+(-1)^{r}y^{r+s}\in I_{B_1}$, creating the twin region $T$. Using the condition $x^{u}y^m,x^my^{u}\in I_{B_1}$ and multiplying the previous relation by $x^m$ and $y^m$ we get $x^{s+u},y^{s+u}\in I_{B_1}$, which implies that $T$ do not grow further that these monomials. 

Examining the representations $\rho^+_{0}$ and $\rho^-_{0}$, we see that $x^{m+s}y^{m-r}+(-1)^{m-r}x^{m-r}y^{m+s}\in I_{B_1}$ (since always $1\in\rho_0^+$) and $x^{m+s}y^{m-r}-(-1)^{m-r}x^{m-r}y^{m+s}\in I_{B_1}$ (otherwise $x^my^m\in I_{B_1}$ and $\Gamma$ would not have enough elements). Then $x^{m+s}y^{m-r},x^{m-r}y^{m+s}\in I_{B_1}$. All these equations together with the $G$-invariant relations are enough to determine a $G$-graph, hence $I$. 

In Figure \ref{GtypeB1} we show the $G$-graph and the underlying relation between the representations of $G$ for the type $B_{1}$ case. Let $I$, $II$ and $III$ be the regions described in Figure \ref{GtypeB1}. Then the number of new elements are
\[
\#T + \#I + \#II + \#III = m^2-(u-m)^2+2r(u-r)+(u-m)^2=u^2-r^2
\]
which is the same number of elements in the overlap. 

It remains to show that in $\Gamma$ we have the right number of elements in every irreducible representation of $G$. In this case we only have the twin region $T$ and, as in the proof of Proposition \ref{TypeA}, the twin relation $x^{r+s}+(-1)^{r}y^{r+s}\in I_{B_1}$ give the correct number of elements in the representations contained in $T$. Again regions $I$ and $II$ contain only 2-dimensional representations, where the element coming from $O$ and the new element from the extended region make the number of elements match with the dimension of the representation. For the region $III$, representations in the diagonal $x=y$ are 1-dimensional and the rest are 2-dimensional. As shown in the diagram below, elements located in the same place in both blocks have the same $H$-character but opposite $G/H$-character. This implies that every representation $\rho_{j(a+1)}^\pm$ has precisely one element in $\Gamma$.

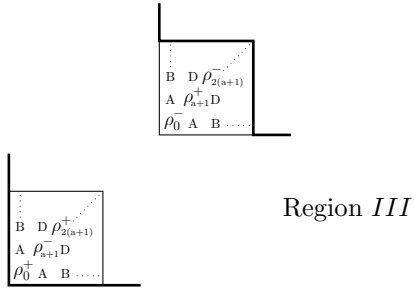
\begin{center}
\begin{pspicture}(0,0.5)(7,4.25)
{\scalebox{0.5}{
	\psset{dotsize=8pt}
\rput(0,1){
	\psline[linewidth=2pt](2,3.5)(2,0)(5.5,0)
	\psline(2,3.5)(2,2.5)(4.5,2.5)(4.5,0)(5.5,0)
	\rput(2.4,0.4){\Large $\rho_0^+$}
	\rput(3,1){\Large $\rho_{\!\text{\scriptsize a+1}}^-$}
	\rput(3.7,1.55){\Large $\rho_{\text{\scriptsize 2(a+1)}}^+$}
	\rput(2.9,0.3){A}\rput(3.5,0.3){B}
	\rput(2.3,0.95){A}\rput(2.3,1.55){B}
	\rput(2.9,1.55){D}\rput(3.5,0.95){D}
	\psline[linestyle=dotted](3.8,0.25)(4.4,0.25)
	\psline[linestyle=dotted](2.3,1.8)(2.3,2.4)
	\psline[linestyle=dotted](3.7,1.7)(4.5,2.5)}

\rput(4,5){
	\psline(2,3.5)(2,0)(5.5,0)
	\psline[linewidth=2pt](2,3.5)(2,2.5)(4.5,2.5)(4.5,0)(5.5,0)
	\rput(2.4,0.4){\Large $\rho_0^-$}
	\rput(3,1){\Large $\rho_{\!\text{\scriptsize a+1}}^+$}
	\rput(3.7,1.55){\Large $\rho_{\text{\scriptsize 2(a+1)}}^-$}
	\rput(2.9,0.3){A}\rput(3.5,0.3){B}
	\rput(2.3,0.95){A}\rput(2.3,1.55){B}
	\rput(2.9,1.55){D}\rput(3.5,0.95){D}
	\psline[linestyle=dotted](3.8,0.25)(4.4,0.25)
	\psline[linestyle=dotted](2.3,1.8)(2.3,2.4)
	\psline[linestyle=dotted](3.7,1.7)(4.5,2.5)}

}}
	\rput(5.5,1.5){Region $III$}
\end{pspicture}
\end{center}

For the 2-dimensional representations we have one element in each block, which give two elements in total as required.
\end{proof}

\begin{prop}\label{TypeB2} Every free $G/H$-orbit of an $H$-cluster defined by an $H$-graph of type $B_2$ is represented by an ideal of the form
\[
I_{B_{2}}=\langle x^{2m}y^{2m}, x^{s+m}, y^{s+m}, x^{u}y^{m}, x^{m}y^{u} \rangle.
\] 
$G$-graphs represented by ideals of the form $I_{B_2}$ are said to be of type $B_2$ and are denoted by $\Gamma_{B_{2}}(r,s;u,v)$.
\end{prop}
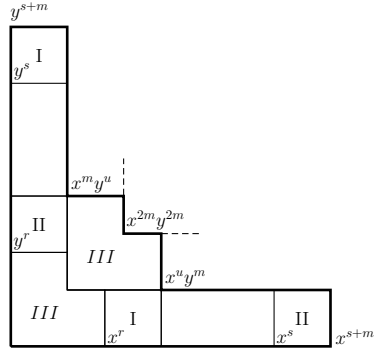
\begin{figure}[h]
\begin{center}
\begin{pspicture}(0,0)(5,4.25)
\scalebox{0.5}{
	\psline[linewidth=2pt](0,0)(8.5,0)(8.5,1.5)(4,1.5)(4,3)(3,3)(3,4)(1.5,4)(1.5,8.5)(0,8.5)(0,0)
	\psline(0,0)(4,0)(4,1.5)(1.5,1.5)(1.5,4)(0,4)(0,0)
	\psline[linestyle=dashed](5,3)(3,3)(3,5)
	\rput(0,0){\psline(0,0)(2.5,0)(2.5,1.5)(1.5,1.5)(1.5,2.5)(0,2.5)(0,0)\rput(0.9,0.9){\Large $III$}}
	\rput(1.5,1.5){\psline(0,0)(2.5,0)(2.5,1.5)(1.5,1.5)(1.5,2.5)(0,2.5)(0,0)\rput(0.9,0.9){\Large $III$}}
	
	\rput(2.5,0){\psline(0,0)(1.5,0)(1.5,1.5)(0,1.5)(0,0)\rput(0.75,0.75){\Large I}}
	\rput(0,7){\psline(0,0)(1.5,0)(1.5,1.5)(0,1.5)(0,0)\rput(0.75,0.75){\Large I}}
	\rput(0,2.5){\psline(0,0)(1.5,0)(1.5,1.5)(0,1.5)(0,0)\rput(0.75,0.75){\Large II}}
	\rput(7,0){\psline(0,0)(1.5,0)(1.5,1.5)(0,1.5)(0,0)\rput(0.75,0.75){\Large II}}

	\rput(3.85,3.35){\Large $x^{2m}y^{2m}$}	
	\rput(2.8,0.25){\Large $x^r$}\rput(0.3,2.8){\Large $y^r$}	
	\rput(7.3,0.25){\Large $x^s$}\rput(0.3,7.3){\Large $y^s$}	
	\rput(9.15,0.25){\Large $x^{s+m}$}\rput(0.5,8.9){\Large $y^{s+m}$}
	\rput(4.65,1.8){\Large $x^{u}y^m$}\rput(2.15,4.3){\Large $x^my^{u}$}
	}
\end{pspicture}
\caption{$G$-graph of type $B_{2}$.}
\label{GtypeB2}
\end{center}
\end{figure}

\begin{proof} Let $e_i=\frac{1}{2n}(r,s)$ and $e_{i+1}=\frac{1}{2n}(u,v)$ be two consecutive lattice points. First observe that 
\[
\text{$(x^{r-m}, y^{r-m})$, $(y^{s+m}, (-1)^{s+m}x^{s+m})$ and $(x^{r}y^m, -x^my^{r})$}
\]
belong to the irreducible representation $V_{2r-u}$. The pair $(x^{r-m}, y^{r-m})\in\Gamma$ since $x^{r-m}$ and $y^{r-m}$ are both elements in the $qG$-graph, thus we have two possibilities for the second element in $\Gamma$ belonging to the representation $V_{2r-u}$.

First assume that $(x^ry^m,x^my^r)\in\Gamma$, that is $x^{s+m},y^{s+m}\in I_{B_2}$. Then, by the assumption $x^{u}y^m,x^my^{u}\in I_{B_2}$ and the $G$-invariant relation $x^{2m}y^{2m}\in I_{B_2}$ we obtain the ideal $I_{B_2}=\langle x^{2m}y^{2m}, x^{s+m}, y^{s+m}, x^{u}y^{m}, x^{m}y^{u} \rangle$ as in the statement with the representation shown in Figure \ref{GtypeB2}. Note that the overlap in this case is extended without twin regions, so the number of elements that we add is
\[
\#I + \#II + \#III = 2m^2+r^2-(r-m)^2=u^2-r^2
\]
which is the number of elements in the overlap $O$. 

Let $\rho\in\Irr G$ be a representation of $G$ contained in $O$. If $\dim(\rho)=1$ then it belongs to the region $III$ and it is located in the diagonal $x=y$. As in the proof of \ref{TypeB1}, the representations located at the diagonal alternate the character of $\beta$ thus we have one element in each 1-dimensional representation. If $\dim(\rho)=2$ then in $\rho$ we have an element coming from $O$ and a second element coming from the extension, and we therefore have the required number of elements in each representation.  

Now assume that $(y^{s+m},x^{s+m})\in\Gamma$, that is $x^ry^m,x^my^r\in I_{B_2}$. Then the pairs 
\[
(x^{c-m},y^{c-m}), (x^{d+m},(-1)^{d+m}x^{d+m}), (x^cy^m,-x^my^c)
\]
are in the same representation $V_{c-m}$, where $c=2r-u$ and $d=2s-v$. By Lemma \ref{B1beforeB2}, the $qG$-graph defined by the previous consecutive pair $e_{i-1},e_i$ is of type $B$, which implies that $b_i=2$. Since $e_{i-1}+e_{i+1}=b_ie_i$ we have that $c$ and $d$ are the entries of the lattice point $e_{i-1}=\frac{1}{2n}(c,d)$. 

Following the same argument as before, the pair $(x^{c-m},y^{c-m})\in\Gamma$ since $x^{c-m}$ and $y^{c-m}$ are both elements in the $qG$-graph, and again we have two choices for the second basis element in $V_{c-m}$: either $(x^cy^m,-x^my^c)\in\Gamma$ or $(x^{d+m},(-1)^{d+m}x^{d+m})\in\Gamma$. 

If $(x^cy^m,-x^my^c)\in\Gamma$ then we get the ideal $I_{B_2}=\langle x^{2m}y^{2m},x^{d+m},y^{d+m},x^ry^m,x^my^r \rangle$, which is the ideal representing the $G$-graph $\Gamma_{B_2}(c,d;r,s)$ obtained with the pair $e_{i-1},e_i$. If $(x^{d+m},(-1)^{d+m}x^{d+m})\in\Gamma$ instead, we reach two new choices as before but now with the consecutive lattice points $e_{i-2}=\frac{1}{2n}(2c-r,2d-s)$ and $e_{i-1}$, and we continue the process in the same way. 

Therefore at every step we obtain an ideal of the desired form $\langle x^{2m}y^{2m}, x^{s_k+m}, y^{s_k+m}, x^{u_k}y^{m}, x^{m}y^{u_k} \rangle$ where $m=s_k-v_k=u_k-r_k$, corresponding to a pair $e_k=\frac{1}{2n}(r_k,s_k),e_{k-1}=\frac{1}{2n}(u_k,v_k)$ of consecutive lattice points defining a $qG$-graph of type $B_2$. The proof concludes by showing that this inductive process terminates in a finite number of steps. Indeed, by Lemma \ref{B1beforeB2} we eventually arrive to a pair $e_j,e_{j+1}$ for some $j<i$ defining a $qG$-graph of type $B_1$, for which the ideal is described in Proposition \ref{TypeB1}. 
\end{proof}

\noindent\underline{Case 2:} Suppose that $(x^uy^m,(-1)^mx^my^u)\in\Gamma$, that is $y^s,x^s\in I$. 
By considering the following three elements:
\[
(x^{u},y^{u}), (y^{v},(-1)^{v}x^{v}), (x^{2u-r}y^m,-x^my^{2u-r})\in V_u,
\]
we can see that $x^{2u-r}y^m,x^my^{2u-r}\in I_{B}$ because $x^{u},y^{u},y^{v}$ and $x^{v}$ are elements of the $qG$-graph $\widetilde{\Gamma}:=\widetilde{\Gamma}_B(r,s;u,v)$. Since the $qG$-graph $\widetilde{\Gamma}$ is of type $B$ we have that $b_i=2$, and we can write $t=2u-r$ and $w=2v-s$ which coincide with the entries of the next lattice point $e_{i+1}=\frac{1}{2n}(t,w)$. 
We obtain in this way the ideal $I_{B_2}=\langle x^{2m}y^{2m},x^{v+m},y^{v+m},x^ty^m,x^my^t \rangle$ with $m=v-w=t-u$ which represents the $G$-graph $\Gamma_{B_2}(u,v;t,w)$ of type $B_2$. Notice that even if $\widetilde{\Gamma}$ is of type $B_1$, by Lemma \ref{B1beforeB2} the $qG$-graph $\widetilde{\Gamma}_B(u,v;t,w)$ is of type $B_2$. As in the proof of Proposition \ref{TypeB2} the overlap $O=I\cup II\cup III$ in the $qG$-graph is extended without twin regions, and the total number of elements is the required. 

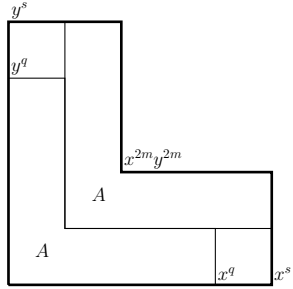
\begin{figure}[h]
\begin{center}
\begin{pspicture}(0,0.1)(3.5,3.75)
\scalebox{0.5}{
	\psline[linewidth=2pt](0,0)(7,0)(7,3)(3,3)(3,7)(0,7)(0,0)
	\psline(0,0)(7,0)(7,1.5)(1.5,1.5)(1.5,7)(0,7)(0,0)
	\psline[linestyle=dashed](5,3)(3,3)(3,5)
	\rput(0,0){\psline(0,0)(5.5,0)(5.5,1.5)(1.5,1.5)(1.5,5.5)(0,5.5)(0,0)
	\rput(0.9,0.9){\Large $A$}}
	\rput(2.4,2.4){\Large $A$}
	\rput(3.85,3.35){\Large $x^{2m}y^{2m}$}	
	\rput(5.8,0.3){\Large $x^q$}\rput(0.3,5.8){\Large $y^q$}	
	\rput(7.3,0.3){\Large $x^s$}\rput(0.3,7.3){\Large $y^s$}	
	}
\end{pspicture}
\caption{The last $G$-graph of type $B$.}
\label{GtypeBLast}
\end{center}
\end{figure}

When $\Gamma:=\Gamma_B(r,s;q,q)$ is the $G$-graph corresponding to the last $qG$-graph $\widetilde{\Gamma}_{h-1}$, we obtain the representation of $\Gamma$ shown in Figure \ref{GtypeBLast}, where the ideal in this case is $I=\langle x^{2m}y^{2m},x^s,y^s \rangle$. We do not treat this case because the corresponding open set $U_\Gamma$ covers the strict transform of the curve $E_h\subset\Hilb{H}{\C^2}$ with two fixed points by $G/H$, and every point in $U_\Gamma$ is already covered by the open sets corresponding to $G$-graphs at the blowup of the singular points in $\Hilb{H}{\C^2}/(G/H)$. These are called $G$-graphs of type $C$ and are explained in the next section.

\subsection{Remaining $G$-graphs: types C and D}\label{typesCyD}

For a given $qG$-graph $\widetilde{\Gamma}$, the corresponding $G$-graph $\Gamma$ was constructed in previous sections by adding some suitable elements to $\widetilde{\Gamma}$. This procedure gives almost all possible $G$-graphs. Indeed, the involution of the middle rational curve $E_h\subset\Hilb{H}{\C^2}$ to itself by $\beta$ gives two isolated fixed points, and therefore two more rational curves $E^+$ and $E^-$ in the exceptional locus for the resolution of $\mathbb{C}^2/G$ (see diagram below). This part of the exceptional locus cannot be recovered with $qG$-graphs. 

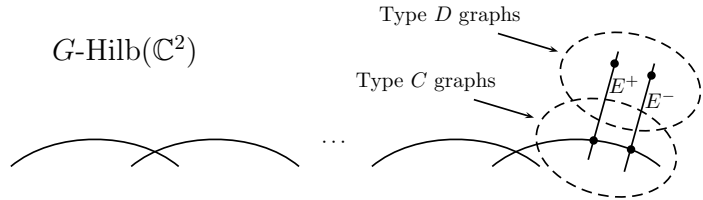
\begin{center}
\begin{pspicture}(0,-0.25)(9,2.5)
\scalebox{0.8}{
	\psset{nodesep=3pt,arcangle=40,dotsize=4pt}
	\rput(0,0){\rnode{P1}{}}\rput(3,0){\rnode{P2}{}}	
	\rput(2,0){\rnode{P3}{}}\rput(5,0){\rnode{P4}{}}
	\rput(6,0){\rnode{P5}{}}\rput(9,0){\rnode{P6}{}}
	\rput(8,0){\rnode{P7}{}}\rput(11,0){\rnode{P8}{}}
	\ncarc{-}{P1}{P2}\ncarc{-}{P3}{P4}
	\rput(5.5,0.5){$\cdots$}
	\ncarc{-}{P5}{P6}\ncarc{-}{P7}{P8}
	\psline(9.7,0.2)(10.2,2)\psline(10.3,0)(10.8,1.8)
	\psdots(9.79,0.52)(10.41,0.37)(10.14,1.8)(10.75,1.6)
	\rput[b]{-15}(10,0.4){\psellipse[linestyle=dashed](0,0)(1.2,0.8)}
	\rput[b]{-15}(10.4,1.5){\psellipse[linestyle=dashed](0,0)(1.2,0.8)}
	\rput(7.4,2.6){Type $D$ graphs}\rput(8.2,2.3){\psline{->}(0,0)(1,-0.3)}
	\rput(7,1.5){Type $C$ graphs}\rput(7.8,1.2){\psline{->}(0,0)(1,-0.3)}
	\rput(10.3,1.45){$E^+$}
	\rput(10.9,1.2){$E^-$}
	\rput(2,2){\Large $\Hilb{G}{\C^2}$}
	}
\end{pspicture}
\end{center}

In this section we construct the $G$-graphs which correspond to the neighbourhood of $E^+$ and $E^-$. They are called {\em type C} and {\em type D}, and each of them have two cases, $C^+$ and $C^-$ ($D^+$ and $D^-$ respectively). As mentioned in Section \ref{Sect:qGraphs}, these $G$-graphs arise from the $\beta$-fixed points in $\Hilb{H}{\C^2}$. Since these fixed points are contained in the middle rational curve $E_h$, the new $G$-graphs must contain the $qG$-graph formed by union of the two middle $H$-graphs $\Gamma_{h-1}$ and $\Gamma_h$. Thus to construct $G$-graphs of types $C$ and $D$ we start from the data given by $\widetilde{\Gamma}_{h-1}(r,s;q,q)$.

Let $e_h=\frac{1}{2n}(r,s)$ and $e_{h+1}=\frac{1}{2n}(q,q)$ be the lattice points giving the last $qG$-graph $\widetilde{\Gamma}_{h-1}(r,s;q,q)$, and define 
\[ m_{1}:=s-q ~\text{ and }~ m_{2}:=q-r. \]
Note that if the last $qG$-graph is of type $B$ then $m_{1}=m_{2}$, and if it is of type $A$ then $m_{1}>m_{2}$. 

\begin{rem}\label{GInvC&D}The monomial $x^{s-r}y^{s-r}$ is $G$-invariant, and therefore is in $I$. Indeed, $\alpha(x^{s-r}y^{s-r})= \beta(x^{s-r}y^{s-r})=x^{s-r}y^{s-r}$ since $r$ and $s$ have the same parity, and $ar\equiv s$ and $as\equiv r$ (mod $2n$). Similarly, the pairs of monomials $x^{s-r+jq}y^{s-r-jq}+(-1)^{jq}x^{s-r-jq}y^{s-r+jq}$ are also $G$-invariant for every $j$. 
\end{rem} 

The new $G$-graphs depend on a choice of basis between the following polynomials in the irreducible representations $\rho_q^+$ and $\rho_q^-$:
\begin{align*}
x^q+(-i)^qy^q,~~x^sy^{m_{2}}+(-1)^ri^qx^{m_{2}}y^s &\in\rho_{q}^+ \\
x^q-(-i)^qy^q,~~x^sy^{m_{2}}-(-1)^ri^qx^{m_{2}}y^s & \in\rho_{q}^- 
\end{align*}
The choice of basis from these polynomials is sufficient since they generate the $S^G$-modules $S_{\rho_q^\pm}$. In fact, $\rho^\pm_q$ are the {\em special} irreducible representations corresponding to the curves $E^+$ and $E^-$ (compare with \cite[Theorems 6.2 and 10.1]{IW08}, together with Section \ref{SpecialRepr}). 

Note that at least one of the polynomials $x^q+(-i)^qy^q$ and $x^q-(-i)^qy^q$ must be in the basis $\Gamma$ of $\C[x,y]/I$. Indeed, if $x^q+(-i)^qy^q, x^q-(-i)^qy^q\in I$ then also $x^q,y^q\in I$, which is a contradiction since they are elements of $\widetilde{\Gamma}_{h-1}$. Thus we have three possibilities: $x^q+(-i)^qy^q\in I$, $x^q-(-i)^qy^q\in I$ and both $x^q\pm(-i)^qy^q\in\Gamma$.

\subsubsection*{$G$-graphs of type $D$}

Suppose that $x^q+(-i)^qy^q\in I$. Then $x^sy^{m_{2}}+(-1)^ri^qx^{m_{2}}y^s$ must be in the $G$-graph $\Gamma$, and the basis polynomials in $\Gamma$ are $x^q-(-i)^qy^q$ and $x^sy^{m_{2}}+(-1)^ri^qx^{m_{2}}y^s$. Similarly, we can choose $x^q-(-i)^qy^q\in I$, which now implies $x^q+(-i)^qy^q, x^sy^{m_{2}}-(-1)^ri^qx^{m_{2}}y^s\in\Gamma$. The first case corresponds to type $D^+$ and the second to type $D^-$. 

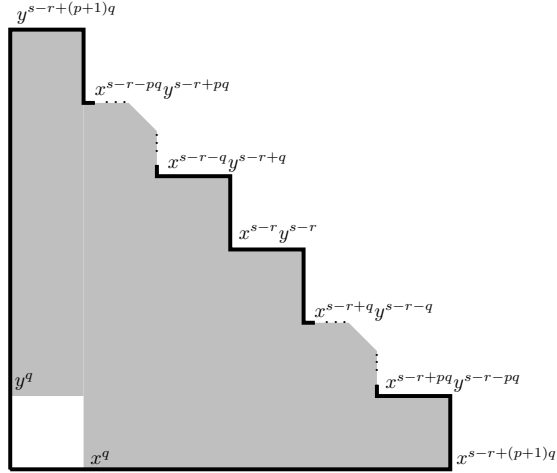
\begin{figure}[h]
\begin{center}
\begin{pspicture}(0,0)(6,6)
	\scalebox{0.75}{
	\psline[linestyle=none,fillstyle=solid,fillcolor=lightgray](1.3,1.3)(1.3,6.5)(2.1,6.5)(2.6,6)(2.6,5.2)(3.9,5.2)(3.9,3.9)(5.2,3.9)(5.2,2.6)(6,2.6)(6.5,2.1)(6.5,1.3)(1.3,1.3)
	\psline[linestyle=none,fillstyle=solid,fillcolor=lightgray](0,1.3)(1.3,1.3)(1.3,7.8)(0,7.8)(0,0)
	\psline[linestyle=none,fillstyle=solid,fillcolor=lightgray](1.3,0)(1.3,1.3)(7.8,1.3)(7.8,0)(0,0)	
	\psline(6.5,1.5)(6.5,1.3)(7.8,1.3)(7.8,0)(0,0)(0,7.8)(1.3,7.8)(1.3,6.5)(1.5,6.5)
	\psline[linewidth=2pt](2.6,5.4)(2.6,5.2)(3.9,5.2)(3.9,3.9)(5.2,3.9)(5.2,2.6)(5.4,2.6)
	\rput(6.5,2){$\vdots$}\rput(5.8,2.6){$\cdots$}
	\rput(1.9,6.5){$\cdots$}\rput(2.6,5.9){$\vdots$}
	\rput(1.6,0.2){$x^q$}\rput(0.25,1.55){$y^q$}
	\rput(8.8,0.25){$x^{s-r+(p+1)q}$}\rput(1,8.1){$y^{s-r+(p+1)q}$}
	\rput(7.85,1.55){$x^{s-r+pq}y^{s-r-pq}$}\rput(2.7,6.75){$x^{s-r-pq}y^{s-r+pq}$}
	\rput(6.45,2.8){$x^{s-r+q}y^{s-r-q}$}\rput(3.85,5.45){$x^{s-r-q}y^{s-r+q}$}
	\rput(4.75,4.2){$x^{s-r}y^{s-r}$}
	\psline[linewidth=2pt](0,0)(7.8,0)(7.8,1.3)(6.5,1.3)(6.5,1.5)
	\psline[linewidth=2pt](1.5,6.5)(1.3,6.5)(1.3,7.8)(0,7.8)(0,0)
	}
\end{pspicture}
\caption{$G$-graph of type $D$. The monomials $x^q$ and $y^q$ are identified.}
\label{GtypeD}
\end{center}
\end{figure}

The assumption $x^q+(-i)^qy^q\in I$ (or analogously, $x^q-(-i)^qy^q\in I$ for the Case $D^-$) identifies the monomials $x^q$ with $y^q$ as twin elements in $\mathbb{C}[x,y]/I$, and together with $x^{s-r}y^{s-r}\in I$ characterise completely the shape of the $G$-graph $\Gamma$ (see Figure \ref{GtypeD}). This gives the following proposition.    

\begin{prop}\label{TypeD} The ideals $I_{D^\pm}=\langle x^q\pm(-i)^qy^q, x^{s-r}y^{s-r}\rangle$ define $G$-clusters. The $G$-graphs represented by $I_{D^{\pm}}$ are said to be of type $D^{\pm}$.
\end{prop}

\begin{proof} We give the proof for the case $D^-$, the case $D^+$ is almost identical. We begin by identifying a collection of monomials $Q_\Z$ that will prove useful. Recall that $M$ denotes the lattice of Laurent monomials in the variables $x$ and $y$. Let $P\subset M\otimes\R$ denote the parallelogram whose vertices include the monomials $(0,0)$, $(\frac{s-r}{2},\frac{s-r}{2})$ and $(q,0)$, and let $P^*$ denote its reflection through the line $(x=y)\subset M\otimes\R$. Define $Q_\Z$ to be the set of monomials in the interior of $P\cup P^*$ together with the monomials
\[
1,x,x^2,\ldots,x^q,y,y^2,\ldots,y^{q-1},x^q(xy),\ldots,x^q(xy)^{\frac{s-r}{2}-1}
\]
that is, $(P\cup P^*)\cap M$ with the appropriate boundary monomials removed.

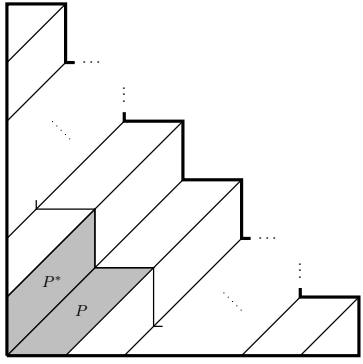
\begin{figure}[h]
\begin{center}
\begin{pspicture}(0,0)(4.5,5)
\scalebox{0.6}{
	\psline[linestyle=none](1.3,1.3)(1.3,6.5)(2.1,6.5)(2.6,6)(2.6,5.2)(3.9,5.2)(3.9,3.9)(5.2,3.9)(5.2,2.6)(6,2.6)(6.5,2.1)(6.5,1.3)(1.3,1.3)
	\psline[linestyle=none](0,1.3)(1.3,1.3)(1.3,7.8)(0,7.8)(0,0)
	\psline[linestyle=none](1.3,0)(1.3,1.3)(7.8,1.3)(7.8,0)(0,0)	
	\psline(6.5,1.5)(6.5,1.3)(7.8,1.3)(7.8,0)(0,0)(0,7.8)(1.3,7.8)(1.3,6.5)(1.5,6.5)
	\psline[linewidth=2pt](2.6,5.4)(2.6,5.2)(3.9,5.2)(3.9,3.9)(5.2,3.9)(5.2,2.6)(5.4,2.6)
	\rput(6.5,2){$\vdots$}\rput(5.8,2.6){$\cdots$}
	\rput(1.9,6.5){$\cdots$}\rput(2.6,5.9){$\vdots$}
	\psline[fillstyle=solid,fillcolor=lightgray](0,0)(1.95,1.95)(3.25,1.95)(1.3,0)
	\psline[fillstyle=solid,fillcolor=lightgray](0,0)(1.95,1.95)(1.95,3.25)(0,1.3)
	\rput(1.65,1){$P$}\rput(1,1.65){$P^*$}
	\rput(1.95,1.95){\psline(0,0)(1.95,1.95)(3.25,1.95)(1.3,0)}
	\psline(0,1.3)(3.9,5.2)
	\psline(0,2.6)(2.6,5.2)
	\psline(0,5.2)(1.3,6.5)
	\psline(0,6.5)(1.3,7.8)
	\psline(2.6,0)(5.2,2.6)
	\psline(5.2,0)(6.5,1.3)
	\psline(6.5,0)(7.8,1.3)
	\psline(0.65,3.45)(0.65,3.25)(1.95,3.25)(1.95,1.95)(3.25,1.95)(3.25,0.65)(3.45,0.65)
	\psline[linestyle=dotted](5.2,1)(4.8,1.4)
	\psline[linestyle=dotted](1,5.2)(1.4,4.8)
	\psline[linewidth=2pt](0,0)(7.8,0)(7.8,1.3)(6.5,1.3)(6.5,1.5)
	\psline[linewidth=2pt](1.5,6.5)(1.3,6.5)(1.3,7.8)(0,7.8)(0,0)
	}
\end{pspicture}
\caption{Parallelograms $P$ and $P^*$ tessellating $M\otimes\R$.}
\label{GtypeD-proof}
\end{center}
\end{figure}

Note that translating $Q_\Z$ by integer combinations of the vectors $(\frac{s-r}{2},\frac{s-r}{2})$ and $(q,-q)$ tessellates $M$. Also note that the Laurent monomials $(xy)^{\frac{s-r}{2}}$ and $x^qy^{-q}$ are $H$-invariant. So for every $\rho\in\Irr H$ there exists a monomial in $Q_\Z$ that belongs to $\rho$. Since $(r,s)$ and $(q,q)$ define an $H$-cluster, $q(r-s)=$ Area$(P\cup P^*)=2n$. Therefore every $\rho\in\Irr H$ contains exactly one monomial ${\bf x}^\rho$ from $Q_\Z$.

Now we use $Q_\Z$ and $(\frac{s-r}{2},\frac{s-r}{2})+Q_\Z$ to give a polynomial in every irreducible representation of $G$. Given $\rho_j\in\Irr H$ for $j\not\equiv aj$ (mod $2n$), take the pairs

\[
({\bf x}^{\rho_j},\beta({\bf x}^{\rho_j})), ((xy)^{\frac{s-r}{2}}{\bf x}^{\rho_j},(xy)^{\frac{s-r}{2}}\beta({\bf x}^{\rho_j}))\in V_j.
\]
Take $\rho_j$ such that $j\equiv aj$ (mod $2n$). Since $q:=\frac{2n}{(a-1,2n)}$, the monomials ${\bf x}^{\rho_j}$ are either of the form $x^q(xy)^l$ or $(xy)^l$, for some $j\in\N$. If ${\bf x}^{\rho_j}=(xy)^l$ take
\begin{align*}
(xy)^l &\in \rho_j^{(-1)^l} \\
(xy)^{l+\frac{s-r}{2}} &\in \rho_j^{(-1)^{l+1}}. 
\end{align*}
If ${\bf x}^{\rho_j}=x^q(xy)^l$ and $n,q$ are even, then take
\begin{align*}
x^q(xy)^l+(-1)^l\beta(x^q(xy)^l) &\in \rho_j^{(-1)^l} \\
x^q(xy)^{l+\frac{s-r}{2}}+(-1)^{l+1}\beta(x^q(xy)^{l+\frac{s-r}{2}}) &\in \rho_j^{(-1)^{l+1}}. 
\end{align*}
If $n$ and $q$ are odd, multiply the second summand by $i$ as in Equation \ref{eqn:polyrep}. Here we are using the fact that we are in the case $D^-$ and the assumption that the group is small, so $\frac{s-r}{2}=k$ is odd.

Observing that the monomials in $Q_\Z\cup((\frac{s-r}{2},\frac{s-r}{2})+Q_\Z)$ form a basis of the vector space $\C[x,y]/I_{D^+}$ completes the proof. 
\end{proof}

\subsubsection*{$G$-graphs of type $C$}

Now assume $x^q+(-i)^qy^q$ and $x^q-(-i)^qy^q$ are basis elements for $\rho^+_{q}$ and $\rho^-_{q}$, i.e.\  $x^sy^{m_{2}}\pm(-1)^ri^qx^{m_{2}}y^s\in I$, which implies that $x^sy^{m_{2}}, x^{m_{2}}y^s\in I$. 

The pairs $(x^{r}, y^{r})$, $(y^{s}, (-1)^{s}x^{s})$, $(x^{q}y^{m_{1}}, (-1)^{m_{1}}x^{m_{1}}y^{q})\in V_{r}$ (omitting the corresponding signs), and notice that all monomials in them must belong to the basis (otherwise the $G$-graph would have less than $|G|$ elements). Therefore there must be an identification between pairs of the same degree. We take as basic elements $(x^{r},y^{r})$ and a linear combination of $(y^{s},(-1)^{s}x^{s})$ and $(x^qy^{m_{1}},(-1)^{m_{1}}x^{m_{1}}y^q)$. We have two possibilities:
\[
\begin{array}{cl}
\textbf{Case $C^+$}:&  y^{m_{1}}(x^q+(-i)^qy^{q}), x^{m_{1}}(x^{q}+(-i)^qy^q)\in I \text{ and } \\
	& (x^{r}, y^{r}), (y^{m_{1}}(x^q-(-i)^qy^{q}), x^{m_{1}}(x^{q}-(-i)^qy^q))\in\Gamma, \\
\textbf{Case $C^-$}:& y^{m_{1}}(x^q-(-i)^qy^{q}), x^{m_{1}}(x^{q}-(-i)^qy^q)\in I \text{ and } \\
	& (x^{r}, y^{r}), (y^{m_{1}}(x^q+(-i)^qy^{q}), x^{m_{1}}(x^{q}+(-i)^qy^q))\in\Gamma.
\end{array}
\]
Note also that $x^{2q}+(-1)^qy^{2q},x^qy^q\in\rho_{2q}^+$, while we need only one element in $\Gamma$ corresponding to $\rho_{2q}^+$. On the other hand both of them must belong to $\Gamma$, otherwise $\Gamma$ does not contain $|G|$ elements. This implies that we have to take a combination of them as our basis for $\rho_{2q}^+$. We take 
\begin{align*}
(x^q+(-i)^qy^q)^2&\in I~~\text{for type $C^+$,  and} \\
(x^q-(-i)^qy^q)^2&\in I~~\text{for type $C^-$}  
\end{align*}

In both cases, if the last $qG$-graph $\Gamma_{h-1}$ is of type $B$ then this last equation is redundant, so it is only needed in the type $A$ case. Indeed, substituting the value of $s=2q-r$ into the equations for case $C^+$ we get that $(-i)^qy^{2q-r}+x^qy^{q-r}, x^{2q-r}+(-i)^qx^{q-r}y^q\in I$. Now multiplying the first polynomial by $(-i)^qy^r$ and the second polynomial by $x^r$, and adding them together we get $x^{2q}+(-1)^qy^{2q}+2(-i)^qx^qy^q$ as desired. 

When the last $G$-graph $\Gamma_{h-1}$ is of type $A$, we do need the equation $x^{2q}+(-1)^qy^{2q}\pm2(-i)^qx^qy^q$ together with the $G$-invariant $x^{m_{1}}y^{m_{2}}+(-1)^{m_{2}}x^{m_{2}}y^{m_{1}}$. In this case $x^{s-r}y^{s-r}$ can be obtained using the remaining identities. In other words, $x^{s-r}y^{s-r}\in\langle x^{2q}+(-1)^qy^{2q}+2(-i)^qx^qy^q, x^sy^{m_{2}}-(-i)^qx^{m_{2}}y^s,x^{m_{1}}y^{m_{2}}+(-1)^{m_{2}}x^{m_{2}}y^{m_{1}}\rangle$. Indeed, using the second and third generators of the ideal we have that 
\begin{align*}
x^{s-r-q}y^{s-r+q}  & = (-i)^{q}x^{2(s-q)}y^{2(q-r)} =-(-1)^{q-r}(-i)^qx^{r-s}y^{r-s}, \text{ and} \\
x^{s-r+q}y^{s-r-q}  &= (-i)^{q}x^{2(q-r)}y^{2(s-q)}  =-(-1)^{q-r}(-i)^qx^{r-s}y^{r-s}
\end{align*}
Multiplying the first generator by $x^{s-r-q}y^{s-r-q}$ getting the polynomial $-(-1)^{q-r}x^{s-r}y^{s-r}-(-1)^rx^{s-r}y^{s-r}+2x^{s-r}y^{s-r}$ (note that $s<2(s-q)$). The only possibility for this to be identically is zero only if $q$ and $r$ are both even at the same time. But this is impossible because it would imply that every boundary lattice point in the lattice $L$ is even, a contradiction. Therefore, the sum above is not identically zero, which implies that $x^{s-r}y^{s-r}\in I$, and we obtain the following proposition.

\begin{prop}\label{TypeC} The ideals 
\begin{align*}
I_{C^\pm_{A}}&=\langle (x^{q}\pm(-i)^qy^{q})^2, x^sy^{m_{2}}\pm(-1)^ri^qx^{m_{2}}y^s, x^{m_{1}}y^{m_{2}}+(-1)^{m_{2}}x^{m_{2}}y^{m_{1}}\rangle \text{ if $\Gamma_{h-1}$ is of type $A$}, or \\
I_{C^\pm_{B}}&=\langle y^{m}(x^q\pm(-i)^qy^{q}), x^{m}(x^q\pm(-i)^qy^q),x^{s-r}y^{s-r},x^sy^{m}, x^{m}y^s\rangle \text{ if $\Gamma_{h-1}$ is of type $B$}. 
\end{align*}
define $G$-clusters. The $G$-graphs represented by $I_{C^{\pm}}$ are said to be of type $C^{\pm}$.
\end{prop}

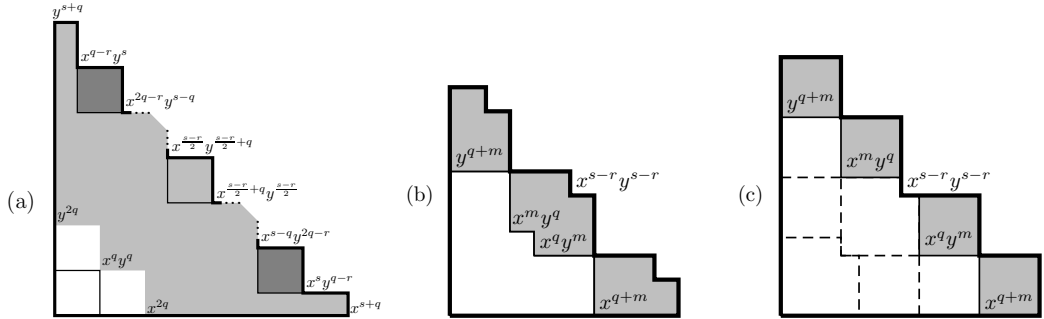
\begin{figure}[h]
\begin{center}
\begin{pspicture}(0,0)(12,4)
\scalebox{0.8}{
\rput(-1,-0.3725){
	\scalebox{0.75}{
	\psline[linestyle=none,fillstyle=solid,fillcolor=lightgray](2.5,0.5)(7,0.5)(7,1)(6,1)(6,2)(5,2)(5,2.6)(4.6,3)(4,3)(4,4)(3,4)(3,4.6)(2.6,5)(2,5)(2,6)(1,6)(1,7)(0.5,7)(0.5,2.5)(1.5,2.5)(1.5,1.5)(2.5,1.5)(2.5,0.5)
	\psline[fillstyle=solid,fillcolor=gray](1,6)(1,5)(2,5)(2,6)(1,6)
	\psline(3,4)(3,3)(4,3)
	\psline[fillstyle=solid,fillcolor=gray](6,1)(5,1)(5,2)(6,2)(6,1)
	\psline(0.5,1.5)(1.5,1.5)(1.5,0.5)
	\psline[linewidth=2pt](3,4.2)(3,4)(4,4)(4,3)(4.2,3)
	\psline[linewidth=2pt](5,2.2)(5,2)(6,2)(6,1)(7,1)(7,0.5)(0.5,0.5)(0.5,7)(1,7)(1,6)(2,6)(2,5)(2.2,5)
	\rput(5,2.55){\Large $\vdots$}\rput(4.45,3){\Large $...$}
	\rput(2.45,5){\Large $...$}\rput(3,4.55){\Large $\vdots$}
	\rput(1.9,1.7){$x^qy^q$}
	\rput(2.8,0.7){$x^{2q}$}\rput(0.8,2.75){$y^{2q}$}
	\rput(7.4,0.7){$x^{s+q}$}\rput(0.8,7.25){$y^{s+q}$}
	\rput(5,3.25){$x^{\frac{s-r}{2}+q}y^{\frac{s-r}{2}}$}\rput(4,4.3){$x^{\frac{s-r}{2}}y^{\frac{s-r}{2}+q}$}
	\rput(6.6,1.25){$x^{s}y^{q-r}$}\rput(5.85,2.25){$x^{s-q}y^{2q-r}$}
	\rput(2.85,5.25){$x^{2q-r}y^{s-q}$}\rput(1.6,6.25){$x^{q-r}y^{s}$}
	\rput(-0.25,3){\Large (a)}
	}}
\rput(6,0){
	\psline(0,0)(2.4,0)(2.4,1)(1.4,1)(1.4,1.4)(1,1.4)(1,2.4)(0,2.4)(0,0)
	\psline[fillstyle=solid,fillcolor=lightgray](1.4,1)(2.4,1)(2.4,2)(2,2)(2,2.4)(1,2.4)(1,1.4)(1.4,1.4)(1.4,1)
	\rput(1.9,1.2){$x^qy^m$}\rput(1.45,1.6){$x^my^q$}
	\rput(2.4,0){\psline[fillstyle=solid,fillcolor=lightgray](0,0)(1.4,0)(1.4,0.6)(1,0.6)(1,1)(0,1)(0,0)
	\rput(0.5,0.25){$x^{q+m}$}}
	\rput(0,2.4){\psline[fillstyle=solid,fillcolor=lightgray](0,0)(1,0)(1,1)(0.6,1)(0.6,1.4)(0,1.4)(0,0)
	\rput(0.5,0.25){$y^{q+m}$}}
	\rput(2.8,2.25){$x^{s-r}y^{s-r}$}
	\rput(-0.5,2){(b)}
	\psline[linewidth=2pt](0,0)(3.8,0)(3.8,0.6)(3.4,0.6)(3.4,1)(2.4,1)(2.4,2)(2,2)(2,2.4)(1,2.4)(1,3.4)(0.6,3.4)(0.6,3.8)(0,3.8)(0,0)	
	}
\rput(11.5,0){
	\psline(0,0)(4.3,0)(4.3,1)(3.3,1)(3.3,2)(2,2)(2,3.3)(1,3.3)(1,4.3)(0,4.3)(0,0)
	\rput(0,0){\psline[linestyle=dashed](0,0)(1.3,0)(1.3,1)(1,1)(1,1.3)(0,1.3)(0,0)}
	\rput(1,1){\psline[linestyle=dashed](0,0)(1.3,0)(1.3,1)(1,1)(1,1.3)(0,1.3)(0,0)}
	\psline[linestyle=dashed](2.3,0)(2.3,1)
	\psline[linestyle=dashed](0,2.3)(1,2.3)
	\rput(3.3,0){\psline[fillstyle=solid,fillcolor=lightgray](0,0)(1,0)(1,1)(0,1)(0,0)
	\rput(0.5,0.25){$x^{q+m}$}}
	\rput(2.3,1){\psline[fillstyle=solid,fillcolor=lightgray](0,0)(1,0)(1,1)(0,1)(0,0)
	\rput(0.5,0.25){$x^{q}y^m$}}
	\rput(1,2.3){\psline[fillstyle=solid,fillcolor=lightgray](0,0)(1,0)(1,1)(0,1)(0,0)
	\rput(0.5,0.25){$x^{m}y^q$}}
	\rput(0,3.3){\psline[fillstyle=solid,fillcolor=lightgray](0,0)(1,0)(1,1)(0,1)(0,0)
	\rput(0.5,0.25){$y^{q+m}$}}
	\rput(2.8,2.25){$x^{s-r}y^{s-r}$}
	\rput(-0.5,2){(c)}
	\psline[linewidth=2pt](0,0)(4.3,0)(4.3,1)(3.3,1)(3.3,2)(2,2)(2,3.3)(1,3.3)(1,4.3)(0,4.3)(0,0)	
	}}
\end{pspicture}
\caption{$G$-graph of type $C$ when $\Gamma_{h-1}$ is (a) of type $A$, (b) of type $B_1$, and (c) of type $B_2$.}
\label{GtypeC}
\end{center}
\end{figure}

The shape of the $G$-graph is the same in both cases, and as in the type $D$ case, it has a stair shape. In this case, the conditions $x^sy^{m_{2}},x^{m_{2}}y^s\in I$ make the stair smaller than in type $D$. See Figure \ref{GtypeC}. 

\begin{proof} We treat the case when $\Gamma_{h-1}$ is of type $A$, the type $B$ case is similar. Consider the parallelograms $P$ and $P^*$, and the set of monomials $Q_\Z$ tessellating $M\otimes\R$ as in the proof of Proposition \ref{TypeD}. Denote by $L_0:(x=y)\subset M\otimes\R$ the diagonal, and by $L_c$ the parallel line to $L_0$ passing through the point $(cq,0)$ where $c\in\Z$. Let also $v_c\in\Aut(M\otimes\R)$ be the refection through the line $L_c$. Note that $P^*=v_0(P)$.

Let $C_0$ be the square in $M\otimes\R$ of vertices $(0,0)$, $(q-1,0)$, $(q-1,q-1)$ and $(0,q-1)$, and define $C$ to be the square $(\frac{s-r}{2},\frac{s-r}{2})+C_0$. Note that reflecting the polygon $P\cup P^*\cup C$ along the lines $L_c$ and taking only monomials in the positive quadrant we obtain the representation of $\Gamma$ (see Figure \ref{GtypeC-Proof}). 

Consider the region $\mathcal{R}:=(P\cup P^*\cup v_1(P)\cup v_{-1}(P^*))\cap\Box$, where $\Box$ denotes the positive quadrant, and define $\widetilde{Q}_\Z$ to be the set of monomials in the interior of $\mathcal{R}$ together with every monomial in $C$ and
\[
1,x,\ldots,x^{2q},y,\ldots,y^{2q-1},x^{2q}(xy),\ldots,x^{2q}(xy)^{\frac{s-r}{2}-q-1}.
\]
As in the proof of Proposition \ref{TypeD} we have that $\widetilde{Q}_\Z$ tessellates $M$ by translating it this time with the vectors $(2q,-2q)$ and $(\frac{s-r}{2},\frac{s-r}{2})$. Notice also that the monomials $x^{2q}y^{-2q}$ and $x^{\frac{s-r}{2}}y^{\frac{s-r}{2}}$ are $H$-invariant. Moreover, since $Q_\Z$ contains one monomial ${\bf x}^\rho$ in each $\rho\in\Irr H$ and $Q_\Z$ tessellates $M$, then $\widetilde{Q}_\Z$ contains exactly two monomials ${\bf x}^\rho$ and ${\bf y}^\rho$ in each $\rho\in\Irr H$.

\begin{figure}[h]
\begin{center}
\begin{pspicture}(0,1.25)(6,5.5)
\scalebox{0.7}{
	\psline[fillstyle=solid,fillcolor=lightgray](1.5,1.5)(3.5,1.5)(5,3)(5,5)(3,5)(1.5,3.5)(1.5,1.5)
	\psline[linewidth=2pt](6.8,3)(7,3)(7,2)(8,2)(8,1.5)(1.5,1.5)(1.5,8)(2,8)(2,7)(3,7)(3,6.8)

	\psline(1.5,1.5)(4,4)(5,4)(2.5,1.5)\rput(3.5,3){$P$}
	\psline(1.5,1.5)(4,4)(4,5)(1.5,2.5)\rput(3,3.5){$P^*$}
	\psline(2.5,1.5)(5,4)(5,3)(3.5,1.5)(2.5,1.5)\rput(3.75,2.25){\small $v_1(P)$}
	\psline(1.5,2.5)(4,5)(3,5)(1.5,3.5)(1.5,2.5)\rput(2.25,3.75){\small $v_{\text{-}\!1}\!(\!P^*\!)$}
	
	\psline[linestyle=dashed](5,5)(6,6)\rput(6.25,6.25){$L_0$}
	\psline[linestyle=dashed](5,4)(6,5)\rput(6.25,5.25){$L_1$}
	\psline[linestyle=dashed](4,5)(5,6)\rput(5.25,6.25){$L_{-1}$}

	\psline(7,1.5)(7,2)(6,2)(6,2.2)
	\psline(1.5,7)(2,7)(2,6)(2.2,6)
	\rput(4.5,4.5){$C$}

	\psline(5,3)(5.2,3)
	\psline(3,5)(3,5.2)
	\psline(5.5,1.5)(6,2)	
	\psline(1.5,5.5)(2,6)	
	\psline(6.5,1.5)(7,2)	
	\psline(1.5,6.5)(2,7)	
	
	\psline[linewidth=2pt](3.8,6)(4,6)(4,5)(5,5)(5,4)(6,4)(6,3.8)
	\psline[linestyle=dotted](2.5,5)(2,5.5)
	\psline[linestyle=dotted](5,2.5)(5.5,2)
	\psline[linestyle=dotted](3.25,5.75)(2.75,6.25)
	\psline[linestyle=dotted](6.25,2.75)(5.75,3.25)

}
\end{pspicture}
\caption{The region $\mathcal{R}\cup C$ containing the set of monomials $\widetilde{Q}_\Z$.}
\label{GtypeC-Proof}
\end{center}
\end{figure}
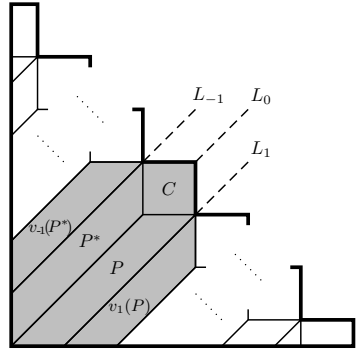

Now using $\widetilde{Q}_\Z$ we produce the correct number of polynomials in $\Irr G$. We consider the case $C^-$, the case $C^+$ is almost identical. Let $\rho_j\in\Irr H$ and suppose that $j\equiv aj$ (mod $2n$). Then we take 
\[
({\bf x}^{\rho_j},\beta({\bf x}^{\rho_j})), ({\bf y}^{\rho_j},\beta({\bf y}^{\rho_j}))\in V_j
\]
Notice that if ${\bf x}^{\rho_j}\in C_0$ then ${\bf y}^{\rho_j}=(xy)^{\frac{s-r}{2}}{\bf x}^{\rho_j}$, and if ${\bf x}^{\rho_j}\notin C_0\cup C$ then either ${\bf x}^{\rho_j}\in P$ and ${\bf y}^{\rho_j}\in v_{-1}(P^*)$, or ${\bf x}^{\rho_j}\in P^*$ and ${\bf y}^{\rho_j}\in v_{1}(P)$.

If $j\not\equiv aj$ (mod $2n$) then ${\bf x}^{\rho_j}$ is either of the form $(xy)^l$, $x^q(xy)^l$ or $x^{2q}(xy)^l$ for some $l\in\N$. If ${\bf x}^{\rho_j}\in C_0$ then ${\bf x}^{\rho_j}=(xy)^l$ so take 
\begin{align*}
(xy)^l &\in \rho_j^{(-1)^l} \\
(xy)^{l+\frac{s-r}{2}} &\in \rho_j^{(-1)^{+1}l} 
\end{align*}
If ${\bf x}^{\rho_j}\notin C_0$ then we have two cases:
\begin{itemize}
	\item[(i)] ${\bf x}^{\rho_j}=(xy)^l$ and ${\bf y}^{\rho_j}=x^{2q}(xy)^{l-q}$. Then take
		\begin{align*}
		(xy)^l &\in \rho_j^{(-1)^l} \\
		(xy)^{l-q}(x^{2q}-(-1)^q\beta(x^{2q})) &\in \rho_j^{(-1)^{+1}l}
		\end{align*}
	\item[(ii)] ${\bf x}^{\rho_j}=x^q(xy)^l$ and ${\bf y}^{\rho_j}=\beta({\bf x}^{\rho_j})$, in which case we take the same elements as in the case $D$, which are 
\begin{align*}
x^q(xy)^l+(-1)^l\beta(x^q(xy)^l) &\in \rho_j^{(-1)^l} \\
x^q(xy)^{l+\frac{s-r}{2}}+(-1)^{l+1}\beta(x^q(xy)^{l+\frac{s-r}{2}}) &\in \rho_j^{(-1)^{l+1}} 
\end{align*}
if  $n,q$ are even, and multiply the second summand by $i$ as in Equation \ref{eqn:polyrep} if $n$ and $q$ are odd.
\end{itemize}
Observing that the monomials in $\widetilde{Q}_\Z$ form a basis of the vector space $\C[x,y]/I_{C^-}$ completes the proof.
\end{proof}

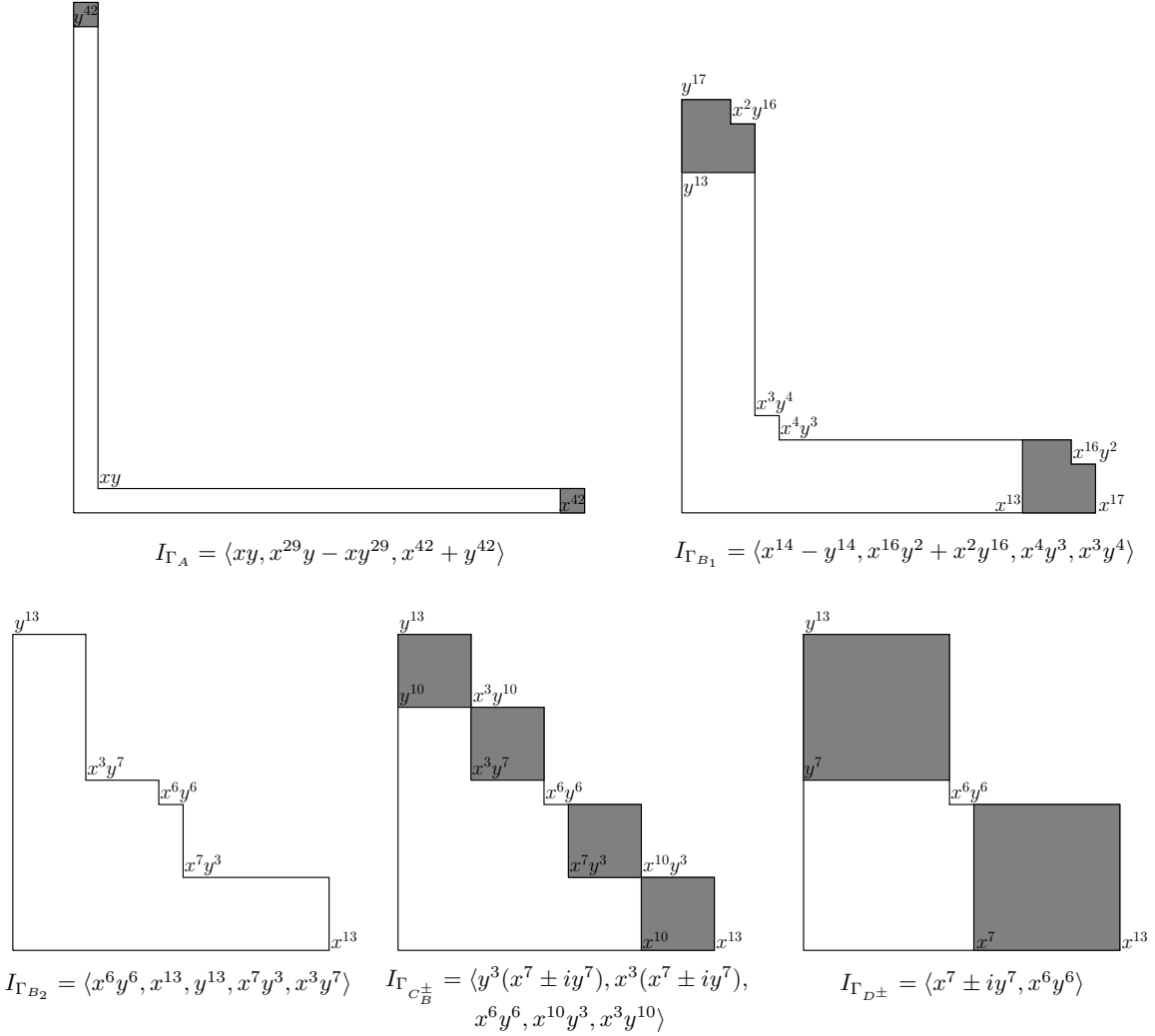
\begin{figure}[h]
\begin{center}
\begin{pspicture}(0,-1.3)(13,12.25)
\scalebox{0.9}{
	\scalebox{0.6}{
	\rput(-0.5,10){
	\psline(0,0)(12.6,0)(12.6,0.6)(0.6,0.6)(0.6,12.6)(0,12.6)(0,0)
	\psline[fillstyle=solid,fillcolor=gray](12,0)(12.6,0)(12.6,0.6)(12,0.6)(12,0)
	\psline[fillstyle=solid,fillcolor=gray](0,12)(0,12.6)(0.6,12.6)(0.6,12)(0,12)
	\rput(0.9,0.8){\Large $xy$}\rput(12.3,0.25){\Large $x^{42}$}\rput(0.3,12.25){\Large $y^{42}$}
	}
	\rput(14.5,10){
	\psline(0,0)(10.2,0)(10.2,1.2)(9.6,1.2)(9.6,1.8)(2.4,1.8)(2.4,2.4)(1.8,2.4)(1.8,9.6)(1.2,9.6)(1.2,10.2)(0,10.2)(0,0)
	\psline[fillstyle=solid,fillcolor=gray](8.4,0)(10.2,0)(10.2,1.2)(9.6,1.2)(9.6,1.8)(8.4,1.8)(8.4,0)
	\psline[fillstyle=solid,fillcolor=gray](0,8.4)(0,10.2)(1.2,10.2)(1.2,9.6)(1.8,9.6)(1.8,8.4)(0,8.4)
	\rput(10.6,0.25){\Large $x^{17}$}\rput(0.3,10.55){\Large $y^{17}$}
	\rput(8.05,0.25){\Large $x^{13}$}\rput(0.35,8.05){\Large $y^{13}$}
	\rput(10.2,1.5){\Large $x^{16}y^2$}\rput(1.8,9.95){\Large $x^2y^{16}$}
	\rput(2.9,2.1){\Large $x^4y^3$}\rput(2.3,2.725){\Large $x^3y^4$}
	}
	\rput(-2,-0.8){
	\psline(0,0)(7.8,0)(7.8,1.8)(4.2,1.8)(4.2,3.6)(3.6,3.6)(3.6,4.2)(1.8,4.2)(1.8,7.8)(0,7.8)(0,0)
	\rput(8.175,0.25){\Large $x^{13}$}\rput(0.35,8.1){\Large $y^{13}$}
	\rput(4.7,2.1){\Large $x^7y^3$}\rput(2.3,4.5){\Large $x^3y^7$}
	\rput(4.1,3.9){\Large $x^6y^6$}
	}
	\rput(7.5,-0.8){
	\psline(0,0)(7.8,0)(7.8,1.8)(6,1.8)(6,3.6)(3.6,3.6)(3.6,6)(1.8,6)(1.8,7.8)(0,7.8)(0,0)
	\psline[fillstyle=solid,fillcolor=gray](6,0)(7.8,0)(7.8,1.8)(6,1.8)(6,0)
	\psline[fillstyle=solid,fillcolor=gray](0,6)(0,7.8)(1.8,7.8)(1.8,6)(0,6)
	\psline[fillstyle=solid,fillcolor=gray](4.2,1.8)(6,1.8)(6,3.6)(4.2,3.6)(4.2,1.8)
	\psline[fillstyle=solid,fillcolor=gray](1.8,4.2)(1.8,6)(3.6,6)(3.6,4.2)(1.8,4.2)
	\rput(8.175,0.25){\Large $x^{13}$}\rput(0.35,8.1){\Large $y^{13}$}
	\rput(4.7,2.1){\Large $x^7y^3$}\rput(2.3,4.5){\Large $x^3y^7$}
	\rput(6.6,2.1){\Large $x^{10}y^3$}\rput(2.4,6.3){\Large $x^3y^{10}$}
	\rput(6.35,0.25){\Large $x^{10}$}\rput(0.35,6.3){\Large $y^{10}$}
	\rput(4.1,3.9){\Large $x^6y^6$}
	}
	\rput(17.5,-0.8){
	\psline(0,0)(7.8,0)(7.8,3.6)(3.6,3.6)(3.6,7.8)(0,7.8)(0,0)
	\psline[fillstyle=solid,fillcolor=gray](4.2,0)(7.8,0)(7.8,3.6)(4.2,3.6)(4.2,0)
	\psline[fillstyle=solid,fillcolor=gray](0,4.2)(0,7.8)(3.6,7.8)(3.6,4.2)(0,4.2)
	\rput(8.175,0.25){\Large $x^{13}$}\rput(0.35,8.1){\Large $y^{13}$}
	\rput(4.5,0.25){\Large $x^7$}\rput(0.25,4.5){\Large $y^7$}
	\rput(4.1,3.9){\Large $x^6y^6$}
	}}
	\rput(3.5,5.4){$I_{\Gamma_{A}}=\langle xy,x^{29}y-xy^{29},x^{42}+y^{42}\rangle$}
	\rput(12,5.4){$I_{\Gamma_{B_{1}}}=\langle x^{14}-y^{14},x^{16}y^2+x^2y^{16},x^{4}y^3,x^3y^{4}\rangle$}
	\rput(1.25,-1){$I_{\Gamma_{B_{2}}}=\langle x^6y^6,x^{13},y^{13},x^7y^3,x^3y^7\rangle$}
	\rput(7,-1){$I_{\Gamma_{C^\pm_{B}}}=\langle y^3(x^7\pm iy^7),x^3(x^7\pm iy^7),$}
	\rput(7.05,-1.5){$x^6y^6,x^{10}y^3,x^3y^{10}\rangle$}
	\rput(12.85,-1){$I_{\Gamma_{D^\pm}}=\langle x^{7}\pm iy^7,x^6y^6\rangle$}
	}
\end{pspicture}
\caption{The $G$-graphs for the group $\BD_{42}(13)$. The shaded areas represent the twin regions.}
\label{BD42}
\end{center}
\end{figure}

\begin{exa} Consider the group $\BD_{42}(13)$. We have that $\frac{42}{13}=[4,2,2,2,4]$ and the lattice points in the Newton polygon that we need to consider are $e_{0}=\frac{1}{42}(0,42)$, $e_{1}=\frac{1}{42}(1,13)$, $e_{2}=\frac{1}{42}(4,10)$ and $e_{3}=\frac{1}{42}(7,7)$. Therefore we have 7 distinguished $\BD_{42}(13)$-graphs shown in Figure \ref{BD42} together with their corresponding ideals. 
\end{exa}

\section{Walking along the exceptional divisor}\label{Sect:Walk}

In this section we prove that every $G$-cluster in $\Hilb{G}{\C^2}$ corresponds to an ideal with $G$-graph of type either $A$, $B$, $C$ or $D$. We start by giving a one-parameter family of ideals which connects any two consecutive $G$-graphs, obtaining every $G$-cluster at the exceptional divisor $E\subset\Hilb{G}{\C^2}$ by passing through all $G$-graphs for a given $\BD_{2n}(a)$ group. 

\begin{thm} \label{walking} Let $G=\BD_{2n}(a)$ and let $\Gamma_{0}, \Gamma_{1}, \ldots, \Gamma_{h-1}, C^+, C^-, D^+, D^-$ be the sequence of $G$-graphs with $I_{\Gamma_{0}}, I_{\Gamma_{1}}, \ldots, I_{\Gamma_{h-1}}, I_{C^+}, I_{C^-}, I_{D^+}, I_{D^-}$ the corresponding defining ideals. Then,
\begin{enumerate}
\item[(i)] For any two consecutive $G$-graphs, $\Gamma_{i}$ and $\Gamma_{i+1}$, there exists a family of ideals $J_{(\xi_{i},\eta_{i})}$ with $(\xi_{i}:\eta_{i})\in\mathbb{P}^1$ such that $J_{(0:1)}=I_{\Gamma_{i}}$, $J_{(1:0)}=I_{\Gamma_{i+1}}$ and $J_{(\xi_{i},\eta_{i})}$ defines a $G$-cluster.
\item[(ii)] There exists a family of ideals $J^+_{(\gamma_{+},\delta_{+})}$ (respectively $J^-_{(\gamma_{-},\delta_{-})}$) with $(\gamma_{+},\delta_{+})\in\mathbb{P}^1$ such that $J^+_{(0:1)}=I_{C^+}$, $J^+_{(1:0)}=I_{D^+}$ and $J^+_{(\gamma_{+},\delta_{+})}$ defines a $G$-cluster (Similarly for $J^-_{(\gamma_{-},\delta_{-})}$). 
\item[(iii)] There exists a family of ideals $J^C_{(\tau,\mu)}$ with $(\tau,\mu)\in\mathbb{P}^1$ such that $J^C_{(0:1)}=I_{C^+}$, $J^C_{(1:0)}=I_{C^-}$, $J^C_{(1:1)}=I_{\Gamma_{h-1}}$ and $J^C_{(\tau,\mu)}$ defines a $G$-cluster.
\end{enumerate}
\end{thm}

\begin{proof}
Let $\Gamma_1:=\Gamma(r,s;u,v)$ and $\Gamma_2:=\Gamma(u,v;t,w)$ be two consecutive $G$-graphs with $I_{\Gamma_1}$ and $I_{\Gamma_2}$ the ideals representing $\Gamma_1$ and $\Gamma_2$ from Propositions \ref{TypeA}, \ref{TypeB1}, \ref{TypeB2}, \ref{TypeD} and \ref{TypeC}. The ideals $I_{\Gamma_1}$ and $I_{\Gamma_2}$ define the $G$-clusters corresponding to the intersection points of two consecutive exceptional curves in $\Hilb{G}{\C^2}$ (see Remark \ref{coeffeq0}). 

Depending on the types of $\Gamma_1$ and $\Gamma_2$ there are 8 possible cases. In every one of them we give a family of ideals $J_{(a:b)}$ parametrized by a $\mathbb{P}^1$ such that $J_{(1:0)}=I_{\Gamma_1}$, $J_{(0:1)}=I_{\Gamma_2}$, and for any $(a:b)\in\mathbb{P}^1$ the ideal $J_{(a:b)}$ defines a $G$-cluster. In other words, we are giving by definition a family of $G$-clusters parametrized by $\mathbb{P}^1_{(a:b)}$ connecting $Z_1$ and $Z_2$. Since $\Hilb{G}{\C^2}$ is a fine moduli space and it is the minimal resolution, the target of the map $\mathbb{P}^1\to\Hilb{G}{\C^2}$ is contained in the exceptional divisor $E\subset\Hilb{G}{\C^2}$. Therefore, by considering all $G$-graphs for a given group $G$ we cover the whole of the exceptional divisor $E$ in this way.

Notice that every generator of $I_{\Gamma_1}$ is contained in a representation of the subset $\{\rho,\rho',\rho_0\} \subset \Irr G$, and every generator of $I_{\Gamma_2}$ is contained in a representation of the subset $\{\rho,\rho'',\rho_0\} \subset \Irr G$, for some $\rho,\rho',\rho''\in\Irr G$ and $\rho_0$ denotes the trivial representation.

Let $f_\rho,g_\rho\in\rho$ be given generators of $I_{\Gamma_1}$ and $I_{\Gamma_2}$ respectively. Notice that $f_\rho\notin I_{\Gamma_2}$ and $g_\rho\notin I_{\Gamma_1}$. The proposed family $J_{(a:b)}$ is given by the union of the two sets of generators with $f_\rho$ and $g_\rho$ replaced by $\sigma_\rho:=af_\rho-bg_\rho$, for $a,b\in\C$. The subset of irreducible representations of $G$ involved in $J_{(a:b)}$ is therefore $\{\rho,\rho',\rho'',\rho_0\}\subset\Irr G$.

We now proof the theorem case by case using the following argument. Let us denote by $I_{\Gamma_1}:=\Span{f_\rho,f_{\rho'},f_{\rho_0}}$ and $I_{\Gamma_2}:=\Span{g_\rho,g_{\rho''},g_{\rho_0}}$ the ideals defining $\Gamma_1$ and $\Gamma_2$. Since we can always choose $1\in\rho_0$ to be in the $G$-graph, notice that $f_{\rho_0}\in I_{\Gamma_2}$ and $g_{\rho_0}\in I_{\Gamma_1}$. Thus in the case $a=0$, to show that $J_{(0:1)}=I_{\Gamma_2}$ it suffices to show that $f_{\rho'}\in I_{\Gamma_2}$. Similarly when $b=0$, to show that $J_{(1:0)}=I_{\Gamma_1}$ it suffices to show that $g_{\rho''}\in I_{\Gamma_1}$. For the rest of the values, if for instance we take $a\neq0$ (or equivalently $a=1$) we show that $\C[x,y]/J_{(1:b)}$ admits $\Gamma_1$ as basis for any value of $b\in\C$. Since $f_\rho\equiv bg_\rho$ mod $J_{(1:b)}$ we can take $g_\rho\in\Gamma_1$ in the basis of $\C[x,y]/J_{(1:b)}$, and it follows that every polynomial in $\Gamma_1$ is not in $J_{(1:b)}$. It remains to show that every polynomial not contained in the representation of $\Gamma_1$ can we written in terms of elements in $\Gamma_1$ modulo $J_{(1:b)}$. The same argument applies if we instead take $b=1$.

\noindent{\em\underline{Case 1}}: $\Gamma_1:=\Gamma_A(r,s;u,v)\to\Gamma_2:=\Gamma_A(u,v;t,w)$: Then the family of ideals $J_{(a:b)}$ given by the polynomials 
\[
\begin{array}{ll}
	\sigma_{\rho}:=ax^{u}y^{u} - b(x^{u+v}+(-1)^{u}y^{u+v}), & \\
	f_{\rho'}:=x^{r+s}+(-1)^ry^{r+s}, &  g_{\rho''}:=x^ty^t, \\
	f_{\rho_0}:=x^{s-v}y^{u-r}+(-1)^{u-r}x^{u-r}y^{s-v}, & g_{\rho_0}:=x^{v-w}y^{t-u}+(-1)^{t-u}x^{t-u}y^{v-w},
\end{array}
\]
which defines a 1-parameter family of $G$-clusters parametrised by a $\mathbb{P}^1$ with coordinates $a$ and $b$. Recall that we have $r<u<t\leq w<v<s$ together with the type $A$ conditions $u<s-v$ and $t<v-w$. 

Let $a=0$. Then by the proof of Lemma \ref{lem:TwinT2} we have that $x^{t+v},y^{t+v}\in I_{\Gamma_2}$, which implies that $x^{r+s},y^{r+s}\in I_{\Gamma_2}$ since $t+v<r+s$. Thus $J_{(0:1)}=I_{\Gamma_2}$. When $b=0$, since $u<t$ we have $g_{\rho''}\in I_{\Gamma_1}$ thus $J_{(1:0)}=I_{\Gamma_1}$.

If $a=1$, we show that $\C[x,y]/J_{(1:b)}$ admits $\Gamma_1$ as basis for any value of $b\in\C$. To show that every polynomial contained in some representation $\rho$ but not in $\Gamma_1$ can we written in terms of elements in $\Gamma_1$ modulo $J_{(1:b)}$, it is enough to check it for monomials $x^jy^u$, $x^uy^j$, $x^ly^{u-r}$, $x^{u-r}y^l$ and $x^{s+u}$, $y^{s+u}$ with $u\leq j< r+s-v$ and $r+s-v\leq l< s+u$ (recall Figure \ref{FigTypeA}).

Indeed, first notice that combining $f_{\rho_0}$ and $g_{\rho''}$ we see that $x^{s-v+t}y^{u-r},x^{u-r}y^{s-v+t}\in I_{(1:b)}$. In addition, using $f_{\rho'}$ we have that $x^{s+u},y^{s+u}\in I_{(1:b)}$. For monomials $x^ly^{u-r}$, by the use of $f_{\rho_0}$ it suffices to show that $x^{u}y^{s-v}$ can be written in terms of elements in $\Gamma_1$ (similarly for $x^{u-r}y^l$). In this case
\[
x^{u}y^{s-v}\equiv-b(x^{u+v}y^{s-u-v}+(-1)^uy^s)\equiv b(bx^{u+2v}y^{s-2u-v}-(-1)^uy^s) \text{ mod }J_{(1:b)}.
\]
By induction, there exists $c>0$ such that $u+cv\geq s-v+t$, thus $x^{u+2v}y^{s-2u-v}$ is either in $J_{(1:b)}$ or in $\Gamma_1$ mod $J_{(1:b)}$.

In the case of monomials $x^jy^u$, by the use of $\sigma_\rho$ it suffices to show that $x^{u}y^{u+v}$ can be written in terms of elements in $\Gamma_1$ (similarly for $x^{u}y^j$). In this case $x^uy^{u+v}\equiv-(-1)^uby^{u+2v}$ which is either in $\Gamma_1$ if $2v<s$, or in $J_{(1:b)}$ otherwise, and we are done.

\noindent{\em\underline{Case 2}}: $\Gamma_1:=\Gamma_A(r,s;u,v)\to\Gamma_2:=\Gamma_{B_1}(u,v;t,w)$: In this case the family of $G$-clusters is defined by the ideals $J_{(a:b)}$ generated by the polynomials 
\[
\begin{array}{lll}
\sigma_{\rho}:=ax^{u}y^{u} - b(x^{u+v}+(-1)^{u}y^{u+v}), & g_{1,\rho''}:=x^ty^m, & f_{\rho_0}:=x^{s-v}y^{u-r}+(-1)^{u-r}x^{u-r}y^{s-v}, \\
f_{\rho'}:=x^{r+s}+(-1)y^{r+s},  & g_{2,\rho''}:=x^my^t, & g_{\rho_0}:=x^{m+v}y^{m-u}+(-1)^{m-u}x^{m-u}y^{m+v}. \\
\end{array}
\]
Notice that $\rho''\in\Irr G$ is 2-dimensional so $I_{\Gamma_2}$ contains two polynomials $g_{1,\rho''},g_{2,\rho''}$ in its set of generators. Recall that in this case $u<s-v$ and $m=v-w=t-u$.
 
As in the previous case, if $a=0$ then $x^{t+v},y^{t+v}\in I_{\Gamma_2}$ and it follows that $f_{\rho'}\in I_{\Gamma_2}$, and when $b=0$ the fact that $u<m$ gives us $g_{1,\rho''},g_{2,\rho''}\in I_{\Gamma_1}$. Then $J_{(0:1)}=I_{\Gamma_2}$ and $J_{(1:0)}=I_{\Gamma_1}$.

The case $a=1$ is identical as Case 1 since $m<t$ so that $g_{1,\rho''},g_{2,\rho''}$ impose this time stronger conditions. 

\noindent{\em\underline{Case 3}}: $\Gamma_1:=\Gamma_{B_1}(r,s;u,v)\to\Gamma_2:=\Gamma_{B_2}(u,v;t,w)$:
In this case the family of ideals $J_{(a:b)}$ is given by:
\[
\begin{array}{lll}
\sigma_{1,\rho}:=ax^{u}y^{m} - by^{s}, & g_{1,\rho''}:=x^ty^m, & f_{\rho_0}:=x^{m+s}y^{m-r}+(-1)^{m-r}x^{m-r}y^{m+s}, \\
\sigma_{2,\rho}:=ax^{m}y^{u} + (-1)^vbx^s, & g_{2,\rho''}:=x^my^t, & g_{\rho_0}:=x^{2m}y^{2m}. \\
f_{\rho'}:=x^{r+s}+(-1)y^{r+s}, &  & \\
\end{array}
\]
In this case $\rho$ and $\rho''$ are 2-dimensional, and for inequalities we have that $m=s-v=u-r=v-w=t-u$ and $u<2m\leq t$. If $a=0$ then clearly $f_{\rho'}\in I_{\Gamma_2}$ and and if $b=0$ then $g_{1,\rho''},g_{2,\rho''}\in I_{\Gamma_1}$, as desired. 

For the rest of the points we consider this time the case $b=1$, and show that we can always take $\Gamma_2$ as basis for $\C[x,y]/J_{(a:1)}$. Since $x^{t}y^m, x^my^{t}\in J_{(a:1)}$, it is enough to check that $x^{s}y^j$, $x^jy^s$ for $0\leq j<m$, can be written in terms of elements in the representation of $\Gamma_2$ (recall Figure \ref{GtypeB2}). But this is obvious since $x^s\equiv(-1)^v ax^my^{u}$ and $y^s\equiv ax^{u}y^{m}$ mod $J_{(a:1)}$, thus $x^sy^{j}\equiv(-1)^v ax^my^{u+j}$ mod $J_{(a:1)}$ for $1\leq j<m$, where $x^my^{u+j}$ is contained in the representation of $\Gamma_2$ (similarly for $x^sy^{j}$ for $1\leq j<m$), and we are done.

\noindent{\em\underline{Case 4}}: $\Gamma_1:=\Gamma_{B_2}(r,s;u,v)\to\Gamma_2:=\Gamma_{B_2}(u,v;t,w)$: The family of ideals $J_{(a:b)}$ is given by the generators 
\[
\begin{array}{llll}
\sigma_{1,\rho}:=ax^{u}y^{m} - by^{s}, & f_{1,\rho'}:=y^{s+m}, & g_{1,\rho''}:=x^ty^m, & f_{\rho_0}:=x^{2m}y^{2m}, \\
\sigma_{2,\rho}:=ax^{m}y^{u} + (-1)^vbx^s, & f_{2,\rho'}:=x^{s+m}, & g_{2,\rho''}:=x^my^t, 
\end{array}
\]
In this case, $\rho,\rho'$ and $\rho''$ are 2-dimensional, $f_{\rho_0}=g_{\rho_0}$ and $m$ is as in the previous case but $u,t\geq 2m$. It is immediate to see in this case that if $a=0$ then $f_{1,\rho'},f_{2,\rho'}\in I_{\Gamma_2}$, and if $b=0$ then $g_{1,\rho''},g_{2,\rho''}\in I_{\Gamma_1}$, as desired. The case $b=1$ is analogous to Case 3.  \\

From now on suppose that the last $G$-graph is $\Gamma_{h-1}:=\Gamma(r,s;q,q)$. Depending on the type of $\Gamma_{h-1}$ we have a different $G$-graph of type $C$, which we will denote by a subindex, $A$ or $B$.

\noindent{\em\underline{Case 5}}: $\Gamma_1:=\Gamma_{C_{A}^-}(r,s;q,q)\to\Gamma_2:=\Gamma_{C_{A}^+}(r,s;q,q)$: In this case the family of ideals $J_{(a:b)}$ is given by the generators:
\[
\begin{array}{ll}
\sigma_\rho:=a(x^q-(-i)^qy^q)^2 - b(x^q+(-i)^qy^q)^2, & g_{\rho''}:=x^sy^{m_{2}} + (-1)^ri^qx^{m_{2}}y^{s} \\
f_{\rho'}:=x^sy^{m_{2}} - (-1)^ri^qx^{m_{2}}y^{s}, & f_{\rho_0}:=x^{m_{1}}y^{m_{2}}+(-1)^{m_{2}}x^{m_{2}}y^{m_{1}}, 
\end{array}
\]
Notice that $f_{\rho_0}=g_{\rho_0}$ and recall that $r<q<s$ and $b_hq=r+s$.

Suppose that $a=0$ so that $\sigma_\rho=x^{2q}+2(-i)^qx^qy^q+(-1)^qy^{2q}$. We show that $f_{\rho'}\in I_{\Gamma_2}$. Define the polynomials $F_j:=x^{s-jq}y^{(j+1)q-r}+cx^{(j+1)q-r}y^{s-jq}$ for $j=0,\ldots,b_h-1$ and $c:=-(-1)^{r+jq}i^q$. Note that $f_{\rho'}=F_0$ (recall Figure \ref{GtypeC} (a)). Notice that $F_j+2(-i)^qF_{j+1}+(-1)^qF_{j+2}\in\Span{\sigma_\rho}$ and that $F_j$ and $F_{b_h-1-j}$ are the same up to a constant. Then $f_{\rho'}\in I_{\Gamma_2}$ if and only if $F_{\frac{b_h}{2}-2}, F_{\frac{b_h}{2}-1}\in I_{\Gamma_2}$, which it is true since $F_{\frac{b_h}{2}-2}+(2(-i)^q+c)F_{\frac{b_h}{2}-1},(1+2(-i)^q)cF_{\frac{b_h}{2}-1}+(-1)^qF_{\frac{b_h}{2}-2}\in I_{\Gamma_2}$. The case when $b=0$ is done identically up to a change of sign.

Now suppose that $a=b=1$. Then $\sigma_\rho=x^qy^q$ and $f_{\rho'},g_{\rho''}\in\Span{\sigma_\rho,f_{\rho_0}}$, which imply that $J_{(1:1)}$ represents the $G$-graph $\Gamma_{h-1}$ of type $A$. It follows that using the same argument as in Case 1, fixing the value $a=1$ we obtain that $\Gamma_{h-1}$ can be chosen to be the basis for $\C[x,y]/J_{(1:b)}$ for every $b\in\C$, and we are done.

\noindent{\em\underline{Case 6}}: $\Gamma_1:=\Gamma_{C_{B}^+}(r,s;q,q)\to\Gamma_2:=\Gamma_{C_{B}^-}(r,s;q,q)$: If the last $G$-graph is of type $B_{1}$ then the family is given by the ideal $J_{(a:b)}$ generated by 
\[
\begin{array}{lll}
\sigma_{1,\rho}:=k_{-}(x^s+(-i)^qx^my^{q}) - k_{+}(x^s-(-i)^qx^my^{q}), & f_{1,\rho'}:=x^{s}y^{m}, & f_{\rho_0}:=x^{2m}y^{2m},  \\
\sigma_{2,\rho}:=k_{-}(y^s+i^qx^qy^{m}) - k_{+}(y^s+i^qx^qy^{m}), & f_{2,\rho'}:=x^{m}y^s. & \\
\end{array}
\] 
In this case $\rho,\rho'$ and $\rho''$ are 2-dimensional. By the equalities $f_{1,\rho'}=g_{1,\rho''}$, $f_{2,\rho'}=g_{2,\rho''}$ and $f_{\rho_0}=g_{\rho_0}$ it is clear that $J_{(0:1)}=I_{\Gamma_2}$ and $J_{(1:0)}=I_{\Gamma_1}$.

As in the previous case, the ideal $J_{(1:1)}$ represents the $G$-graph $\Gamma_{h-1}$ of type $B$, and using the same argument as in Case 3, it follows that if $b=1$ we can always find a basis for $\C[x,y]/J_{(a:1)}$ which is a $G$-graph of type $B$.

The last two cases correspond to the families at the two ``horns'' of the exceptional divisor.

\noindent{\em\underline{Case 7}}: $\Gamma_1:=\Gamma_{C_{A}^\pm}(r,s;q,q)\to\Gamma_2:=\Gamma_{D^\pm}(r,s;q,q)$: In this case the family of ideals $J_{(a_\pm:b_\pm)}$ is given by 
\[
\begin{array}{ll}
\sigma_{\rho}:= a_{\pm}(x^sy^{m_{2}}\pm(-1)^ri^qx^{m_{2}}y^s) - b_{\pm}(x^q\pm(-i)^qy^q), & g_{\rho_0}:= x^{s-r}y^{s-r}, \\
f_{\rho'}:= (x^q\pm(-i)^qy^q)^2, & f_{\rho_0}:= x^{m_{1}}y^{m_{2}}+(-1)^{m_{2}}x^{m_{2}}y^{m_{1}}.
\end{array}
\]

It is straightforward to see that $J_{(0:1)}=I_{\Gamma_2}$ and $J_{(1:0)}=I_{\Gamma_1}$ corresponding to the cases $a=0$ and $b=0$ respectively.

Now suppose that $a_+=1$. The case $a_-=1$ is done similarly. First notice that using $f_{\rho_0}$ we have that $x^{2q-2r}y^{2s-2q},x^{2s-2q}y^{2q-2r}\in J_{(1:b_+)}$. Then $x^{s-r+q}y^{s-r-q}\equiv\frac{w_+}{z_+}x^{2s-r}y^{s-2r}\equiv0$ mod $J_{(1:b_+)}$, where the last equivalence is obtained because $(b_h-1)q<s$ and $b_h\neq2$ since the last $G$-graph is of type $A$. By the same process we see that for $0\leq j\leq k+1$ we have that $x^{s-r+jq}y^{s-r-jq},x^{s-r-jq}y^{s-r+jq}\in J_{(1:b_+)}$, thus every monomial not contained in the representation of the $G$-graph of type $D^+$ belongs to the ideal $J_{(1:b_+)}$, and we are done.

\noindent{\em\underline{Case 8}}: $\Gamma_{C_{B}^\pm}\to\Gamma_{D^\pm}$: In this case the family of ideals $J_{(a_\pm:b_\pm)}$ is given by 
\[
\begin{array}{ll}
\sigma_\rho:= a_{\pm}(x^sy^{m}\pm(-1)^ri^qx^{m}y^s) - b_{\pm}(x^q\pm(-i)^qy^q), & f_{1,\rho'}:=x^m(x^q\pm(-i)^qy^q), \\
f_{\rho_0}:= x^{s-r}y^{s-r}, & f_{2,\rho'}:=y^m(x^q\pm(-i)^qy^q),\\
\end{array}
\]
where $m=s-q=q-r$. In this case $\rho'$ is 2-dimensional and $f_{\rho_0}=g_{\rho_0}$. 

This case is analogous as the previous one. Again it is straightforward to see that $J_{(0:1)}=I_{\Gamma_2}$ and $J_{(1:0)}=I_{\Gamma_1}$ corresponding to the cases $a=0$ and $b=0$ respectively. Finally suppose that $a_+=1$ (the case $a_-=1$ is done similarly). Then $J_{(1,b_+)}=\Span{\sigma_\rho,f_{\rho_0}}$, and in the same way as in the previous case we have that $x^{s-r+jq}y^{s-r-jq},x^{s-r-jq}y^{s-r+jq}\in J_{(1:b_+)}$ for $0\leq j\leq k+1$, and we are done.
\end{proof}

\begin{thm} \label{ABCD} Let $G=\BD_{2n}(a)$ be small and let $P\in G$-$\HILB(\C^2)$ be defined by the ideal $I$. Then we can always choose a basis for $\C[x,y]/I$ from one of the following list:
\[ \Gamma_{A}, \Gamma_{B}, \Gamma_{C^+}, \Gamma_{C^-}, \Gamma_{D^+}, \Gamma_{D^-} \]
\end{thm}

\begin{proof} By construction, every point in $\Hilb{G}{\C^2}$ away from the ``horns'' corresponds to a pair of $H$-clusters in $\Hilb{H}{\C^2}$. Therefore, we can chose for these $H$-clusters the $H$-graphs $\widetilde{\Gamma}(r,s;u,v)$ and $\widetilde{\Gamma}(v,u;s,r)$ identified by $\beta$ for some $(r,s)$ and $(u,v)$ boundary lattice points in the Newton polygon for $\frac{1}{2n}(1,a)$, and we can take $\Gamma(r,s;u,v)$ to be the $G$-graph (of type $A$ or $B$) for our $G$-cluster. 

For the clusters in the exceptional ``horns'' $E^+$ and $E^-$, we know by Theorem \ref{walking} that these exceptional curves are covered by the ideals $J^+_{(\gamma_{+},\delta_{+})}$ and $J^-_{(\gamma_{-},\delta_{-})}$, which correspond to $G$-graphs of type $C^\pm$ and $D^\pm$, and we are done. 
\end{proof}

Let $U_{\Gamma}$ the open set in $\Hilb{G}{\C^2}$ which consists of all $G$-clusters $\mathcal{Z}$ such that $\Gamma$ is a basis of $\mathcal{O}_{\mathcal{Z}}$. As a corollary of the previous theorem we have that the $G$-graphs for a $\BD_{2n}(a)$ group gives us the open set for the covering of $\Hilb{G}{\C^2}$ that we are looking for.
 
\begin{col}\label{cover}
Let $G=\BD_{2n}(a)$ a small binary dihedral group and let $\Gamma_{0}, \Gamma_{1}, \ldots, \Gamma_{h-1},  \Gamma_{C^+}, \Gamma_{C^-}$, $\Gamma_{D^+}, \Gamma_{D^-}$ the list of $G$-graphs. Then 
\[
U_{\Gamma_{0}}, U_{\Gamma_{1}}, \ldots, U_{\Gamma_{h-1}},  U_{\Gamma_{C^+}}, U_{\Gamma_{C^-}}, U_{\Gamma_{D^+}}, U_{\Gamma_{D^-}}
\]
form an open cover of $G$-$\HILB(\C^2)$.
\end{col}

\begin{rem} By deforming the $G$-graph $\Gamma_{i}$ located at the origin, we can calculate the explicit equation of the open set $U_{\Gamma_{i}}$ as it is done in \cite{Leng} for binary dihedral subgroups in $\SL(2,\C)$. Extending the groups to $\GL(2,\C)$ increases the amount of choices for generators of these ideals, which makes this approach much harder in practice. Fortunately, one can associate to any $G$-graph an open set of the moduli space $\mathcal{M}_{\theta}(Q,R)$ of $\theta$-stable quiver representations of the McKay quiver with relations $(Q,R)$, which coincides with $\Hilb{G}{\C^2}$, and where the calculation of the open set turns out to be much easier. See \cite{NdC2}. 
\end{rem}

\section{Special representations} \label{SpecialRepr}

For a finite small subgroup $G\subset\GL(2,\C)$, the special McKay correspondence states that there is a one-to-one correspondence between exceptional divisors $E_{i}$ in the minimal resolution of $\C^2/G$ and the {\em special} irreducible representations $\rho_{i}$ of $G$. In \cite{Ish02} Ishii proves that the minimal resolution is in fact $\Hilb{G}{\C^2}$. 

\begin{thm}[\cite{Ish02}, \S7.1]\label{thm!Ishii} Let $G\subset\GL(2,\C)$ be small and denote by $I_{y}$ the ideal corresponding to $y\in G$-$\HILB(\C^2)$ and by $\mathfrak{m}$ the maximal ideal of $\mathcal{O}_{\C^2}$ corresponding to the origin 0. If $y$ is in the exceptional locus, then we have an isomorphism
\[
I_{y}/\mathfrak{m}I_{y}\cong\left\{
	\begin{array}{ll}
	\rho_{i}\oplus\rho_{0} & \text{if $y\in E_{i}$, and $y\notin E_{j}$ for $j\neq i$}, \\
	\rho_{i}\oplus\rho_{j}\oplus\rho_{0} & \text{if $y\in E_{i}\cap E_{j}$},
	\end{array}\right.
\]
as representations of $G$, where $\rho_{i}$ is the special representation associated with the irreducible  exceptional curve $E_{i}$.
\end{thm}

In other words, for any point in the exceptional divisor of $\Hilb{G}{\C^2}$, only the trivial and the special representations corresponding to the curves which the point lies on are involved in the ideal defining the $G$-cluster. In our case, the explicit description of these ideals is the following:

\begin{prop} \label{0strata}
Let $G=\BD_{2n}(a)$ be small and $y\in\Hilb{G}{\C^2}$ be a point in the exceptional locus. Denote by $I_{y}$ the ideal defining $y$ and by $\Gamma_{y}=\Gamma(r,s;u,v)$ the corresponding $G$-graph.  Then 
{\renewcommand{\arraystretch}{1.25}
\[
I_{y}/\mathfrak{m}I_{y}\cong\left\{
	\begin{array}{ll}
	\rho^{(-1)^r}_{r+s}\oplus\rho^{(-1)^{u}}_{u+v}\oplus\rho^+_{0} & \text{if $\Gamma_{y}$ is of type $A$,} \\
	\rho^{(-1)^r}_{r+s}\oplus V_{r}\oplus\rho^+_{0} & \text{if $\Gamma_{y}$ is of type $B_{1}$,} \\
	V_{2r-u}\oplus V_{r}\oplus\rho^+_{0} & \text{if $\Gamma_{y}$ is of type $B_{2}$,} \\
	\rho^{(-1)^q}_{2q}\oplus\rho^{\pm}_{q}\oplus\rho^+_{0} & \text{if $\Gamma_{y}$ is of type $C^{\pm}$ and $\Gamma_{h-1}$ is of type $A$,} \\
	V_{r}\oplus\rho^{\pm}_{q}\oplus\rho^+_{0} & \text{if $\Gamma_{y}$ is of type $C^{\pm}$ and $\Gamma_{h-1}$ is of type $B$,} \\
	\rho^{\pm}_{q}\oplus\rho^+_{0} & \text{if $\Gamma_{y}$ is of type $D^{\pm}$,} 
	\end{array}\right.
\]}
where $\mathfrak{m}$ the maximal ideal of $\mathcal{O}_{\C^2}$ corresponding to the origin 0.
\end{prop} 

\begin{proof} Reformulating Propositions \ref{TypeA}, \ref{TypeB1}, \ref{TypeB2}, \ref{TypeD} and \ref{TypeC} in the language of Theorem \ref{thm!Ishii}, we see that the representations involved in the generators of each of the ideals are the ones presented above. By Theorem \ref{walking} every point in the exceptional divisor $E\subset\Hilb{G}{\C^2}$ is defined by one of those ideals, so the result follows.
\end{proof}

As a consequence of both \ref{thm!Ishii} and \ref{0strata} we obtain the special representations for any group $\BD_{2n}(a)$ in terms of the continued fraction $\frac{2n}{a}$, which we list in the following theorem.

\begin{thm} \label{thm:Special}
Let $\Gamma_{0}, \Gamma_{1}, \ldots, \Gamma_{h-1}$ the sequence of $qG$-graphs given by $e_{0}=\frac{1}{2n}(0,2n), e_{1}=\frac{1}{2n}(1,a), e_{2}=\frac{1}{2n}(c_{1},d_{1}), \ldots, e_{m-1}=\frac{1}{2n}(c_{m-2},d_{m-2}), e_{m}=\frac{1}{2n}(q,q)$, where 
\begin{align*}
\Gamma_{0}, \ldots, \Gamma_{i} & \text{ are of type $A$} \\
\Gamma_{i+1} & \text{ is of type $B_{1}$} \\
\Gamma_{i+2}, \ldots, \Gamma_{h-1} & \text{ are of type $B_{2}$} 
\end{align*}
Then, the special representations are
\begin{align*}
\rho^{-}_{1+a}, \rho^{(-1)^{c_{1}}}_{c_{1}+d_{1}}, \rho^{(-1)^{c_{2}}}_{c_{2}+d_{2}}\ldots, \rho^{(-1)^{c_{i+1}}}_{c_{i+1}+d_{i+1}} & \text{ from type $A$,} \\
V_{c_{i}} & \text{ from type $B_{1}$,} \\
V_{c_{i+1}}, \ldots, V_{c_{h-2}} & \text{ form type $B_{2}$ and} \\
\rho^+_{q}, \rho^-_{q} & \text{ from types $C$ and $D$} 
\end{align*}
\end{thm}

\begin{rem} We want to note that the same result holds for groups of the form $\BD_{2n}(a,q)$ since the $G$-graphs are not affected by the change in the generator $\beta$. In other words, Theorem \ref{thm:Special} is valid also for any small binary dihedral subgroup $G\subset\GL(2,\C)$ with maximal cyclic subgroup $H=\Span{\frac{1}{2n}(1,a)}$.
\end{rem}

There is also a relation between $G$-graphs of types $A$ and $B$, and the dimension of the corresponding irreducible special representations. Let $E=\bigcup E_{i}\subset \BD_{2n}(a)$-Hilb$(\mathbb{C}^2)$ be the exceptional locus. By \cite{Wun88}, the dimension of the special representations $\rho_{i}$ corresponding to $E_{i}$ is equal to the coefficient of $E_{i}$ in the fundamental cycle $Z_{\text{fund}}$ (the smallest effective divisor such that $Z_{\text{fund}}\cdot E_{i}\leq0$). Let 
$-2,-2,-a_{m},\ldots,-a_{2},-a_{1}$ be the selfintersections along the minimal resolution of $G/\BD_{2n}(a)$, where the first two $-2s$ correspond to the ``horns'' of the Dynkin diagram.  

\begin{col} \label{ColSpecial}
The special irreducible representations of $\BD_{2n}(a)$ are all 1-dimensional if and only if the middle entry $b_{h}$ in the continued fraction $\frac{2n}{a}=[b_{1},\ldots,b_{h},\ldots,b_{1}]$ is not 2. 
\end{col}

\begin{proof} The result follows from Corollary \ref{NoTypeB} (ii) and Theorem \ref{thm:Special}. 
\end{proof}

Thus we have a one-to-one correspondence between $qG$-graphs of type $A$ and 1-dimensional special representations (except the two corresponding to the ``horns'' which are also 1-dimensional and they are covered by the $G$-graphs of type $C$ and $D$), and another correspondence between $qG$-graphs of type $B$ and 2-dimensional special representations.

\bibliographystyle{plain}

\begin{thebibliography}{10}

\bibitem{BKR01}
Tom Bridgeland, Alastair King, and Miles Reid.
\newblock The {M}c{K}ay correspondence as an equivalence of derived categories.
\newblock {\em J. Amer. Math. Soc.}, 14(3):535--554 (electronic), 2001.

\bibitem{Bri}
Egbert Brieskorn.
\newblock Rationale {S}ingularit\"aten komplexer {F}l\"achen.
\newblock {\em Invent. Math.}, 4:336--358, 1967/1968.

\bibitem{CR02}
Alastair Craw and Miles Reid.
\newblock How to calculate {$A$}-{H}ilb {$\mathbb C\sp 3$}.
\newblock In {\em Geometry of toric varieties}, volume~6 of {\em S\'emin.
  Congr.}, pages 129--154. Soc. Math. France, Paris, 2002.

\bibitem{Ish02}
Akira Ishii.
\newblock On the {M}c{K}ay correspondence for a finite small subgroup of {${\rm
  GL}(2,\Bbb C)$}.
\newblock {\em J. Reine Angew. Math.}, 549:221--233, 2002.

\bibitem{Ito02}
Yukari Ito.
\newblock Special {M}c{K}ay correspondence.
\newblock In {\em Geometry of toric varieties}, volume~6 of {\em S\'emin.
  Congr.}, pages 213--225. Soc. Math. France, Paris, 2002.

\bibitem{IN96}
Yukari Ito and Iku Nakamura.
\newblock Mc{K}ay correspondence and {H}ilbert schemes.
\newblock {\em Proc. Japan Acad. Ser. A Math. Sci.}, 72(7):135--138, 1996.

\bibitem{IW08}
Osamu Iyama and Michael Wemyss.
\newblock The classification of special {CM} modules.
\newblock {\em Math. Z.}, 265(1):41Ð83, 2010.

\bibitem{Kidoh}
Rie Kidoh.
\newblock Hilbert schemes and cyclic quotient surface singularities.
\newblock {\em Hokkaido Math. J.}, 30(1):91--103, 2001.

\bibitem{Leng}
Becky Leng.
\newblock The Mckay correspondence and orbifold Riemann-Roch.
\newblock {\em PhD thesis, University of Warwick}, 2002.

\bibitem{McK80}
John McKay.
\newblock Graphs, singularities, and finite groups.
\newblock In {\em The {S}anta {C}ruz {C}onference on {F}inite {G}roups ({U}niv.
  {C}alifornia, {S}anta {C}ruz, {C}alif., 1979)}, volume~37 of {\em Proc.
  Sympos. Pure Math.}, pages 183--186. Amer. Math. Soc., Providence, R.I.,
  1980.

\bibitem{Nak01}
Iku Nakamura.
\newblock Hilbert schemes of abelian group orbits.
\newblock {\em J. Algebraic Geom.}, 10(4):757--779, 2001.

\bibitem{thesis}
Alvaro Nolla~de Celis.
\newblock Dihedral groups and ${G}$-{H}ilbert schemes.
\newblock {\em PhD thesis, University of Warwick}, 2008.

\bibitem{NdC2}
Alvaro Nolla~de Celis.
\newblock Dihedral ${G}$-{H}ilb via representations of the {M}ckay quiver.
\newblock {\em Proc. Japan Acad. Ser. A}, 88(5):78--83, 2012.

\bibitem{Reid}
Miles Reid.
\newblock Surface cyclic quotient singularities and {H}irzebruch--{J}ung
  resolutions.
\newblock Available at \href{http://homepages.warwick.ac.uk/~masda/surf/more/cyclic.pdf}{http://homepages.warwick.ac.uk/$\sim$masda/surf/more/cyclic.pdf}, accessed 20 June 2012.

\bibitem{Rei02}
Miles Reid.
\newblock La correspondance de {M}c{K}ay.
\newblock {\em Ast\'erisque}, 276:53--72, 2002.
\newblock S{\'e}minaire Bourbaki, Vol. 1999/2000.

\bibitem{Rie74}
Oswald Riemenschneider.
\newblock Deformationen von {Q}uotientensingularit\"aten (nach zyklischen
  {G}ruppen).
\newblock {\em Math. Ann.}, 209:211--248, 1974.

\bibitem{Rie04}
Oswald Riemenschneider.
\newblock Special representations and the two-dimensional {M}c{K}ay
  correspondence.
\newblock {\em Hokkaido Math. J.}, 32(2):317--333, 2003.

\bibitem{Seb05}
Magda Sebestean.
\newblock Correspondance de McKay et \'equivalences d\'eriv\'es.
\newblock {\em PhD thesis, Universit\'e Paris 7}, 2005.

\bibitem{Ch}
G.~C. Shephard and J.~A. Todd.
\newblock Finite unitary reflection groups.
\newblock {\em Canadian J. Math.}, 6:274--304, 1954.

\bibitem{Wun88}
J{\"u}rgen Wunram.
\newblock Reflexive modules on quotient surface singularities.
\newblock {\em Math. Ann.}, 279(4):583--598, 1988.

\bibitem{Yos}
Yuji Yoshino.
\newblock Cohen-{M}acaulay modules over {C}ohen-{M}acaulay rings, volume
  146 of {\em London Mathematical Society Lecture Note Series}.
\newblock Cambridge University Press, Cambridge, 1990.

\end{thebibliography}

Alvaro Nolla de Celis \\
Graduate School of Mathematics, Nagoya University, Chikusa-ku, Nagoya 464-8602, Japan. \\
{\em email address}: alnolla@gmail.com

\end{document}